\documentclass[12pt,oneside]{article}

\RequirePackage[OT1]{fontenc}
\RequirePackage{amsthm,amsmath}
\RequirePackage{natbib}
\usepackage{authblk}
\usepackage{graphicx}
\usepackage{subcaption}
\numberwithin{equation}{section}
\theoremstyle{plain}

\usepackage{float}

\usepackage{amssymb,amsmath,color}

\usepackage{epsfig}
\usepackage{dsfont}
\usepackage{color}
\usepackage{enumerate}
\usepackage[official]{eurosym}

\usepackage{etoolbox}
\makeatletter
\patchcmd{\@makecaption}
  {\parbox}
  {\advance\@tempdima-\fontdimen2} % decrease the width!
  {}{}
\makeatother

\newtheorem{theorem}{Theorem}

\newtheorem{proposition}[theorem]{Proposition}

\newcommand{\norm}[1]{\left\lVert#1\right\rVert}

\usepackage{setspace,caption}
\usepackage{lipsum}
\title{More than one Author with different Affiliations}
\author[1]{Munir Hiabu	\thanks{munir.hiabu@sydney.edu.au (Corresponding author) \\ Acknowledgement: We are grateful to the authors of \cite{Lin:etal:16}, who provided us with code for the calculation of their estimator.}}
\author[2]{Enno Mammen}%\thanks{mammen@math.uni-heidelberg.de}}
\author[3]{M. Dolores Mart{\'\i}nez-Miranda}%\thanks{mmiranda@ugr.es }}
\author[4]{Jens P. Nielsen}%\thanks{Jens.Nielsen.1@city.ac.uk}}
\affil[1]{School of Mathematics and Statistics, University of Sydney, Camperdown NSW 2006, Australia}
\affil[2]{Institute for Applied Mathematics, Heidelberg University, Im Neuenheimer Feld 205, 69120 Heidelberg, Germany}
\affil[3]{Department of Statistics and Operations Research, University of Granada, Campus Fuentenueva s/n. 18071 Granada, Spain  }
\affil[4]{Cass Business School, City, University of London, 106 Bunhill Row, London, EC1Y 8TZ, United Kingdom}

\title{Smooth backfitting of proportional hazards with multiplicative components}

\begin{document}

  \maketitle

\newpage

\begin{abstract}
Smooth backfitting has proven to have a number of theoretical and practical advantages in structured regression. By projecting the data down onto the structured space of interest smooth backfitting provides  a direct link between data and estimator. This paper introduces the ideas of smooth backfitting to survival analysis in a proportional hazard model, where we assume an underlying conditional hazard with multiplicative components.  We develop asymptotic theory for the estimator.
In a comprehensive simulation study we show that our smooth backfitting estimator successfully circumvents the curse of dimensionality and outperforms existing estimators. This is especially the case in difficult situations like high number of covariates and/or high correlation between the covariates, where other estimators tend to break down. 
We use the smooth backfitter  in a practical application where we extend recent advances of in-sample forecasting methodology by allowing more information to be incorporated, while still obeying the structured requirements of in-sample forecasting.
\end{abstract}

%\begin{keyword}[class=MSC]
%\kwd{	62G20}
%\kwd{62G05}
%\kwd{62P05}
%\end{keyword}
%
%\begin{keyword}
%\kwd{Aalen's multiplicative model}
%\kwd{ Local linear kernel estimation} 
%\kwd{Survival analysis}
%\kwd{ In-sample forecasting}
%\kwd{Smooth backfiitng}
%\end{keyword}
%
%\end{frontmatter}
%

%
%\begin{keywords}
%%Aalen's multiplicative model;  Local linear kernel estimation; Survival data.
%\end{keywords}
\section{Introduction}
Nonparametric models suffer from the curse of dimensionality in high dimensional data spaces. 
Random forests \citep{Breiman:01} circumvent the dimensionality problem by assuming that not all variables are relevant and that the function of interest can be approximated well by piecewise constant functions; see 
\cite{Wright:Ziegler:17} for a recent survival implementation.
An alternative is to introduce  some structure that  stabilizes  the system. Introducing structure  has the additional advantage that it allows to visualize, interpret, extrapolate and forecast the properties of the underlying data.  
In this paper, we concentrate on structured models.
The smooth backfitting algorithm of \citet{Mammen:etal:99} consideres the simplest nonparametric structure in the regression context - the additive structure. It has
many theoretical and practical advantages to earlier approaches of regression backfitting.
 The popular regression backfitting approach of \citet{Hastie:Tibshirani:90} are  numerical  iterating-procedures estimating  one component given the estimates of the rest. In contrast,  smooth backfitting is a direct projection of the data down onto the structured space of interest. This direct relationship between data and estimates gives a more solid grip on what is being estimated and the theoretical properties underlying it, see also \citet{Nielsen:Sperlich:05} and \citet{Huang:Yu:19}.   The purpose of this paper is to introduce smooth backfitting to the  field of survival analysis and nonparametric smooth hazard estimation. While the additive structure is the most natural and most widely used in regression, the multiplicative structure seems more natural in hazard estimation. The omnipresent Cox regression model is a proportional hazard model (also known as Cox proportional hazard model) and many extensions and alternatives to the Cox regression model have been formulated in a multiplicative framework. We will  therefore consider the  multiplicative structure in this paper.  It  can  be used in applications to test Cox regression or other proportional hazard models either visually or quantitatively. But this is beyond the scope of this paper. 
Multiplicative smooth backfitting is theoretically more challenging  than additive smooth backfitting.  The smooth backfitting multiplicative regression structure  was analysed by  \citet{Yu:etal:08} as a special case of generalized additive models. \citet{Yu:etal:08}  showed that the multiplicative structure - in contrast to the simpler additive case - provides asymptotic theory with a number of non-trivial interactions between exposure available in different directions.  Naturally, the asymptotics provided here for smooth backfitting of multiplicative hazards contains similar interactive components in the asymptotic theory. 
The survival projection introduced in this paper is different and less intuitive than the nonparametric regression considered in \citet{Mammen:etal:99, Mammen:etal:01}. As in density estimation, see \citet{Jones:93}, hazard estimation requires a projection of a dirac-delta-sequence related to the jumps of the counting process, see also \citet{Nielsen:98} and \citet{Nielsen:Tanggaard:01}.
 We   provide a simple algorithm first projecting the data down onto an unstructured estimator, and then further projecting the unstructured estimator down onto the multiplicative space of interest. Our numerical algorithms are greatly simplified by a new principle of weighting the projection according to the final estimates.% In this paper, the focus is on the local constant estimator \citep{Nielsen:Linton:95}. The suggested approach is  to project the local constant estimator   down to the space of multiplicative hazard functions.
 
We  consider a  multiplicatively structured proportional hazard model:
 \begin{align}\label{multhazard1}
\alpha(t, z)=\alpha_0(t)\alpha_1(z_1)\cdots\alpha_d(z_d), %\quad (t>0, z_1,\dots, z_d \in \mathbb R)
\end{align}
where $\alpha_k,$ $k=0, \dots, d$, are some smooth and positive one-dimensional functions and $z=(z_1, \dots, z_d)\in \mathbb R^d$ are possibly time $t$-dependent covariates. We do not impose  any further structural assumption on $\alpha_k$.
This is in contrast to the semiparametric Cox proportional hazard model where all  components, with exception of the baseline hazard $\alpha_0$, are assumed to take a   log-linear shape. We do expect that smooth backfitting of proportional hazard models can be generalized to much the same way as the  Cox proportional hazard model is generalized. As one recent example that could be interesting to treat as a smooth backfitting problem, see \cite{Hsu:etal:18}.

Estimators for model \eqref{multhazard1} can be  categorized in four groups:
(i) \citet{Therneau:etal:90} and \citet{Grambsch:etal:95} start with  the Cox model and investigate smoothed residual plots;
  (ii) \citet{Hastie:Tibshirani:90,Hastie:Tibshirani:90b}, \citet{OSullivan:88,OSullivan:93},   \citet{Sleeper:Harrington:90} and \citet{Huang:99}  consider 
  splines via penalized partial likelihood;
 (iii) \citet{Linton:etal:03}, \citet{Honda:05} build on marginal integration \citep{Linton:Nielsen:95} and 
 (iv) \citet{Lin:etal:16} use kernel smoothers starting from a  global partial likelihood criterion.
 
\citet{Lin:etal:16} prove asymptotic efficiency of their estimator and they show by a detailed simulation study that their estimator outperforms   the proposals in \citet{Huang:99}, \citet{Linton:etal:03}, \citet{Honda:05}. For this reason in this paper we take their estimator as benchmark. 
It has been argued that in additive regression models smooth backfitting is less affected by sparseness of high-dimensional data and by strong correlated covariables, see  \citet{Nielsen:Sperlich:05}. In this paper we will show that this also holds for our smooth backfitter in a  multiplicatively structured proportional hazard model. For this purpose we study the smooth backfitter in settings which include challenging high dimensional data and highly correlated covariates. The smooth backfitting approach turns out to show a very good performance compared to all other estimators and to be very robust to both, high dimensional and correlated data. In particular in high-dimensional settings it strongly outperforms the approach of \citet{Lin:etal:16} that already  for low dimensions runs into instabilities and numerical problems. 
In this paper, we do not treat the problem of variable selection. It would in particular be interesting to investigate  Lasso approaches or  nonparametric tests on the significance of one component.
For the development of tests our theory may be used as a first step. 

The paper is structured as follows.
 Section 2 contains the mathematics of the underlying survival model. In Section 3  the smooth backfitting estimator is defined as a projection of unstructured hazard estimators. This is done for unstructured hazard estimators that can be written as a ratio of smooth occurrence and smooth exposure. An example is given by the local constant Nadaraya-Watson estimator. In Section 4  asymptotic properties are outlined for the smooth backfitting estimator. Details are explained in the Supplementary Material of this paper. Section 5 contains our finite sample study illustrating the strong performance of smooth backfitting. In Section 6 we consider a sophisticated version of in-sample forecasting generalising earlier approaches via our proportional hazard model, see \citet{Mammen:etal:15}, \cite{Hiabu:etal:16} \citet{Lee:etal:15,Lee:etal:17,Lee:etal:18a}. A smooth extension of the popular actuarial chain ladder model that is used in virtually all non-life insurance companies in the world while estimating outstanding liabilities. In-sample forecasting is possible because of the imposed multiplicative structure. 

%%%%%%%%%%%%%%%%%%%%%%%%%%%%%%%%%%%%%%%%%%%%%%%%%%%%%%%%%%%%%%%%%%%%

\section{Aalen's multiplicative intensity model}\label{sec:aalen}
We consider  a counting process formulation satisfying Aalen's multiplicative intensity model. It allows for very general observations schemes.
	It covers filtered observations arising from left truncation and right censoring and also more complicated patterns of occurrence and exposure.
In the next section we describe how to embed left truncation and right censoring into 
this framework.
In contrast to \citet{Linton:etal:03} we will hereby allow the filtering to be correlated to the survival time
and be represented in the covariate process.
We briefly summarize the general  model we are assuming. \\
We observe $n$ $iid$ copies of the stochastic processes $(N(t),Y(t),Z(t)), \ t\in[0,R_0], \, R_0>0.$ Here, $N$ denotes a right-continuous counting process
which is zero at time zero and has jumps of size one.
The process $Y$ is  left-continuous and takes values in $\{0,1\}$ where the value $1$ indicates that the  individual is under risk.
Finally, $Z$ is a $d$-dimensional left-continuous covariate process with values in a hyperrectangle $\prod_{j=1}^d [0,R_j]\subset \mathbb R^d$.%, j=1, \dots, d$.
The multivariate process $((N_1,Y_1,Z_1), \dots, (N_n,Y_n,Z_n))$, $\ i=1,\dots,n$, is adapted to the filtration $\mathcal F_t$ which satisfies  les conditions habituelles (the usual conditions), see \cite{Andersen:etal:93} (pp. 60). Now we assume that $N_i$ satisfies Aalen's multiplicative intensity model, that is
\begin{align}\label{eq:aalen}
\lambda_i(t)=\lim_{h \downarrow 0} h^{-1}E[N_i((t+h)-)-N_i(t-)|\ \mathcal F_{t-}]=\alpha(t,Z_i(t))Y_i(t).
\end{align}
The deterministic function $\alpha(t,z)$ is called hazard function and it is the failure rate  of an individual at time $t$ given the covariate $Z(t)=z$.%$Z_i(t)=z$. 

\subsection{Left truncation and right censoring time as covariates}
\label{sec:model}

The most prominent example for Aalen's multiplicative intensity model is filtered observation due to left truncation and right censoring.
We now show how to embed model \textcolor{red}{\eqref{multhazard1}} with  covariate, $Z$,  possibly carrying truncation and censoring information into Aalen's multiplicative intensity model. Every covariate coordinate can carry individual truncation information as long as it corresponds to left truncation.
That is, we observe $(T,Z)$ if and only if $(T, Z(T))\in\mathcal I$,
where the set $\mathcal I$ is compact and it holds that if
 $(t_1,Z(t_1))\in \mathcal I$  \text{and}   $t_2\geq t_1, \text{then}$  $(t_2,Z(t_2)) \in \mathcal I, \ a.s..$
 The set $\mathcal I$  is allowed to be random but is independent of $T$ given the covariate process $Z$.
Furthermore, $T$ can be subject to right censoring with censoring time $C$.
We assume that  also $T$ and $C$ are conditional  independent given the covariate process $Z$.
This includes the case where the censoring time equals  one covariate coordinate.
In conclusion, we observe $n$ $iid$ copies of $(\widetilde T
, Z^{*}, \mathcal I,\delta)$, where 
$\delta=\mathds 1 (T^{*} <C), \ \widetilde T= \min (T^{*},C)$, and $(T^{*}, Z^{*})$ is the truncated version of $(T,Z)$,  i.e, $(T^*, Z^*)$  arises from $(T,Z)$ by conditioning on the event  $(T,Z(T))\in \mathcal I$.

Then, for each subject, $i=1,\ldots, n$, we can define a counting process $N_i$ as
$
N_i(t)=\mathds 1 \left\{\widetilde T_i \leq t,\ \delta_i=1\right\},
$
with respect to the filtration
$
\mathcal F_{i,t}=\sigma \left( \bigg\{\widetilde T_i\leq s,\  Z^   *_i(s), \ \mathcal I_i, \ \delta_i : \ s\leq t\bigg\} \cup \mathcal N\right),
$
where $\mathcal N$ is a class of null-sets that completes the filtration.
After straightforward computations one can conclude that
under the setting above, Aalen's multiplicative intensity model \eqref{eq:aalen} is satisfied with
\begin{align*}
 \alpha_z(t)&=\alpha(t,z) =\lim_{h \downarrow 0} h^{-1}\mathrm{Pr}\{  T_i\in [t, t+h)| \   T_i\geq t, \ Z_i(t)=z\}, \\
Y_i(t)&=  \mathds 1 \big\{(t,Z^*_i(t))\in \mathcal I_i, \ t\leq \widetilde T_i\big\}.
\end{align*}

\section{Estimation}
\label{sec:estimation}
We tackle the problem in two steps. First the data are projected down onto an unstructured space resulting in an unstructured estimator of the  $d+1-$dimensional hazard function $\alpha(t,z)$; see \eqref{multhazard1}.   In the second step, the unstructured  estimator is projected further down onto the multiplicative space of interest resulting in $d+1$ one-dimensional smooth backfitting estimators of  the multiplicative components, $(\alpha_0(t), \alpha_1(z_1),\dots,\alpha_d(z_d))$.   
For the first step we assume  to have an unstructured estimator with simple ratio of smoothed occurrence and smoothed exposure. 
%This resembles the simple structure known from the local constant Nadaraya-Watson estimator in regression, however, local linear regression does not satisfy this simple structure.
%This is the background for our approach of being able to derive general underlying conditions for our smooth backfitter to work the encompass both local constant and local linear estimators.
Our theory will encompass all estimators with this simple ratio structure. 
We discuss in Section \ref{sec:prop} and in the accompanied Supplementary Material that  our estimation procedure works under quite general assumptions. In particular  we do not need the unstructured estimator of the first step to be consistent. 
The final structured estimator circumvents the curse of dimensionality even if consistency is not assured in the first step.
This is reassuring noting that the  unstructured estimator will most probably have exponentially deteriorating performance with growing dimension $d$.

We introduce the notation $X_i(t)=\left(t,Z_i(t)\right)$. We also set $x=(t,z)$, with coordinates $x_0=t, \  x_1=z_1, \dots, x_d=z_d$, and write the hazard as $\alpha(t,z)=\alpha(x)$.

\subsection{First step: The unstructured estimator}
\label{first step}
%\section{ The smooth backfitting estimator of multiplicative hazards} 
%
% including all local polynomial kernel hazard estimators and other situations, where the unconstrained estimator can be expressed as such simple ratio.
%
%\subsection{First step: projecting the data down onto the unstructured space resulting in an unconstrained estimator}
%
%In this section we concentrate on the local constant  as a projection of the data down onto the unconstrained space. We notice that the estimator can be expressed as a ratio of a smoothed occurrence and a smoothed exposure. This will be important in the next section where the unconstrained estimator is projected further down to the multiplicative space of interest.
%\label{unstructured}

To estimate the components of the structured hazard in \eqref{multstructure} below, we will need an unstructured pilot estimator of the hazard $\alpha$ first.
We propose the local constant kernel estimator, $\widetilde \alpha^{LC}$.
Its value in $x$  is defined  as 
\begin{align}\label{ll}
\widetilde \alpha^{LC}(x)
=\lim_{\varepsilon \to 0} \arg\min_{\theta_0\in \mathbb R}    \sum_{i=1}^n \int &\left\{\frac 1 \varepsilon \int_{s}^{s +\varepsilon} \mathrm{d}N_i(u) -   \theta_0  \right\}^2  \\ & \quad \times \kappa_n(X_i(s)) K_b(x-X_i(s)) Y_i(s)\ \mathrm {d}s. \notag
\end{align}
In the following, we restrict ourselves to a multiplicative kernel, i.e., for  ,
$(u_0,\dots,u_d)\in \mathbb R^{d+1}$, $K(u_0,\dots,u_d)=\prod_{j=0}^d k(u_j)$, and a one-dimensional bandwidth  $b$ with boundary correction
$\kappa_n(x) = \prod_{j=0}^d \left(\int K_b(x_j-u_j)\ du_j\right )^{-1}$ and $K_b(u)=\prod_{j=0}^d b^{-1} k(b^{-1} u_j)$, where for simplicity of notation the bandwidth $b>0$ does not depend on $j$. 
More general choices would have been possible with the cost of extra notation.
%The local linear estimator includes boundary corrections so that the bias is of same order at the boundary as in the interior of the support, namely $O(b^2)$, or for the more general case of varying bandwidths we do not consider here, \ $O(\max_{1 \leq j \leq d+1} b_j^2)$.
%The local constant estimator achieves only slower rates at the boundary region and local polynomial estimators of higher order, like in regression, have the usual drawback known from higher order kernels, that they perform poorly as long as sample sizes are not very large.

%Minimising \eqref{locpoly} for $p=0,1$, straight forward calculations lead to estimators being a ratio  of  smooth estimators of   the number of occurrence and the exposure. 
%In the local linear case, for $p=1$ we get
%\begin{align*}
%\widehat O^{LL}(x)&=n^{-1}\sum_{i=1}^n \int \left\{1-(x-X_i(s))D(x)^{-1}c_1(x)\right\} K_b(x-X_i(s)) \mathrm d N_i(s), \\
%\widehat E^{LL}(x)&=n^{-1}\sum_{i=1}^n \int \left\{1-(x-X_i(s))D(x)^{-1}c_1(x)\right\} K_b(x-X_i(s)) Y_i(s) \mathrm ds,
%\end{align*}
%where $D(x)=(d_{jk}(x))_{jk}$ is a $(d+1)\times (d+1)-$dimensional matrix defined in the precious section.
%The local linear estimator is then defined as $\widehat \alpha^{LL}(x)= \widehat O^{LL}(x)/\widehat E^{LL}(x)$. Compare this estimator with 
We get $\widetilde  \alpha^{LC}(x)= \widehat O^{LC}(x)/\widehat E^{LC}(x)$, with smoothed occurrence and smoothed exposure given by
 \begin{align*}
\widehat O^{LC}(x)&=  \sum_{i=1}^n \int \kappa_n(X_i(s)) K_b(x-X_i(s)) \mathrm d N_i(s), \\
\widehat E^{LC}(x)&=  \sum_{i=1}^n  \int  \kappa_n(X_i(s)) K_b(x-X_i(s)) Y_i(s) \mathrm ds.
\end{align*}

Under standard smoothing conditions, if $b$ is  chosen of order $ n^{-1/(4+d+1)}$, then  the bias of  $\widetilde \alpha^{LC}(x)$ is of order $n^{-2/(4+d+1)}$ and the variance is of order $n^{-4/(4+d+1)}$,
 which is the optimal rate of convergence in the corresponding regression problem, see \citet{Stone:82}.
For an asymptotic theory of these estimators see \citet{Linton:etal:03}.

\subsection{Second step: The structured smooth backfitting estimator}
In this section we will project the unstructured estimator of the previous section down onto the multiplicative space of interest. Other choices that have a simple ratio structure of occurrence and exposure are possible.
Due to filtering, observations are assumed to be only available on a subset of the full support, $\mathcal X \subseteq \mathcal  R=\prod_{j=0}^d[0, R_j]$.
Our estimators are restricted to this set and detailed assumptions on $\mathcal X$ and the data generating functions are given in the Supplementary Material.
  Our calculations simplify via a new principle we call solution-weighted minimization. We assume that we have the solution and use it strategically in the least squares weighting. While the definition is not explicit, it is made feasible by defining it as an iterative procedure.
In the sequel we will assume a multiplicative structure of the hazard $\alpha$, i.e.,
\begin{align}\label{multstructure}
\alpha(x)=\alpha^* \prod_{j=0}^d \alpha_j(x_j),
\end{align}
where $\alpha_j$, $j=0,\dots,d,$ are some functions and $\alpha^*$ is a constant.
For  identifiability of the components, we  make the following further assumption:
\begin{align*}%\label{qprop1}
\int \alpha_j(x_j) w_j(x_j)\ \mathrm{d} x_j=1, \quad j=0,\dots,d,
\end{align*}
where the $w_j$'s  are some weight functions. 

We also need the following notation:
\begin{align*}
F_t(z)=Pr\left(Z_1(t) \leq z |\ Y_1(t)=1\right), \quad
y(t)=E[Y_1(t)]. 
\end{align*}
By denoting $f_t(z)$ the density corresponding to $F_t(z)$ with respect to the Lebesgue measure, we also define
$
E(x)=f_t(z)y(t)
$
and $O(x)=E(x)\alpha(x)$.

We define the estimators $\widehat \alpha^*$ and $\widehat \alpha=( \widehat \alpha_0,\dots, \widehat \alpha_d)$ of the  hazard components  in \eqref{multstructure} as solution of the following system of equations: 
\begin{align}\label{backfitting2}
&\widehat \alpha_k(x_k)=\frac{\int_{\mathcal X_{x_k}}\widehat O(x)\mathrm dx_{-k}}
{\int_{\mathcal X_{x_k}} \widehat \alpha^*  \prod_{j\neq k}\widehat \alpha_j(x_j)\widehat E(x)\mathrm dx_{-k}},  \quad k=0,\dots, d, \\
&\label{qprop1} \int \widehat \alpha_k(x_k) w_k(x_k)\ \mathrm{d} x_k=1, \quad k=0,\dots,d.
\end{align}

Here $\mathcal X_{x_k}$  denotes the set $\{(x_0, \dots, x_{k-1},x_{k+1},\dots, x_d)|\ (x_0,\dots,x_d) \in \mathcal X\}$, and\\  $x_{-k}=(x_0, \dots, x_{k-1},x_{k+1},\dots, x_d)$. Furthermore, $\widehat E$ and $\widehat O$ are some full-dimensional estimators of $E$ and $O$ -- not necessarily the one provided in the previous section. We will discuss  in the Supplementary Material that the system above has a solution with probability tending to one. In the next section and  in the Supplementary Material we will show asymptotic properties of the estimator. We will see that we do not require that the full-dimensional estimators $\widehat E$ and $\widehat O$ are consistent. We will only need asymptotic consistency of marginal averages of the estimators. This already highlights that our estimator efficiently circumvents  the curse of dimensionality.

In practice, system \eqref{backfitting2} can be solved by the following iterative procedure:
\begin{align}\label{iteration}
\widehat \alpha_k^{(r+1)}(x_k)= \frac{\int_{\mathcal X_{x_k}}\widehat O(x)\mathrm dx_{-k}}
{\int_{\mathcal X_{x_k}} \prod_{j=0}^{k-1}\widehat \alpha_j^{(r+1)}(x_j) \prod_{j=k+1}^{d+1}\widehat \alpha_j^{(r)}(x_j)\widehat E(x)\mathrm dx_{-k}},  \quad k=0,\dots, d
\end{align}
After a finite number of cycles or after a termination criterion applies, the last values of $\widehat \alpha_k^{(r+1)}(x_k)$, $k=0,...,d$, are multiplied by a factor such that the constraint \eqref{qprop1} is fulfilled.
This can always be achieved by multiplication with constants. This gives the backfitting approximations of $\widehat \alpha_k(x_k)$ for  $k=0,...,d$.
\subsection{Interpretation as direct projection}
The strength of our smooth backfitting estimator is that 
it can be motivated directly from a least squares criterium on the data without ad-hoc adjustment.
The estimator $\widehat \alpha$ can be motivated as solution of 
\begin{align}\label{locpolyproj}
%\begin{pmatrix} \widehat \theta_0 \\ \widehat \theta_1\\ \end{pmatrix}=
\lim_{\varepsilon \to 0} \arg\min_{\theta}    \sum_{i=1}^n \int \int &\Bigg\{\frac 1 \varepsilon \int_{s}^{s +\varepsilon} \mathrm{d}N_i(u) -   \theta(x)\Bigg\}^2 
&  \times  K_b(x-X_i(s)) Y_i(s)\ \mathrm {d}s \ \mathrm d\nu(x),\end{align}
where $\theta$ runs over some space of smooth multiplicative functions of the form
$\theta=\prod_{j=0}^d \theta_j(x_j)$.
 To see this, consider the estimator $\overline \alpha=(\overline \alpha^*,\overline \alpha_0, \dots, \overline \alpha_d)$ that minimizes
\begin{align}\label{min}
\arg \min_{\overline \alpha} \int_{\mathcal X} \left\{ \widetilde \alpha^{LC}(x) - \overline \alpha^*\prod_{j=0}^d \overline \alpha_j(x_j)\right\}^2 w(x) \mathrm d x.
\end{align}
For $\nu(x)=w(x)/\widehat E^{LC}(x)dx$ the solution of \eqref{min} is exactly \eqref{locpolyproj}.
With that choice, we get $$\overline \alpha^*=\frac{\int_{\mathcal X} \widetilde \alpha(x) \prod_{j=0}^d\overline\alpha_j(x_j)w(x)\mathrm dx}
{\int_{\mathcal X} \left \{ \prod_{j=0}^d\overline \alpha_j(x_j)\right \}^2 w(x)\mathrm dx},$$
and $( \overline \alpha_0,\dots, \overline \alpha_d)$ can be described via the backfitting equation
\begin{align}\label{backfitting1}
&\overline \alpha_k(x_k)=\frac{\int_{\mathcal X_{x_k}} \widetilde \alpha(x) \prod_{j\neq k}\overline\alpha_j(x_j)w(x)\mathrm dx_{-k}}
{\int_{\mathcal X_{x_k}}\overline \alpha^* \left \{ \prod_{j\neq k}\overline \alpha_j(x_j)\right \}^2 w(x)\mathrm dx_{-k}},  \quad k=0,\dots, d.
\end{align}
 The asymptotic variance of kernel estimators of $\alpha$ is proportional to $\alpha(x)/E(x)$, see e.g.\ \citet{Linton:Nielsen:95}. This motivates the choice 
 $w(x)=E(x)/\alpha(x)$. However,
this choice is not possible because  $E(x)$ and $\alpha(x)$ are unknown.  One could use $w(x)= \check E(x)/\check \alpha(x)$ where $\check E(x)$ and $\check \alpha(x)$ are some pilot estimators of $E$ and $\alpha$. We follow another idea and we propose to weight the minimization \eqref{min} with its solution.
We choose \begin{align}\label{wq}
w(x)&= \frac{\widehat E(x)}{\prod_i \widehat \alpha_i(x)},
\end{align}
and heuristically, by putting $\overline \alpha_j = \widehat \alpha_j$ and  by plugging \eqref{wq}
 into \eqref{backfitting1},  we get \eqref{backfitting2}.
\section{Asymptotic properties of the smooth backfitter of multiplicative hazards}
\label{sec:prop}
In the Supplementary Material we show that $\widehat \alpha_j$ converges to the true $\alpha_j$
with optimal one dimensional nonparametric rate  $n^{-2/5}$, given that the bandwidth $b$ is chosen of order $n^{-1/5}$.
This means in particular that the asymptotic rate does not depend on the dimension $d$.
Under regularity assumptions, in Theorem 3 in the Supplementary Material we show that
\begin{align*}
 n^{2/5}\{(\widehat \alpha_j- \alpha_j)(x_j) -  \alpha_j(x_j) B_j(x_j))\} \rightarrow \mbox{N}(0, \alpha^2_j(x_j)\sigma_j^2(x_j)),\end{align*}
where  $B_j(x_j)$ can be formally defined as the $j$-th component of a projection of the bias  of the unconstrained estimator onto the multiplicative space
 and  $\sigma_j^2(x_j)$ is the variance of a weighted average of  the unconstrained estimator  where the other components are integrated out.

In the  simulation study of the next section we show that estimation seems to work well even when $d$ is of similar order as $n$. This is, we believe, a major strength of our smooth backfitting estimator.

A few notes on the proof.
The estimator $\widehat \alpha_j$ is defined as solution of a nonlinear operator equation. 
Asymptotic properties are derived by approximating the estimators of this equation by a linear equation that can be interpreted as one equation that arises in nonparametric additive regression models \citep{Mammen:etal:99}, and then one shows that the solution of the linear equation  approximates the estimation error in the multiplicative model. The linear equation and its solution is well understood from the theory of additive models. This is  our essential step to arrive at an asymptotic understanding of our estimator $\widehat \alpha_j$. Assumptions   are of standard nature in marker dependent hazard papers and they can be verified for the local constant  estimators we are interested in, see in particular \citet{Nielsen:Linton:95}, \citet{Nielsen:98} and \citet{Linton:etal:03} for related arguments. However, the conditions are formulated more general and they are not restricted to the local constant smoothers. They are not even tight to kernel smoothers. Any smoother could be used as long as it obeys the structure of being a ratio of a smoothed occurrence and a smoothed exposure. 
\section{Simulation study}\label{sec:sim}
In this section we present detailed simulations of the smooth backfitter. The Supplementary Material contains additional simulation results for  the setting discussed in \cite{Honda:05} and  \cite{Lin:etal:16}. In these latter settings, the performance of the smooth backfitter is similar to the estimator of \cite{Lin:etal:16}. And both estimators outperform the estimators considered in \cite{Honda:05}.

In the simulation study below, and with additional results in The Supplementary Material, we compare our estimator with the estimator of \cite{Lin:etal:16}. The models of these simulations contain high-dimensional and correlated covariates and in particular, we will show that our estimator -- in contrast to \cite{Lin:etal:16} -- works in high dimensions where the dimension of the covariates is of similar order as the sample size and also in cases with higher correlation between the covariates.
To highlight the impact of increasing dimension we will consider  the dimensions $d=2$ and $d=9$.
Further results for dimensions $d=30$ and $d=99$ can be found in the Supplementary Material.
The low dimensional setting, $d=2$, with uncorrelated covariates is similar to the setting of \cite{Honda:05} and  \cite{Lin:etal:16}, and the difference between our estimator and the estimator of \cite{Lin:etal:16} are indeed marginal in  this case.
However, below and in  the Supplementary Material we show that this changes drastically with increasing dimension and/or correlation.

\subsection{The setting}
Since the estimator of \cite{Lin:etal:16} is based on a partial likelihood approach, it does not estimate the baseline hazard, $\alpha_0$.  We will consider the sub-model
\[
\alpha(x)=\prod_{k=0}^{d}\alpha_k(x_k)=\exp\left\{\sum_{k=1}^{d} \eta_k(x_k)\right\},
\]
i.e., we assume a constant baseline hazard, $\alpha_0\equiv1$.
More specifically, we assume that the survival times $T_i$  follow an exponential distribution with parameter value $\alpha(X_i), \quad (i=1,\dots,n)$.
We add right censoring  with censoring variables $C_i$ that follow an exponential distribution with parameter 
$\frac 4 7 \alpha(X_i), \quad (i=1,\dots,n)$.
We will compare the estimators for $\eta_k, k=1,\dots,d$. Our proposed estimator is derived as $\widehat \eta_k^{SBF}=\log(\widehat \alpha_k)$ and we compare it to, $\widehat \eta_k^{\text{Lin et al.}}$, proposed in \cite{Lin:etal:16}.

We used the following two models:
\begin{align*}
&\text{Model 1: } \eta_k(z_k)=\begin{cases}-z_k \quad &\text{if k is odd,}\\ 2z_k \quad &\text{if k is even.}\end{cases}\\
&\text{Model 2: } \eta_k(z_k)=\begin{cases}2\sin(\pi z_k) \quad &\text{if k is odd,}\\ 2z_k \quad &\text{if k is even.}\end{cases}
\end{align*}
The covariates $(Z_{i1}, \dots, Z_{id})$ are generated as follows.
We first  simulate $(\widetilde Z_{i1}, \dots, \widetilde Z_{id})$ from a $d$-dimensional  multi-normal distribution with mean equal
0 and $\textrm{Corr}(Z_{ij},Z_{il})=\rho$ if $j\neq l,$ else 1.
Afterwards we set
\[
        Z_{ik}=2.5\pi^{-1}\text{arctan}(\widetilde Z_{ik}).
\]
This has been independently repeated for every individual $i=1,\dots,n$.
%We will consider combinations of the following parameter settings:
%\begin{align*}
%n= 200, 500,  \quad d=2,99  \quad \rho=0,0.5,0.8.
%\end{align*}
After trying several bandwidths,  if not said otherwise, we present the results for a bandwidth $b=0.3$ for both estimators in every model. As kernel function k, we used the Epanechnikov kernel. Performance is measured via the integrated squared error evaluated at the observed points:
\begin{align*}
ISE_k=n_{obs}^{-1}\sum_{i_{obs}} \left(\eta_{k}(Z_{i_{obs}k})-\widehat \eta_{k}(Z_{i_{obs}k})\right)^2,\\
ISE_{odd}=\{ISE_k| k=odd\},\quad ISE_{even}=\{ISE_k| k=even\},\quad 
\end{align*}
where $i_{obs}$ runs over all observed individuals and $n_{obs}<n$ is the number of observed individuals.
The next two subsections discuss the outcomes of the simulations for $d=2$ and $d=9$.
\subsection{Dimension $d=2$}
\begin{figure}
%\centering
%\begin{subfigure}[t]{0.2\textwidth}
%\centering
\includegraphics[width=16cm]{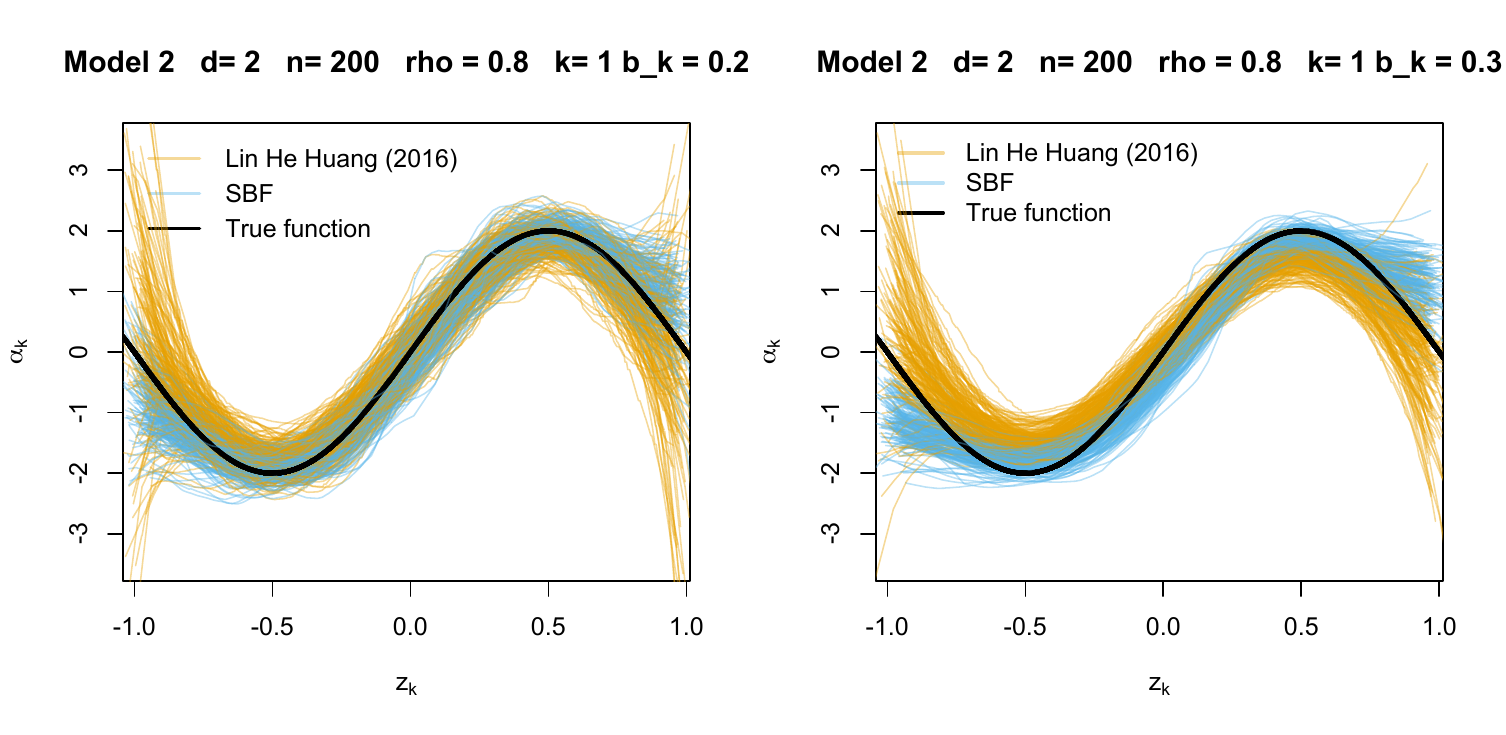} 
%\caption{Generic} \label{fig:timing1}
%\end{subfigure}
\caption{Estimates of $\alpha_1$ from 200 simulations. Left panel uses a bandwidth of 0.2 and right panel a bandwidth of 0.3. }
  \label{fig:sim1}
\end{figure}

\begin{figure}
%\centering
%\begin{subfigure}[t]{0.2\textwidth}
%\centering
\includegraphics[width=14cm]{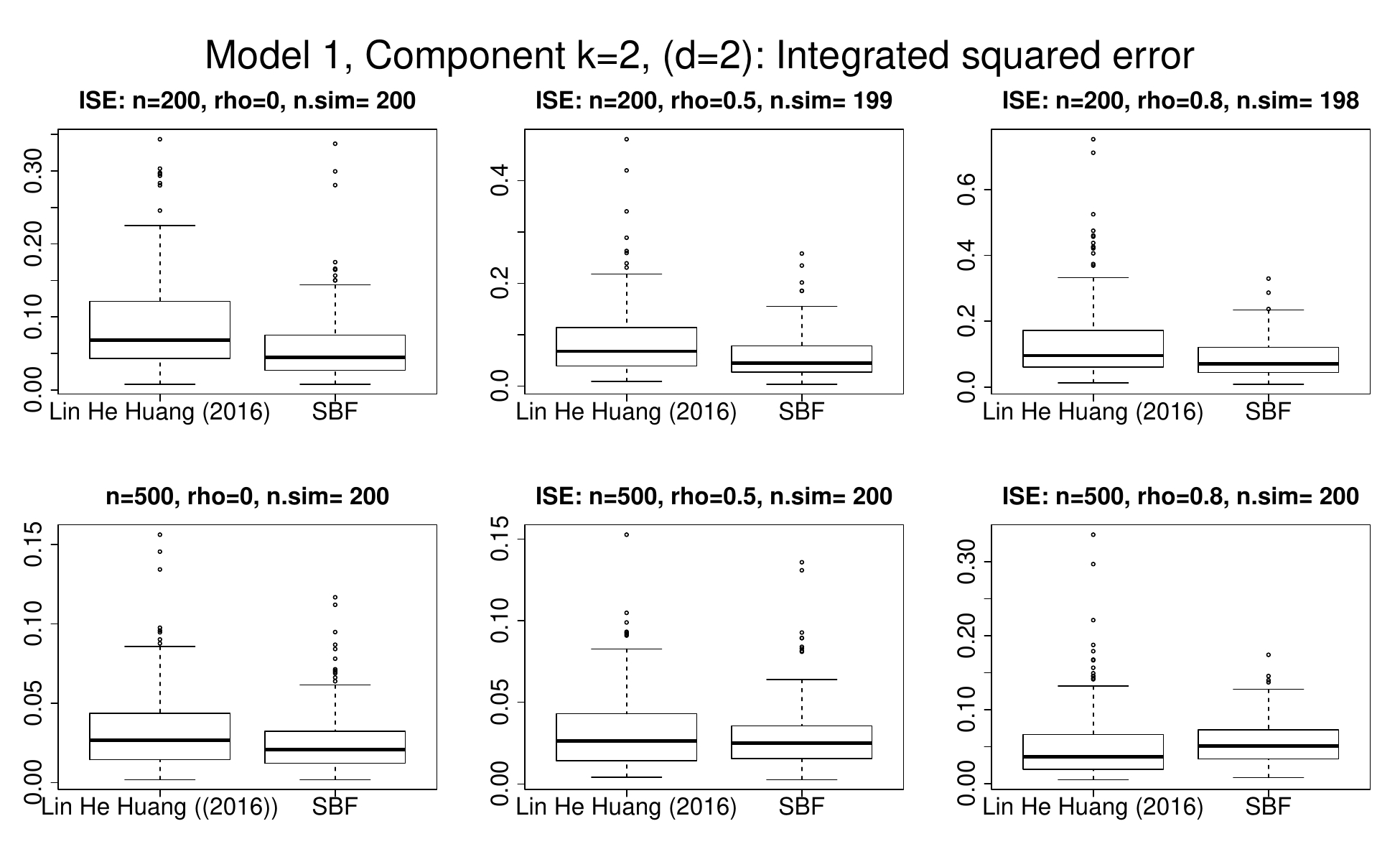} 
%\caption{Generic} \label{fig:timing1}
%\end{subfigure}

%\begin{subfigure}[t]{0.2\textwidth}
%\centering
\includegraphics[width=14cm]{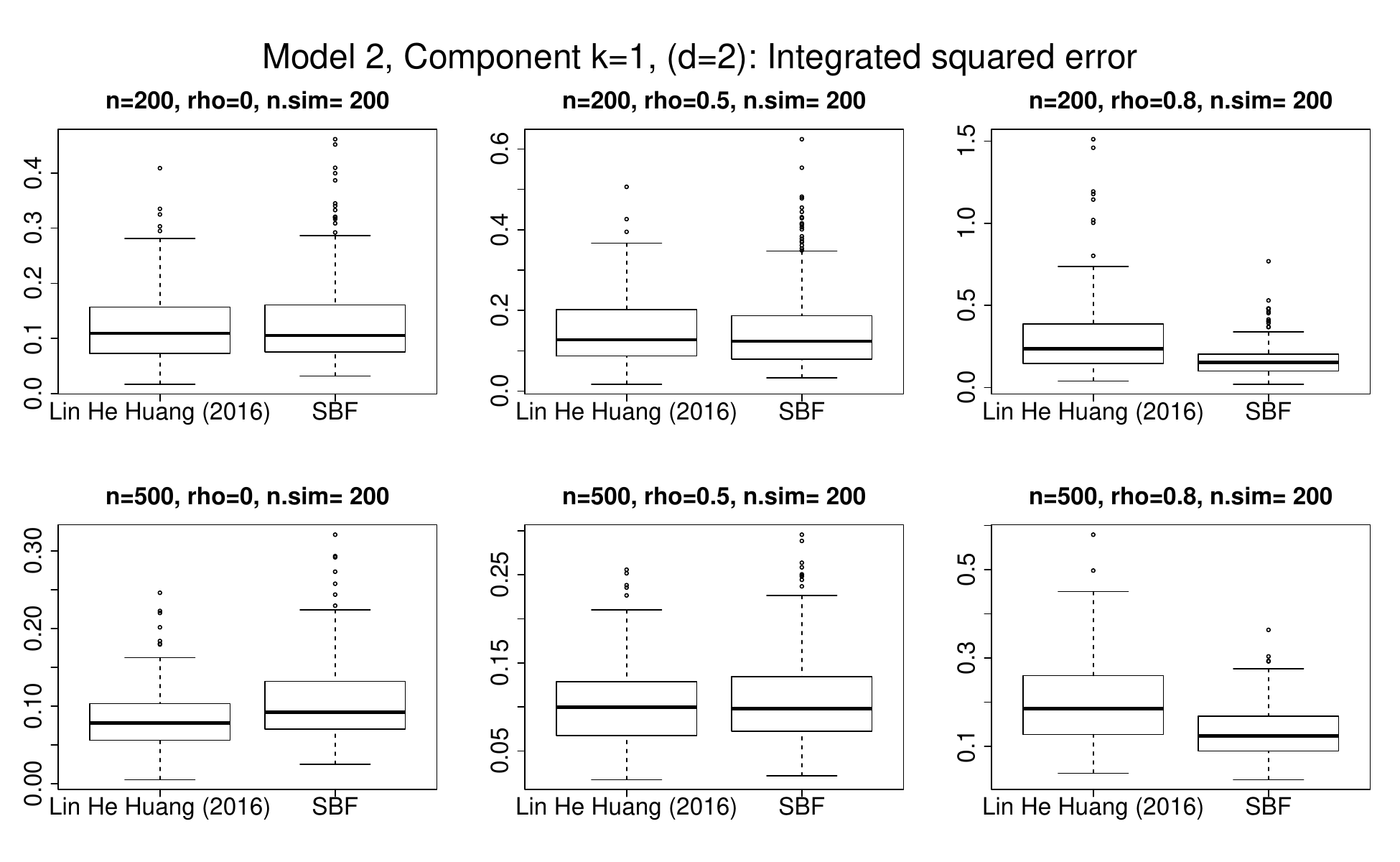} 
%\caption{Competitors} \label{fig:timing2}
%\end{subfigure}
\caption{Boxplots of the integrated squared errors. Simulations are taken out where the algorithm for the calculation of the estimator in  \cite{Lin:etal:16} stopped without calculation of all values of the estimator. The value $n.sim$ is the number of the remaining simulations, i.e.,  200 minus number of break downs.}
\label{fig:boxplot:ise:d3}
\end{figure}

In low dimensions both estimators perform satisfactorily.
Already  in dimension $d=2$, for Model 1, it occurred that the algorithm of \cite{Lin:etal:16}  stopped before the final calculation of the estimator.  This happened for three simulation samples out of 1200  ($6$ settings $\times 200$ samples each) runs.  
The algorithm stopped at a step where  a matrix has to be inverted that is nearly singular.
In contrast, our estimator does not need to invert matrices and we did not have any singularity or convergence issues for our estimator.
%Another issue both estimators share is the case when the iterations do not converge.
%This later issue will be covered in a seperate section covering the computational complexity, Section \ref{sec:sim:complex}.
%This did not happen in this low dimensional case.
Performance-wise, we refer to Figure \ref{fig:boxplot:ise:d3}  displaying boxplots of the integrated squared errors of both estimators. The two estimators perform similar with a small advantage towards our smooth backfitting estimator.
The smooth backfitting estimator performs especially better for the smaller sample size $n=200$. 
While  \cite{Lin:etal:16}  performed better in Model 2 with no correlation,
the performance of the smooth backfitting estimator seems more stable when correlation is added.
Figure \ref{fig:sim1}  shows the 200 sample estimates of the first component $\alpha_1$ of Model 2 with $\rho = 0.8$.
The estimator of  \cite{Lin:etal:16}  struggles especially at the boundaries and this is more pronounced with a small bandwidth (left panel). If bandwidth is increased (right panel), the estimate of   \cite{Lin:etal:16}  seems  over-smoothed and it is not able to replicate the full magnitude of the local extrema at $-0.5$ and $-0.5$.
%\begin{figure}
%\begin{subfigure}{.5\textwidth}
%  \centering
%  \includegraphics[width=.8\linewidth]{Figure.pdf}
%  \caption{}
%  \label{fig:sfig1}
%\end{subfigure}%
%\begin{subfigure}{.5\textwidth}
%  \centering
%  \includegraphics[width=.8\linewidth] {Figure1.pdf}
%  \caption{}
%  \label{fig:sfig2}
%\end{subfigure}
%\begin{subfigure}{.5\textwidth}
%  \centering
%  \includegraphics[width=.8\linewidth] {Figure2.pdf}
%  \caption{}
%  \label{fig:sfig2}
%\end{subfigure}
%\begin{subfigure}{.5\textwidth}
%  \centering
%  \includegraphics[width=.8\linewidth] {Figure3.pdf}
%  \caption{}
%  \label{fig:sfig2}
%\end{subfigure}
%\begin{subfigure}{.5\textwidth}
%  \centering
%  \includegraphics[width=.8\linewidth] {Figure4.pdf}
%  \caption{}
%  \label{fig:sfig2}
%\end{subfigure}
%\begin{subfigure}{.5\textwidth}
%  \centering
%  \includegraphics[width=.8\linewidth] {Figure5.pdf}
%  \caption{}
%  \label{fig:sfig2}
%\end{subfigure}
%\caption{Plots of six simulations. Left column: Three simulation are picked: (from top to bottom) the 20th best (a), 100th best (c) and 180th best (e)  performance of the \cite{Lin:etal:16} estimator measured via ISE. In each of those three  simulations
%we compare the performance between the smooth backfitting estimator (SBF) and the estimator of  \cite{Lin:etal:16}. 
%Right column: Analogue for the SBF estimator.
%}
%\label{fig1}
%\end{figure}

%\begin{figure}[h!]
%\includegraphics[width=\textwidth]{boxplot_ise321.pdf} 
%%\caption{Generic} \label{fig:timing1}
%%\end{subfigure}
%
%%\begin{subfigure}[t]{0.2\textwidth}
%%\centering
%\includegraphics[width=\textwidth]{boxplot_ise322.pdf} 
%\caption{Boxplots of the integrated squared errors. Simulations where the algorithm of  Lin He Huang (2016) broke down are taken out. The value $n.sim$ is the number of successful simulations, i.e.,  200 minus number of break downs.}
%\label{fig:boxplot:ise:d32}
%\end{figure}

%\begin{subfigure}[t]{0.2\textwidth}
%\centering

\subsection{Dimension $d=9$}
\begin{figure}[!ht]
%\centering
%\begin{subfigure}[t]{0.2\textwidth}
%\centering
\includegraphics[width=16cm]{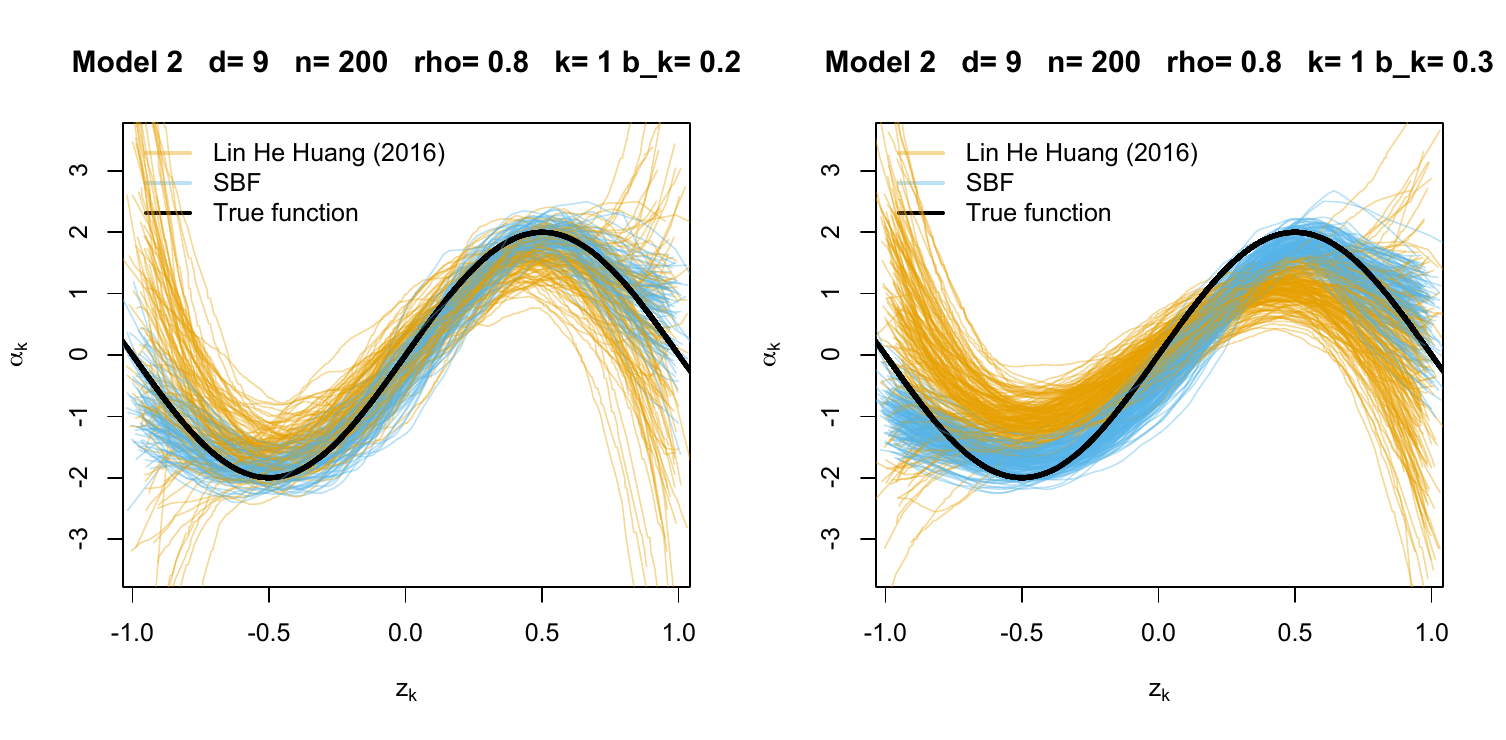} 
%\caption{Generic} \label{fig:timing1}
%\end{subfigure}
\caption{Estimates of $\alpha_1$ from 200 simulations. Left panel uses a bandwidth of 0.2 and right panel a bandwidth of 0.3. }
  \label{fig:sim2}
\end{figure}

\begin{figure}
%\centering
%\begin{subfigure}[t]{0.2\textwidth}
%\centering
\includegraphics[width=14cm]{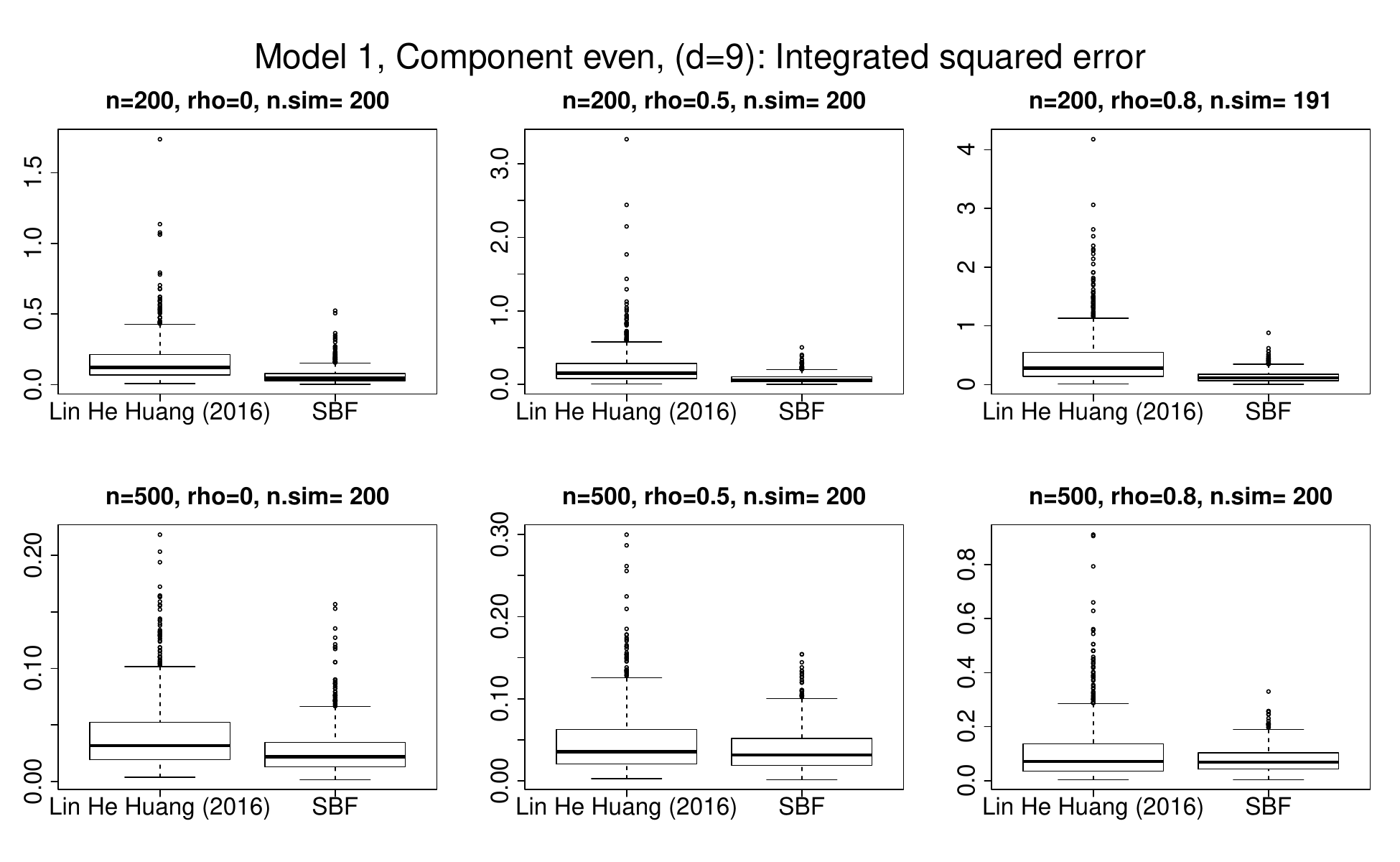} 
%\caption{Generic} \label{fig:timing1}
%\end{subfigure}

%\begin{subfigure}[t]{0.2\textwidth}
%\centering
\includegraphics[width=14cm]{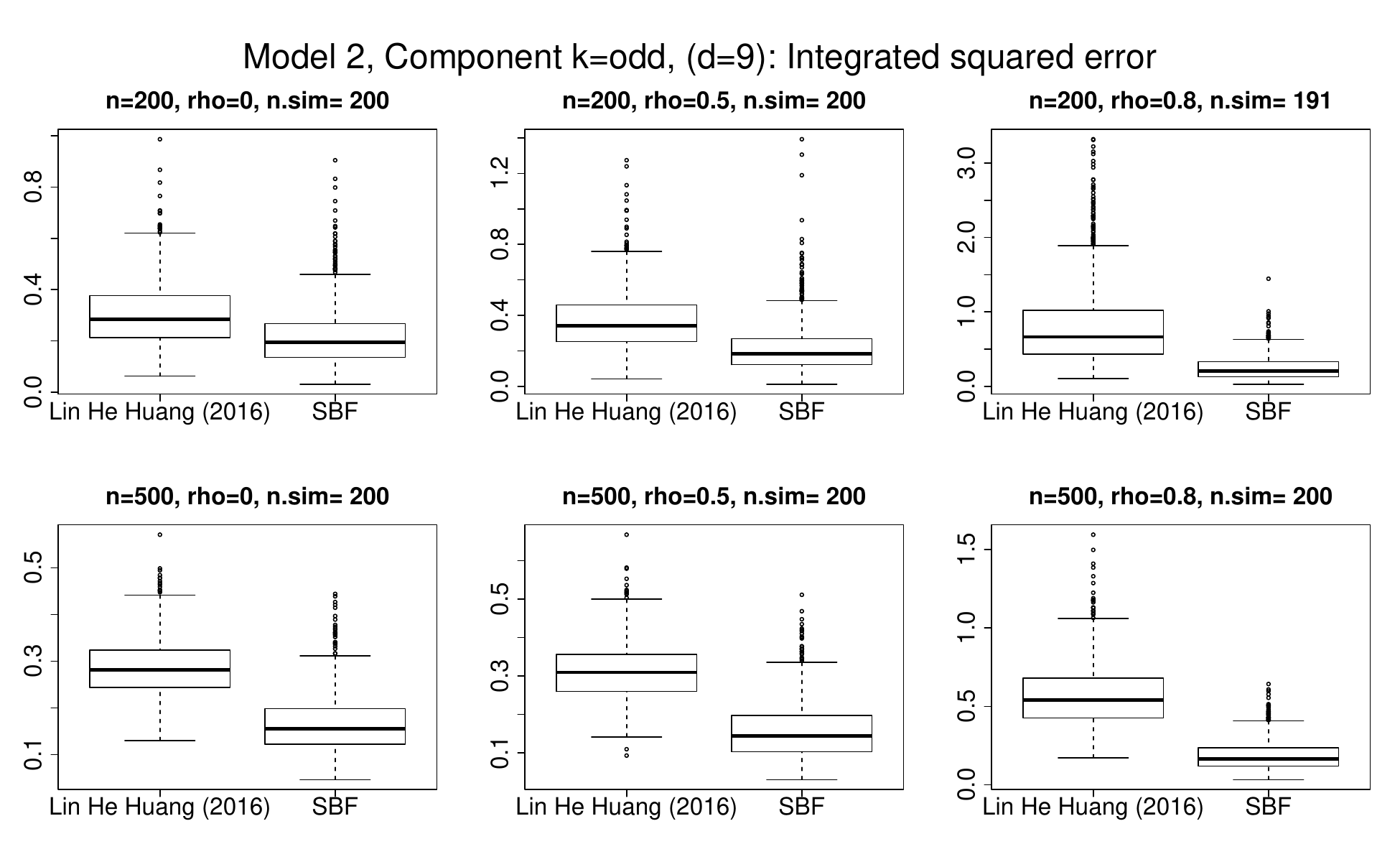} 
%\caption{Competitors} \label{fig:timing2}
%\end{subfigure}
\caption{Boxplots of the integrated squared errors. Simulations are taken out where the algorithm for the calculation of the estimator in  \cite{Lin:etal:16} stopped without calculation of all values of the estimator. The value $n.sim$ is the number of the remaining simulations, i.e.,  200 minus number of break downs.}
\label{fig:boxplot:ise:d9}
\end{figure}

When the dimension $d$ is increased to 9, the estimator of \citet{Lin:etal:16}  breaks down considerably more often than in the case $d=2$, i.e., in 59 (4+5+48+1+1 out of 1200) cases in Model 1 and  9 times in Model 2; see Table \ref{tab:breakdownd:d9}.
Note that in the more extreme cases of $d=99$, provided in the Supplementary Material, nearly all simulations of \citet{Lin:etal:16}  (780 and 607 out of 800 for Model 1 and Model 2, respectively) broke down.
In contrast, our estimator converged in all cases considered.
Performance-wise we refer to Figure \ref{fig:boxplot:ise:d9} displaying boxplots of the integrated squared errors of both estimators. We see that for $d=9$ the smooth backfitter performs better in every setting. The better performance  is more  pronounced when more correlation is present. Figure \ref{fig:sim2}
shows  200 sample estimates of the first component $\alpha_1$ of Model 2 with $\rho = 0.8$.
The results are similar as in the case $d=2$, but now more pronounced: The estimator of  \cite{Lin:etal:16}  struggles at the boundaries and if bandwidth is increased (right panel), the estimate of   \cite{Lin:etal:16}  seems too smooth and is not able to replicate to full magnitude of the local extrema at $-0.5$ and $-0.5$.

\begin{table}
\centering
\begin{tabular}{|c|ccc|ccc|}
\hline
\multicolumn{7}{|c|}{Number of breakdowns in \cite{Lin:etal:16} for $d=9$} \\ \multicolumn{7}{|c|} {(out of 200 simulations)} \\
\hline
 &\multicolumn{3}{c|}{Model 1}&\multicolumn{3}{c|}{Model 2}\\
  \hline
&$\rho=0$ & $\rho=0.5$ & $ \rho=0.8$  &  $\rho=0$ & $\rho=0.5$ & $\rho=0.8$\\
  \hline
n=200 & 4 &5 & 48& 0&0&9 \\ 
  n=500 & 0 &1 & 1& 0&0&0 \\ 
  \hline
\end{tabular}
\caption{Number of breakdowns in  the algorithm of  \cite{Lin:etal:16} out of 200 simulations for dimension $d=9$.}
\label{tab:breakdownd:d9}
\end{table}

%\begin{figure}
%\begin{subfigure}{.5\textwidth}
%  \centering
%  \includegraphics[width=.8\linewidth]{Figure6.pdf}
%  \caption{}
%  \label{fig:sfig1}
%\end{subfigure}%
%\begin{subfigure}{.5\textwidth}
%  \centering
%  \includegraphics[width=.8\linewidth] {Figure7.pdf}
%  \caption{}
%  \label{fig:sfig2}
%\end{subfigure}
%\begin{subfigure}{.5\textwidth}
%  \centering
%  \includegraphics[width=.8\linewidth] {Figure8.pdf}
%  \caption{}
%  \label{fig:sfig2}
%\end{subfigure}
%\begin{subfigure}{.5\textwidth}
%  \centering
%  \includegraphics[width=.8\linewidth] {Figure9.pdf}
%  \caption{}
%  \label{fig:sfig2}
%\end{subfigure}
%\begin{subfigure}{.5\textwidth}
%  \centering
%  \includegraphics[width=.8\linewidth] {Figure10.pdf}
%  \caption{}
%  \label{fig:sfig2}
%\end{subfigure}
%\begin{subfigure}{.5\textwidth}
%  \centering
%  \includegraphics[width=.8\linewidth] {Figure11.pdf}
%  \caption{}
%  \label{fig:sfig2}
%\end{subfigure}
%\caption{Plots of six simulations. Left column: Three simulation are picked: (from top to bottom) the 20th best (a), 100th best (c) and 180th best (e)  performance of the \cite{Lin:etal:16} estimator measured via ISE. In each of those three  simulations 
%we compare the performance between the smooth backfitting estimator (SBF) and the estimator of  \cite{Lin:etal:16}. 
%Right column: Analogue for the SBF estimator.
%}
%\label{fig2}
%\end{figure}

%%%%%%%%%%%%%%%%%%%%%%%%%%%%%%%%%%%
%%%%%%%%%%%%%%%%%%%%%%%%%%%%%%%%%%%
%%%%%%%%%%%%%%%%%%%%%%%%%%%%%%%%%%
%%%%%%%%%%%%%%%%%%%%%%%%%%%%%%%%%%%
%%%%%%%%%%%%%%%%%%%%%%%%%%%%%%%%%%%

\section{In-sample forecasting of outstanding loss liabilities}
\label{sec:reserving}

The so-called chain ladder method is a popular approach to estimate outstanding liabilities. It  started off as a deterministic algorithm, and it is used today for almost every single insurance policy over the world in the business of non-life insurance.  In many developed countries, the non-life insurance industry has revenues amounting to around 5\%. It  is therefore comparable to - but smaller than - the banking industry. In  every single product sold, the chain ladder method (because actuaries hardly use other methods) comes in,  estimating the outstanding liabilities that eventually aggregate to the reserve - the single  biggest  number of most non-life insurers balance sheets. The insurers liabilities often amount to many times the underlying value of the company.  In Europe alone those outstanding  liabilities are estimated to accumulate to  around \EUR$1$trn. It is therefore of obvious importance that this estimate is not too far from the best possible estimate. We describe in this section how the methodology introduced in this paper can be applied to provide a solution to this challenging problem. 

\begin{figure}
\centering
\includegraphics[width=8cm]{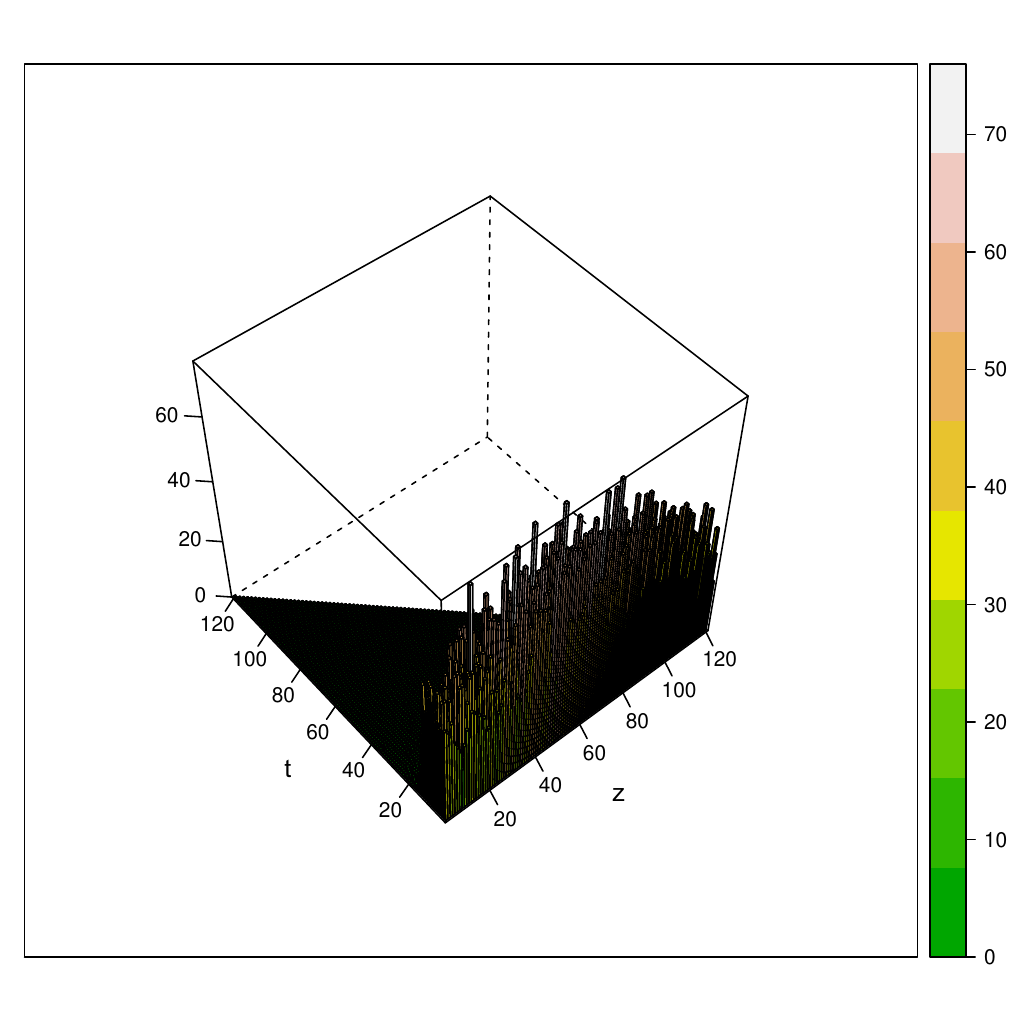} 
\caption{\label{histo} Histogram of claim numbers of a motor business line between 2004  and 2013. Axis $z$ represents the underwriting time (in months) and axis $t$ the reporting delay (in months).}
\end{figure}

We analyze reported  claims from a motor business line in Cyprus.
The same data set has been used by \citet{Hiabu:etal:16} and it consists of the number of claims reported between 2004 and 2013. During these $10$ years (3654 days), $n=58180$  claims were reported.
The data are given as $\{(T_1,Z_1),\ldots, (T_n,Z_n)\}$, where $Z_i$ denotes the underwriting date of claim $i$,  and $T_i$ the time between underwriting date and the date of report of a claim  in days, also called  reporting delay.
Hence, in the notation of the previous sections, the covariate underwriting date, $Z(t)=Z$,  does not depend on time and has dimension $d=1$.
The data  exist on a triangle, with $T_i + Z_i \leq 31 \text{ December } 2013=R_0$, which is a subset of the full support 
$\mathcal R=[0,R_0]^2$ ($0=1 \ \text{January} \  2004$). 
The  aim is to forecast the number of future claims from contracts written in the past which have not been reported yet.   Figure \ref{histo} shows the observed data that lie on a triangle, while the forecasts are required on the triangle that added to the first completes a square. Here it is  implicitly  assumed that the maximum reporting delay of a claim is $10$ years.
Actuaries call this assumption that the triangle is fully run off. In our data set, this is a reasonable assumption looking at  Figure \ref{histo}.

The classical chain ladder method is able to provide a simple solution to the above problem. Recently, \citet{Martinez:etal:13} have pointed out that this method can be viewed as a multiplicative density method:
the  original, un-truncated random variable $(T,Z)$ having density $f(t,z)=f_1(t)f_2(z)$;
and the authors suggested to embed the method in a more standard mathematical statistical vocabulary to engage mathematical statisticians in future developments. In particular, \citet{Martinez:etal:13} showed that one could consider the traditional chain ladder estimator as a multiplicative histogram in a continuous framework, and presented an alternative by projecting an unconstrained local linear density down onto a multiplicative subspace. This approach was called continuous chain ladder and it has been further analyzed by \citet{Mammen:etal:15, Lee:etal:15, Lee:etal:17}, providing full asymptotic theory of the underlying density components. A related approach by \citet{Hiabu:etal:16, Hiabu:17} proposes to transform the two-dimensional multiplicative continuous chain ladder problem to two one-dimensional continuous hazard estimation problems via an elegant time-reverting trick. The application considered in this paper generalizes the most important of these reversed hazards to a two-dimensional multiplicatively structured hazard. In this way the continuous chain ladder is improved and generalized allowing more flexibility for the estimation of outstanding liabilities in the insurance business. 

\begin{figure} 
% ($t$ fixed). }
\centering
\makebox{  \includegraphics[width=13cm]{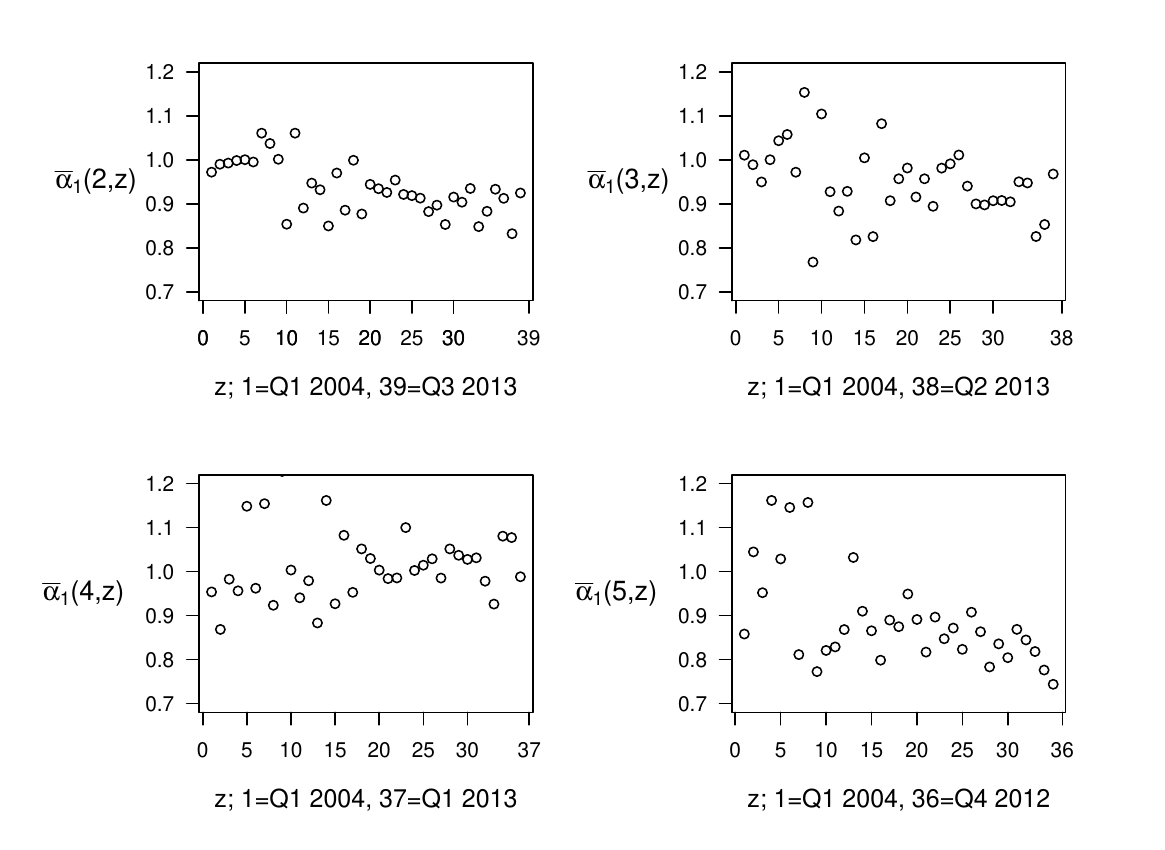}  }
\caption{\label{Fig:devfactors}Scaled quarterly hazard rates of the first four development quarters.}
\end{figure}

In \citet{Hiabu:etal:16} it is assumed that $T$ and $Z$ are independent, which means that the underwriting date of a claim has no effect on the reporting delay.
 We are not going to impose such a strong  restriction.
In order to discuss the  independence assumption, consider Figure \ref{Fig:devfactors}.
The points in the plots are derived by
first transforming the data into a  triangle with dimension $3654\times3654$,
\[
\mathcal N_{t,z}=\sum^n_{i=1}  I\big(T_i=t,\  Z_i =z\big), \quad (t,z)\in \{1,2,\dots,3654\}^2, t+z \leq 3654
\]
and
then aggregating the data into a quarterly triangle, ($\mathcal N^Q_{t,z})$, with dimension $40\times 40$, see also \citet{Hiabu:17}.
Then, for $t=2,\dots,5$, one derives the quarterly hazard rate as ratio of occurrence and exposure, 
$\overline \alpha(t,z)=\mathcal N^Q_{t,z}/{\sum_{l=1}^{z}\mathcal N^Q_{t,l}}$.
These values are then scaled by an eye-picked  norming  factor, $\overline \alpha_0(t), t=2,\dots,5$, letting  $\overline \alpha(t,z)$ start at around $1$ as a function of $z$ with fixed $t$.
The final values, $\overline \alpha_1(t,z)= \overline \alpha_0(t) \overline \alpha(t,z)$, are displayed in Figure \ref{Fig:devfactors}.
We only show plots for $t\leq5$ since almost all claims are reported after five   quarters.

If the independence assumption of \citet{Hiabu:etal:16} is satisfied,  the points should lie around a horizontal line in each plot.
If the multiplicative hazard assumption of this paper is satisfied, then any smooth shape is allowed, but all four graphs must be equal after correction for noise. This is because under the model which will be defined below, the graphs, $\overline \alpha_1(\cdot, \cdot)$, with the first component fixed, mimic a quarterly version of $\alpha_1$.
 
Inspecting the four plots, one can argue to see a negative drift of similar magnitude in each graph,  the values decaying from around $1$ to $0.8$.  
This  indicates that the approach of this paper should give a better fit to the data compared to the model of \cite{Hiabu:etal:16}.
 
%A heuristic idea of the multiplicative assumption is that a correction factor less than 1 for one underwriting year means that claims develop faster with growing development time relative to the baseline hazard. 
%Accordingly slower for a correction factor greater than 1

%To estimate the outstanding loss liabilities we  would like to
%consider two scenarios. 
%First we would like to estimate 
%the conditional hazard given the underwriting year
%$\alpha_y(w)$.
%Secondly we would like to compute an estimate of the hazard conditional only on the calendar time that is
%$\alpha_{(y+w)}(w)$.
%Thus, in the second scenario we assume that only the calendar time
%as part of the information given by the underwriting year has impact on the development time.
%Therefore we will th0
From this discussion we now continue with embedding our observations in the proportional hazard framework. Afterwards we will show how the hazard estimate can be used to forecast the number of outstanding claims. First note that we cannot apply the approach of this paper  directly, since in this application
we only observe $T$ if $T\leq R_0-Z$, which is a right  truncation. Analogue to  \citet{Hiabu:etal:16}, we transform the random variable $T$ to 
$  T^R=R_0-T$. This has the result that the right truncation truncation  becomes left truncation, $T^R \geq Z$.
Thus we consider the random variable $T^R$ as our variable of interest.
With the notation considered in Section \ref{sec:model},  we now have
$T=T^R, d=1, Z(t)=Z, \delta=1$,
$\mathcal I= \{ (t^R,z)\in \mathcal R | 0\leq z \leq t^R\}$. 
We conclude that the counting process 
$
N_i(t^R)=\mathds 1 \left\{T^R_i \leq t^R\right\},
$
satisfies Aalen's multiplicative intensity model with respect to the filtration given in Section \ref{sec:model} and
\begin{align*}
 \alpha_z(t^R)&=\alpha(t^R,z) =\lim_{h \downarrow 0} h^{-1}\mathrm{Pr}\{  T^R\in [t^R, t^R+h)| \   T^R\geq t^R, \ Z=z\}, \\
Y_i(t^R)&=  \mathds 1 \big\{(t^R,Z_i)\in \mathcal I, \ t^R\leq T_i^{R,*}\big\}.
\end{align*}

In Section \ref{first step} we suggested a local constant estimator as pilot. In this application we prefer a local linear estimator because we anticipate high mass at the boundaries.
The local linear estimator \citep{Nielsen:98} has an automated boundary correction.
We  estimate  the unstructured hazard as ratios of occurrence and exposure from a local linear estimation  \citep{Gamiz:etal:13a}: 
\begin{align*}
\widehat O^{LL}(x)&=n^{-1}\sum_{i=1}^n \int \left\{1-(x-X_i(s))D(x)^{-1}c_1(x)\right\} K_b(x-X_i(s)) \mathrm d N_i(s), \\
\widehat E^{LL}(x)&=n^{-1}\sum_{i=1}^n \int \left\{1-(x-X_i(s))D(x)^{-1}c_1(x)\right\} K_b(x-X_i(s)) Y_i(s) \mathrm ds,
\end{align*}
where the components of the $(d+1)$ - dimensional vector $c_1$ are 
\begin{align*}
c_{1j}(x)&=n^{-1}\sum_{i=1}^n\int  K_b(x-X_i(s))(x_j-X_{ij}(s))Y_i(s)ds, \quad j=0,\dots,d,
\end{align*}
and the entries $(d_{jk})$ of the $(d+1)\times (d+1)$ - dimensional matrix $D(x)$ are given by
 \begin{align*}
d_{jk}(x)&=n^{-1}\sum_{i=1}^n \int K_b(x-X_i(s))(x_j-X_{ij}(s))(x_k-X_{ik}(s))Y_i(s) \mathrm ds.
\end{align*}
Note that we have $\mathcal X=\mathcal I$.
The components of the multiplicative conditional hazard are then computed as in \eqref{iteration}.
These estimators require the choice of the bandwidth parameter, which was assumed to be scalar in order to simplify the notation in this paper. In this application we generalize this restriction allowing for different smoothing levels in each dimension, namely reporting delay and underwriting time. The bandwidth parameter is then a vector $b=(b_0,b_1)$ and we estimate it using cross-validation, see further details in the Supplementary Material. To alleviate the computational burden of cross-validation we aggregated the data triangle $\mathcal N_{t,z}$ considering bins of two days when applying a discrete version of the estimators described in the Supplementary Material. After several trials we run the cross-validation minimization over $b_0\in\{1300,1400,1500,1600,1700, 1800\}$ and $b_1\in\{2,3,4,5\}$. The cross-validated bandwidth components were $b_0=1600$ and $b_1=3$ (unit=2 days).

\begin{figure}
\centering
\includegraphics[width=14cm]{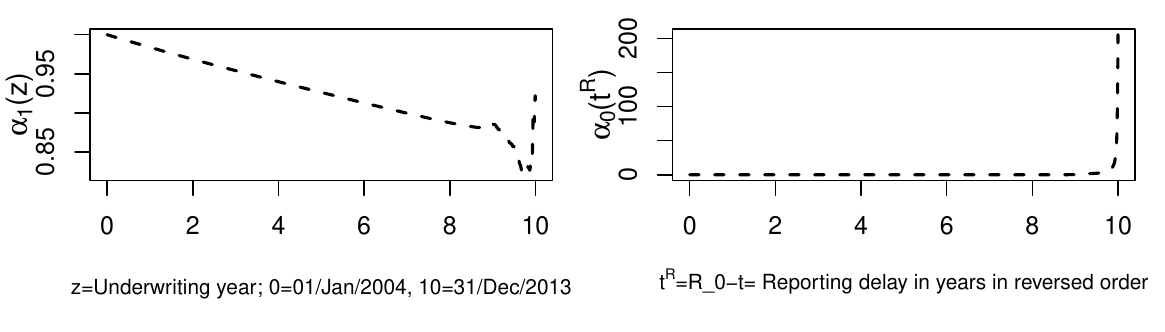} 
\caption{\label{alpha12} The estimated multiplicative hazard components . }
\end{figure}

The result of the estimation procedure is given in Figure \ref{alpha12} which dispays the estimated components of the multiplicatively structured  hazard estimator.
The function $\alpha_2$ which captures the underwriting date effect,
seems, up to a possible boundary effect, linear. This suggests that a semiparametric approach
with nonparametric baseline hazard and linear covariate effect might be suitable. While we do not investigate  this point further in our data illustration, this  particular case illustrates nicely how   our nonparametric approach could be employed for model selection.
%Figure \ref{difference} shows the difference, $\widetilde \alpha(x)-\widehat \alpha(x)$,
%between the structured and the unstructured estimators.

%
%\begin{figure}[!ht]
%\centering
%\includegraphics[width=12cm]{Figure_diff} 
%\caption{\label{difference} Difference between structured and unstructured hazard estimator, $\widetilde \alpha(x)-\widehat \alpha(x)$, on a 2-day grid. }
%\end{figure}
%

Finally the total number of outstanding claims, i.e. the reserve, can be estimated as
%\begin{align}
%\int_0^T R(y) \mathrm dy, \\
%R(y)= \frac{ \int_{T-y}^T \widehat f_y(t) \mathrm dt}{\int_{0}^{T-y} \widehat f_y(t) \mathrm dt} \sum_{i=1}^n K_b( y- Y_i),
%\end{align}
%or
\[
\textrm{Reserve}=\sum_{i=1}^n \frac{ \int_{R_0-Z_i}^{R_0} \widehat f_{Z_i}(t) \mathrm dt}{\int_{0}^{R_0-Z_i} \widehat f_{Z_i}(t) \mathrm dt},
\quad 
\widehat f_z(t)= \widehat \alpha_0(R_0-t) \widehat\alpha_1(z)\exp\left\{-\int_0^{R_0-t}\widehat \alpha_0(s) \widehat\alpha_1(z)\mathrm ds\right\}.
\]
Note that $\widehat f_z(t)$ is an estimator of the conditional density of the survival time $T$.
The reserve can be also decomposed further to provide the 'cash flow' of the next periods.
If the future is divided into $M$ periods, each of them with length $\delta=R_0/M$, then
%one can divide the square $\mathcal R$ into $M$ equally sized bins
%$[t_k,t_{k+1}]\times[z_l,z_{l+1}] (k,l=0,\dots, M)$.
the amount of claims forthcoming in the $a$th $(a=1,2,\dots, M$) period can be  estimated by
\[ \textrm{Reserve}_P(a)=\sum_{i=1}^n 
\frac{ \int_{(R_0-Z_i+a\delta-1)\wedge R_0}^{(R_0-Z_i+a\delta)\wedge R_0} \widehat f_{Z_i}(t) \mathrm dt}
{\int_{0}^{R_0-Z_i} \widehat f_{Z_i}(t) \mathrm dt}.
\]

{\small 
\begin{table}
\caption{\label{tab:reserves} Number of outstanding claims for future quarters; $1=2014 \ Q1, \ldots, \ 39=2022 \ Q3$. The  backfitting approach in this paper (SBF) is compared with the chain ladder method (CLM) and the approach in \citet{Hiabu:etal:16}.}
\centering
\fbox{%
\begin{tabular}{rccccccccccccc} %\hline
Future quarter &1  & 2&3&4&5 &6&7 &8&9&10&11&12 -- 39&Tot. \\  \hline
%CV& 1027&733&465&201&15&5&3&2&1&1&1&0&2452\\
\cite{Hiabu:etal:16}&970&684&422&166&14&5&3&2&1&1&1&0&2270\\
CLM  &948&651&387&148&12&5&3&2&1&1&1&0&2160\\
SBF&872&621&400&130&53&7&4&3&2&1&1&1&2193
\end{tabular}
}
\end{table}
}

Table \ref{tab:reserves} shows the estimated  number of of outstanding claims for future quarters.
We compare the approach of this paper with the results derived by \citet{Hiabu:etal:16} and the traditional chain ladder method. 
The two latter approaches have in common that they assume independence between underwriting date, $Z$, and reporting delay, $T$.
We see that while all approaches estimate a similar total number of outstanding claims (reserve), those two approaches have  distributions over the quarters that are very different from the results obtained by the method proposed in this paper.
It seems that the  violation of the independence assumption has not a big influence on the reserve, since it balances 
the different development patterns arising from different periods out.
However, the problem becomes quite serious if one is interested in more detailed estimates like the cash flow.

%%%%%%%%%%%%%%%%%%%%%%%%%%%%%%%%%%%%%%%%%%%%%%%%

%%%%%%%%%%%%%%%%%%%%%%%%%%%%%%%%%%%%%%%%%%%%%%%%%%%%%%%%%%%%%%%%%%%%%%%%%%%
%\section{Conclusion}\label{sec:conclusion}
%This paper provided a first introduction of smooth backfitting into survival analysis and hazard estimation. The starting point has been the popular proportional hazard model with fully nonparametric components. One could imagine that smooth backfitting could play a role in a long list of structured problems in semiparametric and nonparametric survival analysis. One could for example imagine that smooth backfitting can provide useful extensions of some of the practical dynamic survival models of \citet{Martinussen:Scheike:06} or  and one can also think of applications to many of known extensions of the Cox regression model, see for example  \citet{Su:etal:16, Su:etal:12}, where the understanding of the link between data and estimators is improved via the direct projection approach of smooth backfitting.

%

%The research of the second author was supported by the  Deutsche Forschungsgemeinschaft through the Research Training Group RTG 1953.
%Third author acknowledges the support  from the
%Spanish Ministry of Economy and Competitiveness, through grant number MTM2013-
%41383P, which includes support from the European Regional Development Fund (ERDF). 

\newpage
\section*{Appendix}
\appendix
\section{Asymptotic properties of the smooth backfitter of multiplicative hazards}
\label{sec:prop}
The estimator $\widehat \alpha_j$ is defined as solution of a nonlinear operator equation. We are going to approximate this equation by a linear equation that can be interpreted as equation that arises in nonparametric additive regression models, and then show that the solution of the linear equation approximates $\widehat \alpha_j$. The linear equation and its solution is well understood from the theory of additive models. This will be our essential step to arrive at an asymptotic understanding of our estimator $\widehat \alpha_j$. Assumptions  [B1]-[B7]  below are of standard nature in marker dependent hazard papers and can be verified for the local constant  estimators we are interested, see in particular \citet{Nielsen:Linton:95}, \citet{Nielsen:98} and \citet{Linton:etal:03} for related calculations. However, one should notice that the conditions are not restricted to the local constant smoothers. They are not even tight to kernel smoothers. Any smoother could be used as long as it obeys the structure of being a ratio of a smoothed occurrence and a smoothed exposure.

For our main theorem we make the following assumptions.
We hereby do not make assumptions on the full support  $\mathcal R$ but only on a subset $\mathcal X \subseteq \mathcal R$. We will make use of the following notations: $\mathcal X_{x_k}$  denotes the set \\ $\{(x_0, \dots, x_{k-1},x_{k+1},\dots, x_d)|\ (x_0,\dots,x_d) \in \mathcal X\}$, and\\  $x_{-k}=(x_0, \dots, x_{k-1},x_{k+1},\dots, x_d)$  and $\mathcal X_{x_j,x_k}$  denotes the set $\{(x_l: l\in\{0,...,d\} \backslash \{j,k\}) |\ (x_0,\dots,x_d) \in \mathcal X\}$. Furthermore, we define $\mathcal X_{k}=\{x_k|\ (x_0,\dots,x_d) \in \mathcal X\text { for some values of } (x_l :l \not =k)\}$, $\mathcal X_{j,k}=\{(x_j,x_k)|\ (x_0,\dots,x_d) \in \mathcal X\text { for some values of } (x_l :l \not \in \{j,k\})\}$, and $\mathcal X_{j(x_k)}=\{x_j|\ (x_0,\dots,x_d) \in \mathcal X$ for some values of $ (x_l :l \not \in \{j,k\})\}$.
\begin{description}
  \item[A1]  The function $E(x)$ is two times continuously differentiable  for $x\in \mathcal X$ and $\inf_{x\in \mathcal X} E(x) >0$.
\item[A2] The hazard $\alpha$ is two times continuously differentiable for $x\in \mathcal X$ and $\inf_{x\in \mathcal X} \alpha(x) >0$.
\item[A3] The kernel $K$ has compact support which is without loss of generality supposed to be $[-1,1]$. Furthermore it is symmetric and continuous.
\item [A4] It holds that $nb^{5}\rightarrow c_b$ for a constant $0 < c_b < \infty$ as $n\rightarrow \infty$.
\item [A5] It holds that $$\int_{\mathcal X_{j,k}} \frac {1} {O_j(x_j)O_{k}(x_k)} \mathrm{d} x_j\ \mathrm{d} x_k < \infty$$
for $j,k= 0,...,d$, $j \not = k$, where $O_j(x_j)= \int O(x) \ dx_{-j}$ and $O(x) =\alpha^* \prod_{j=0}^d \alpha_j(x_j) E(x)$.
\item [A6] It holds that the two-dimensional marginal occurences $O_{j,k}(x_j,x_k)= \int O(x) \ dx_{-(j,k)}$
are bounded from above and bounded away from 0. 
% \item[A7]
% The projections of  $\mathcal X$ and $\mathcal R=[0,R_0]\times \prod_{i=1}^d [0,R_i]$ to their $j$th $(j=0,\dots,d)$ coordinate are equal, that is 
%\[ \bigcup_{x_j\in [0, R_j]} \mathcal X_{x_j}= [0,R_0]\times \prod _{k \neq j}[0,R_k], \quad  j=0,\dots,d.
%\]
 \item[A7] For some $\delta > 0$ it holds that for $j,k= 0,...,d$, $j \not = k$
 \begin{align*}
 &\int_{\mathcal X_{j,k} }\frac {1} {O_j^{1+\delta}(x_j)O_{k}(x_k)} \mathrm{d} x_j\ \mathrm{d} x_k < \infty,\\
& \sup_{x_k \in \mathcal X_{k}}  \int_{\mathcal X_{j(x_k)} }\frac {1} {O_j^{1-\delta}(x_j)O_{k}(x_k)} \mathrm{d} x_j < \infty,\\
& \sup_{x_k \in \mathcal X_{k}}  \int_{\mathcal X_{j(x_k)} }\frac {1} {O_j^{1/2}(x_j)O_{k}^{1/2}(x_k)} \mathrm{d} x_j < \infty.
 \end{align*} 
 \end{description}
Note that  assumptions [A1]-[A4] are standard in kernel smoothing theory. In Assumptions [A5] and [A6] we only assume that the two-dimensional marginal occurrences of $O$ are bounded from above and bounded away from 0, but we do not make the assumption that the one-dimensional marginal occurrences have this property. 
This allows that the support of a two-dimensional marginal density $O_{jk}$ has a triangular shape $\{ (x_j,x_k): x_j + x_k \leq c;\  x_j, x_k \geq 0\}$ for some constant $c>0$. This can be easily seen. Suppose for simplicity that $O_{jk}$ is the uniform density on the triangle. Then $O_j(x_j) = 2 c^{-2}(c-x_j)_+$ and  $O_{k}(x_k) = 2 c^{-2}(c-x_k)_+$ and we have  $$\int \frac {1} {O_j(x_j)O_{k}(x_k)} \mathrm{d} x_j\ \mathrm{d} x_k = \int_{x_j + x_k \leq c;\ x_j,x_k \geq 0}\frac 2 {c^2} \frac {1} {(c-x_j)(c-x_k)} \mathrm{d} x_j \ \mathrm{d} x_k = \frac{\pi^2}{3c^2}<\infty.$$
Thus, our assumption [A5] on one-dimensional marginals is fulfilled. One can easily verify that also [A7] holds for this example.  This discussion can be extended to other shapes of two-dimensional marginals that differ from rectangular supports. Note also that [A5] and [A7] trivially hold if the one-dimensional marginal $O_j$ are bounded away from zero.

The estimators $\widehat \alpha_0, \dots, \widehat \alpha_d$ can be rewritten  as solutions of
\begin{align*}
\int_{\mathcal X_{x_k}}\widehat O(x)\mathrm dx_{-k}
-
\int_{\mathcal X_{x_k}}  \widehat \alpha^* \prod_{j=0}^d\widehat \alpha_j(x_j)\widehat E(x)\mathrm dx_{-k}=0,  \quad k=0,\dots, d.
\end{align*}
Since,
$
\int_{\mathcal X_{x_k}} O(x)\mathrm dx_{-k}
-
\int_{\mathcal X_{x_k}}  \alpha^* \prod_{j=0}^d \alpha_j(x_j)E(x)\mathrm dx_{-k}=0,
$
the difference of those two zero-terms is zero as well, and we have
\begin{align}
0&=\widehat \Delta_k(x_k)-\int_{\mathcal X_{x_k}}\left\{\widehat \alpha^* \prod_{j=0}^d \widehat \alpha_j(x_j)-\alpha^*\prod_{j=0}^d  \alpha_j(x_j)\right\} \widehat E(x) \mathrm dx_{-k} \notag\\
&=\widehat \Delta_k(x_k)-\int_{\mathcal X_{x_k}}\left[(1+ \widehat \delta^*)\prod_{j=0}^d  \{1+\widehat \delta_j(x_j)\}-1\right]\alpha^* \prod_{j=0}^d  \alpha_j(x_j) \widehat E(x) \mathrm dx_{-k},  \label{backfitting:deltahat}
\end{align}
where
\begin{align*}
\widehat \Delta_k(x_k)&= \int_{\mathcal X_{x_k}} \left \{\widehat O(x)-O(x)\right\} \mathrm dx_{-k}- \int_{\mathcal X_{x_k}}\alpha^* \prod_{j=0}^d  \alpha_j (x_j)\{\widehat E(x)-E(x)\} \mathrm dx_{-k}, \\
\widehat \delta_j(x_j)&= \frac{\widehat \alpha_j (x_j)- \alpha_j(x_j)}{\alpha_j(x_j)},\\
\widehat \delta^*&= \frac{\widehat \alpha^*- \alpha^*}{\alpha^*}.
\end{align*}
Note that $\widehat \delta$ is defined as a root of a non-linear operator.
Motivated by \eqref{backfitting:deltahat}, we define  an approximation, $\overline \delta^*$ and  $\overline \delta_j(x_j)$ $(0 \leq j \leq d)$, as solution of the linear equation
\begin{align}\label{deltabar}
\int_{\mathcal X_{x_k}}\left [\overline \delta^*+\sum_{j=0}^d  \overline \delta_j(x_j)\right ] \alpha^*\prod_{j=0}^d  \alpha_j(x_j) \widehat E(x) \mathrm dx_{-k}= \widehat \Delta_k(x_k)
\end{align}
under the constraint $\int \overline \delta_k(x_k)  \omega_k(x_k) \mathrm{d} x_k =0$ where 
$$\omega_k(x_k) = \int\prod_{j=0}^d  \alpha_j(x_j)\widehat E(x) \mathrm dx_{-k}.$$
This is equal  to the constraint \eqref{qprop1} in the main paper for the choice $w_k(x_k)= \int\prod_{j=0}^d  \alpha_j(x_j)\widehat E(x) \mathrm dx_{-k} /(\alpha_k(x_k)$ $ \int\prod_{j=0}^d  \alpha_j(x_j)\widehat E(x) \mathrm dx)$. Under this constraint  one has that $$\overline \delta^* = \frac {\int_{\mathcal X} \left \{\widehat O(x)-O(x)\right\} \mathrm dx- \int_{\mathcal X}\alpha^*\prod_{j=0}^d  \alpha_j (x_j)\{\widehat E(x)-E(x)\} \mathrm dx}{\int_{\mathcal X}\alpha^*\prod_{j=0}^d  \alpha_j (x_j)\widehat E(x) \mathrm dx}.$$
Note that the  norming of the constraint  cannot be used in practice because $\alpha$ is unknown but it will simplify the theoretical discussion and the results can be carried over to feasible weighting.

Equation \eqref{deltabar} can be rewritten as an integral equation of the second kind
%\begin{align*}
%&\overline \delta_k(x_k)\int_{\mathcal X_{x_k}}\prod_{j=0}^d  \alpha_j(x_j) \widehat E(x) \mathrm dx_{-k}\\  &+
% \sum_{j\neq k}\int_{\mathcal X_{x_k}} \left\{\int{\mathcal X_{x_k,x_j}}\prod_{j=0}^d   \alpha_j(x_j) \widehat E(x)\mathrm d_{-(k,j)}\right\} \overline \delta_j(x_j)\mathrm dx_{-j}=\widehat \Delta_k(x_k)
%\end{align*}
%or
\begin{align*}
\overline \delta_k(x_k)+ \sum_{j\neq k}  \int_{\mathcal X_{j(x_k)}}  \widehat \pi_{k,j}(x_k,x_j) \overline \delta_j(x_j) \mathrm dx_{j}=\widehat \mu_k(x_k) - \overline \delta^*,
\end{align*}
with
\begin{align*}
&\widetilde O (x) =\alpha^* \prod_{j=0}^d \alpha_j(x_j) \hat E(x), \\
 &\widetilde O _{j,k}(x_j,x_k)= \int \widetilde O (x) \ dx_{-(j,k)},\\
 &\widetilde O _{k}(x_k)= \int \widetilde O (x) \ dx_{-k},\\
&\widehat \pi_{k,j}(x_k,x_j)=\frac{\widetilde O _{j,k}(x_j,x_k)}{\widetilde O _{k}(x_k)}, \\
&
\widehat \mu_k(x_k)= \frac{\widehat \Delta_k(x_k)}{\widetilde O _{k}(x_k)}.
\end{align*}
Note that all these functions depend on $n$. The integral equation can also be simply written as $\overline \delta + \widehat \pi \overline \delta=\widehat \mu -\overline \delta^*$ with $\overline \delta = (\overline \delta_0,...,\overline \delta_d)^\intercal$, where $\widehat \pi$ is the integral operator with kernel $\widehat \pi_{k,j}$. We will show that $\overline \delta$ approximates $\widehat\delta$. Before we come to this point we state a proposition that gives the asymptotics for $\overline \delta$ .

For the next results we need some conditions on the estimators $\widehat E$ and $\widehat O$. We decompose $\widehat \mu_k$  into three terms, $\widehat \mu_k=\widehat \mu^A_k+\widehat \mu ^B_k+\widehat \mu ^C_k$, which depend on $n$. For some deterministic functions $O^*(x)$ and $E^*(x)$ these terms are defined as:
\begin{align*}
\widehat \mu^A_k(x_k)&= \frac {-\int_{\mathcal X_{x_k}}\alpha^* \prod_{j=0}^d  \alpha_j (x_j) \left \{\widehat E(x)-E^*(x)\right\} \mathrm dx_{-k}+\int_{\mathcal X_{x_k}} \left \{\widehat O(x)-O^*(x)\right\} \mathrm dx_{-k}}   {\widetilde O _{k}(x_k)},\\
\widehat  \mu^B_k(x_k)&=  \frac {-\int_{\mathcal X_{x_k}}\alpha^* \prod_{j=0}^d  \alpha_j (x_j)\left \{ E^*(x)-E(x)\right\} \mathrm dx_{-k}+\int_{\mathcal X_{x_k}} \left \{ O^*(x)-O(x)\right\} \mathrm dx_{-k}}{ O_{k}(x_k)},\\
\widehat  \mu^C_k(x_k)&= \left [ \frac { O_{k}(x_k)} {\widetilde O _{k}(x_k)} -1 \right ]\widehat  \mu^B_k(x_k). \end{align*} Note that $\widehat \mu^B_k$ are deterministic functions.
Typical choices of $O^*(x)$ and $E^*(x)$ are the expectations of $\widehat E(x)$ and $\widehat O(x)$. Then $\widehat \mu^A_k(x_k)$ is the stochastic part of a one-dimensional smoother and $\widehat \mu^B_k(x_k)$ is its bias. Both terms are well understood and can be easily treated by standard smoothing theory. We now want to develop an asymptotic theory for the estimators $\hat \alpha_j$ where their asymptotic properties are described by the properties of $\widehat \mu^A_k(x_k)$ and $\widehat \mu^B_k(x_k)$. We use the following normings of these quantities: $\overline \mu_k^{*,r} = \int \widehat \mu^r_k(x_k)  O _{k}(x_k)\ \mathrm{d} x_k$ for $r \in\{A,B,C\}$ and $\overline \mu_k^{*,r,n} = \int \widehat \mu^r_k(x_k) \widetilde O _{k}(x_k)\ \mathrm{d} x_k$ for $r \in\{A,B,C\}$. 

We assume that with \begin{align*}
 \pi_{k,j}(x_k,x_j)&=\frac{\int\prod_{l=0}^d  \alpha_l(x_l) E(x) \mathrm dx_{-(k,j)}}{\int\prod_{l=0}^d  \alpha_l(x_l)E(x) \mathrm dx_{-k}}
\end{align*} the following assumptions hold:

\begin{description}
  \item[B1]  It holds that $\int \widetilde O  (x)^2 dx = O_P(1)$ and $$\widetilde O _{j,k}(x_j,x_k) -  O_{j,k}(x_j,x_k) = o_P((\log n)^{-1})$$ uniformly over $0\leq j < k \leq d$ and $x_j,x_k$,
 where $ O _{j,k}(x_j,x_k)= \int  O (x) \ dx_{-(j,k)}$. 
  \item[B2]  $$\sup_{x_j} | O_j^{1/2}(x_j)  \widehat \mu_j^A(x_j)| = o_P(n^{-1/5})$$ 
  and $$\sup_{x_j} | O_j^{1/2}(x_j)  \widehat \mu_j^B(x_j)| = o_P(n^{-1/5})$$
  for $0\leq j  \leq d$, where $ O _{j}(x_j)= \int  O (x) \ dx_{-j}
$.
   \item[B3]  For $x_j$ with $O_j(x_j) > 0$ it holds that $$ n^{2/5} \widehat \mu_j^A(x_j) \to \mbox{N}(0, \sigma_j^2(x_j))$$ for $0\leq j  \leq d$ with some function $\sigma_j^2(x_j) >0$.
\item[B4] $$\int \widehat \mu_j^A(x_j)^2 O_j(x_j) \ \mathrm{d} x_j = O_P(n^{-4/5})$$ and 
$$\int \widehat \mu_j^B(x_j)^2 O_j(x_j) \ \mathrm{d} x_j = O(n^{-4/5})$$
for $0\leq j  \leq d$.
\item[B5] It holds that 
$$\sup_{x_j \in \mathcal X_{j}}  O_j^{1/2}(x_j)\int_{\mathcal X_{k(x_j)} }\frac {O_{j,k}(x_j,x_k)} {O_j(x_j)}  \widehat \mu_k^A(x_k) \mathrm dx_k = o_P(n^{-2/5}(\log n) ^{-1}).$$

\item[B6] It holds for $0\leq j  \leq d$ that 
\begin{align*}
&\sup_{x_j \in \mathcal X_{j}}  \int_{x \in \mathcal {X} }\frac {1} {\prod_{k \in \{0,...,d\}\backslash \{j\}} O_{k}(x_k)^{1/2}} {\widetilde O(x)}  \mathrm dx_{-j} = O_P(1),\\
&\sup_{x_j \in \mathcal X_{j}}  \int_{x \in \mathcal {X} }\frac {1} {\prod_{k \in \{0,...,d\}\backslash \{j\}} O_{k}(x_k)^{1/2}} { O(x)}  \mathrm dx_{-j} = O(1).\end{align*}

\item[B7] It holds that $\overline \mu^{*,B,n}- \overline \mu^{*,B}= o_p(n^{-2/5})$.

  \end{description}

We shortly discuss these assumptions. Condition [B1] is a mild consistency assumption for a two-dimensional smoother. Also [B2] is a weak condition because in our setting one-dimensional smoothers are typically $O_P(n^{-2/5})$-consistent. [B3] is a standard limit result for many one-dimensional smoothers and [B4] assumes  rates for L$_2$-norms of the stochastic part and the bias part of a nonparametric  one-dimensional smoother that are standard under our smoothness assumptions. For the interpretation of [B5] note that the integral on the left hand side in the formula is a global average of $\widehat \mu_k^A$. Because $\widehat \mu_k^A$ is a local average this is a global weighted average of mean zero random varibales. Thus one expects a $O_P(n^{-1/2})$ rate for the integral. For the supremum of the integrals one expects a $O_P((\log n)^{1/2} n^{-1/2})$ rate which is faster than the required rate. [B6] states a bound for the total number of occurrences. It can be easily verified under the assumption that the one dimensional marginals $O_{l}$ are bounded from below. Furthermore, one can easily check it if the one-dimensional marginals $O_{l}$ are bounded from below for $l \not \in \{j,k\}$ and if the two-dimensional marginal $O_{j,k}$ has the properties discussed in the example after assumption [A7]. The following proposition states a stochastic expansion for $\overline \delta$.
 
\begin{proposition}\label{prop:deltabar}
Make the assumptions [A1]--[A7], [B1]--[B7]. Then the function $\overline \delta =(\overline \delta_0,..., \overline \delta_d)$, introduced in \eqref{deltabar}, exists and is uniquely defined, with probability tending to one.
Moreover, it has the following expansion:
\[
\norm{\overline \delta - \widehat  \mu^A- (I- \pi)^{-1}(\widehat \mu^B- \overline \mu^{*,B})}_{O,\infty} = o_p(n^{-2/5}),
\]
where, for a function $f(x) = (f_0(x_0),...,f_d(x_d))^\intercal$, we define
$$\| f \|_{O,\infty} = \sup_{x\in {\cal X}} \max_{0 \leq j \leq d} |O^{1/2}_j(x_j) f_j(x_j)|.$$
Furthermore, the function $\pi: {\cal L} \to {\cal L}$ is defined as 
$\pi_k(f)(x_k) = \sum_{j \not = k} \int \pi_{kj}( x_k, x_j) f_j(x_j) \mathrm{d} x_j$ for $f \in {\cal L}= {\cal L}_1 \times ... \times {\cal L}_d$ with $ {\cal L}_j =\{ \delta_j: {\mathcal X}_j \to \mathbb{R}: \int_{{\cal X}_j}  \delta_j^2(x_j)O_j(x_j) \ \mathrm{d} x_j  < \infty,  \int_{{\cal X}_j}  \delta_j(x_j)O_j(x_j) $ $ \mathrm{d} x_j =0\}$.

\end{proposition}

From the proposition we get as a corollary the asymptotic distribution of $\overline \delta_j(x_j)$. 

\begin{proposition}\label{prop:expansion}
Make the assumptions [A1]--[A7], [B1]--[B7]. Then for $x_j$ ($0\leq j  \leq d$) with $O_j(x_j) > 0$ it holds that 
\begin{align*}
 n^{2/5}\{\overline \delta_j(x_j) -[ (I- \pi)^{-1}(\widehat \mu^B - \overline \mu^{*,B})]_j(x_j) \}\rightarrow \mbox{N}(0, \sigma_j^2(x_j)),\end{align*}
 in distribution. Under the additional assumption $\widehat \mu_j^B(x_j)= O(n^{-2/5})$ we have that the bias $[ (I- \pi)^{-1}(\widehat \mu^B - \overline \delta^{B,*})]_j(x_j)$ is of order $ O(n^{-2/5})$.
\end{proposition}

The following theorem states that  $\overline \delta$  is indeed a good approximation of the relative estimation error $\widehat \delta$.
\begin{theorem}\label{theoremfbar}
Under assumptions [A1]--[A7], [B1]--[B7] it holds that with probability tending to one  there exists a solution $\widehat \delta^*$ and $\widehat \delta=( \widehat \delta_0,\dots, \widehat \delta_d)$ that solves equation \eqref{backfitting:deltahat} with
\begin{align*}
&\norm {\widehat \delta -\overline \delta}_{O,\infty}=o_p(n^{-2/5}),\\
&\widehat \delta^* -\overline \delta^*=o_p(n^{-2/5}).
\end{align*}
For this solution we get that
\begin{align*}
 n^{2/5}\{(\widehat \alpha_j- \alpha_j)(x_j) -  \alpha_j(x_j) [ (I- \pi)^{-1}(\widehat \mu^B - \overline \mu^{*,B})]_j(x_j)\} \rightarrow \mbox{N}(0, \alpha^2_j(x_j)\sigma_j^2(x_j)),\end{align*}
 in distribution, for $x_j$ ($0\leq j  \leq d$) with $O_j(x_j) > 0$.
\end{theorem}

\section{Proofs}\label{sec:proof}

\subsection{Proof of Proposition \ref{prop:deltabar}} The proof of this proposition follows the lines of the proof of Theorem 1 in \citet{Mammen:etal:99} but it needs some major modifications in the last steps of the proof because we have weaker assumptions than the ones assumed in the latter theorem. We outline that the first part of the proof in \citet{Mammen:etal:99} also goes through under our weaker assumptions and we show how additional arguments can be used in the last part.

Note that under our assumptions [A5], [A6] we get that 
$\int O_{jk}(x_j,x_k)^2 O_{j}(x_j)^{-1} O_{k}(x_k)^{-1}$ $ \mathrm{d} x_j \ \mathrm{d} x_k$ $ < \infty$. 
As in Lemma 1  in \citet{Mammen:etal:99} this implies that for some constants $c,C> 0$
\begin{equation} \label{normbound} c \max_{0 \leq j \leq d} \| \delta_j\| \leq  \| \delta_0 + ... + \delta_d\| \leq C \max_{0 \leq j \leq d} \| \delta_j\|\end{equation}
for $\delta_j \in {\cal L}_j$ where $\|...\|$ denotes the norm $\| m(x)\|^2 = \int m(x)^2 O(x)\ dx$.  Furthermore, one gets that 
$\|T\|= \sup \{ \|T(\delta_0+...+\delta_d)\|:  \delta_j \in {\cal L}_j$ with $ \|\delta_0+...+\delta_d\| < 1 \}< 1$, where here $T$ is the operator $T= \Psi_d\cdot ... \cdot \Psi_0$ with
\begin{align*}\Psi_j(\delta^* + \delta_0+...+\delta_d)(x)&= 
 \delta^* +\delta_0(x_0) +...+\delta_{j-1}(x_{j-1}) \\ & \qquad + \delta_j^*(x_j) + \delta_{j+1}(x_{j+1}) + ...+ \delta_{d}(x_{d}) ,\\
\delta_j^* (x_j) &= - \sum_{k \not = j} \int \delta_k(x_k) \pi_{j,k}(x_j,x_k)\ \mathrm{d} x_k \end{align*}
for $(\delta^*, \delta_0,...,\delta_d) \in  {\cal L}$.
Furthermore, note that 
for $j\neq k$ it holds that 
\begin{align*}
&\norm{\frac{\widetilde  O_j(x_j)-O_j(x_j)}{O_j(x_j)}}=o_P(1), \\
&\int {\left (\frac{\widetilde O_{j,k}(x_j,x_k)} {O_j(x_j)O_{k}(x_k)}- \frac{O_{j,k}(x_j,x_k)} {O_j(x_j)O_{k}(x_k)}\right )^2 O_j(x_j)O_{k}(x_k) \mathrm{d} x_j \mathrm{d} x_k}=o_P(1), \\
&\int {\left (\frac{\widetilde O_{j,k}(x_j,x_k)} {\widetilde O_{j}(x_j)O_{k}(x_k)}- \frac{O_{j,k}(x_j,x_k)} {O_j(x_j)O_{k}(x_k)}\right )^2 O_j(x_j)O_{k}(x_k) \mathrm{d} x_j \mathrm{d} x_k}=o_P(1).
\end{align*}
These equations follow from [A5], [A6] and [B1]. To see this note that [A6] and [B1] imply that, uniformly for $x_j,x_k$ it holds that $\widetilde O_{j,k}(x_j,x_k) - O_{j,k}(x_j,x_k)= o_P((\log n)^{-1}) O_{j,k}(x_j,x_k)$. This gives that 
\begin{equation} \label{equniform} [\widetilde O_{j}(x_j) /O_j(x_j)] - 1 = o_P((\log n)^{-1}), \end{equation} uniformly for $x_j \in {\mathcal X}_j$ and $0 \leq j \leq d$. Together with  [A5] and [B1], this implies the three equations. As in Lemma 2 in \citet{Mammen:etal:99} we conclude from these equations that $$\| \hat T \|_n < \gamma$$
for some $\gamma < 1$ with probability tending to one. Here, we define $\hat T$, $\| ...\|_n$, ${\mathcal X}_{n,j}$, $\widehat \Psi_j$, $\widetilde {\cal L}_j$ and $\widetilde {\cal L}$ as $T$, $\| ...\|$, ${\mathcal X}_{j}$, $\Psi_j$, ${\cal L}_j$ and $ {\cal L}$  but with $O_j, \pi_{jk}$ replaced by $\widetilde O_{j},\widehat \pi_{jk}$ $(0 \leq j,k\leq d; j \not = k)$.
In particular, we put 
$\widetilde  {\cal L}_j =\{ \delta_j: {\mathcal X}_j \to \mathbb{R}: \int_{{\cal X}_j}  \delta_j^2(x_j)\widetilde  O_j(x_j) \ \mathrm{d} x_j  < \infty,  \int_{{\cal X}_j}  \delta_j(x_j)\widetilde  O_j(x_j) \ \mathrm{d} x_j =0\}$,
$\widetilde{\cal L}= \widetilde {\cal L}_1 \times ... \times \widetilde {\cal L}_d$, and
$\|T\|_n= \sup \{ \|T(\delta_0+...+\delta_d)\|_n:  \delta_j \in \widetilde {\cal L}_j$ with $ \|\delta_0+...+\delta_d\|_n < 1 \}$.

Arguing as in the first part of Lemma 3 in \citet{Mammen:etal:99} this gives that 
$\overline \delta_k(x) = \overline \delta_k^A(x)+ \overline \delta_k^B(x) + \overline \delta_k^C(x)$,
where  for $r\in \{A,B,C\}$ the functions $\overline \delta^r_k \in \widetilde {\cal L}_k$ are defined by
$$\overline \delta_0^r(x_0)+...+\overline \delta_d^r(x_d)= \sum_{l=0}^s \widehat T^l \widehat \tau^r (x) + \widehat R^{r,[s]}(x) $$
with $\|\widehat R^{r,[s]}\| \leq C \gamma^s$ with probability tending to one for some constant $C>0$. Here we put
\begin{align*} \widehat \tau^r &= \widehat \Psi_d \cdot ... \cdot \widehat \Psi_1 ( \widehat \mu^r_0 - \overline \mu_0^{*,r,n}) + ... + \widehat \Psi_d  ( \widehat \mu^r_{d-1}- \overline \mu_{d-1}^{*,r,n} ) +  ( \widehat \mu^r_{d}- \overline \mu_d^{*,r,n} ), \\
\widehat R^{r,[s]}(x) &=  \sum_{l=s+1}^\infty \widehat T^l \widehat \tau^r (x). \end{align*} 
Up to this point we followed closely the arguments in the proof of Theorem 1 in \citet{Mammen:etal:99}. 
The arguments of the further parts of the proof of the latter theorem would need  that, in our notation, 
\begin {equation}\label{eq:nottrue} \sup_{x_j\in {\cal X}_j} \int_{{\cal X}_{k(x_j)}} \frac{\widetilde O^2_{j,k}(x_j,x_k)} {\widetilde O^2_j(x_j)O_k(x_k)} \mathrm d x_k \end{equation}
is bounded by a constant, with probability tending to one. This would imply that with probability tending to one for some constant $C>0$ for all functions $g: {\cal X}_{k(x_j)}\to \mathbb R$
\begin{align} \label{eq:not1}  \sup_{x_j\in {\cal X}_j}\bigg | \int_{{\cal X}_{k(x_j)}} \frac{\widetilde O_{j,k}(x_j,x_k)} {\widetilde O_{j}(x_j)} g(x_k)\mathrm d x_k\bigg | \leq C \|g\|,
\end{align}
as can be seen by application of the Cauchy-Schwarz inequality. The proof of Theorem 1 in \citet{Mammen:etal:99} shows that this can be used to show that $\sup_{x\in {\cal X}, 0 \leq j \leq d} |R_j^{r,[s]}(x) | \leq C \gamma^s$ with probability tending to one for some constant $C>0$. 
Unfortunately in our setting \eqref{eq:nottrue} does not hold and thus we cannot follow that 
\eqref{eq:not1} holds in our setting. Indeed, one can check that in general \eqref{eq:not1} does not hold under our assumptions. Consider e.g. the set-up discussed after the statement of assumption [A7]. Thus we do not have that $T$ and $\hat T$ map a function with bounded L$_2$-norm into a function with bounded L$_\infty$-norm. This also does not hold if we replace the L$_\infty$-norm by our weighted norm $\|..\|_{O,\infty}$. We now argue that after twice application of  $T$ or $\hat T$ a function with bounded $\|..\|$-norm is transformed into a function with bounded $\|..\|_{O,\infty}$-norm. This follows from the following two estimates for functions $g: \mathcal X_{k} \to \mathbb R$ with some constant $C>0$
\begin{align} \label{helpdec1}
&\int_{\mathcal X_{j} }\left (\int_{\mathcal X_{k(x_j)} }\frac {O_{j,k}(x_j,x_k)} {O_j(x_j)}g(x_k) \mathrm dx_k\right )^2 O^{1-\delta}_j(x_j)\mathrm dx_j  \leq C \int_{\mathcal X_{k} }O_{k}(x_k)g^2(x_k) \mathrm {d}x_k,\\ \label{helpdec2}
&\sup_{x_j \in \mathcal X_{j}}  O_j^{1/2}(x_j)\left |\int_{\mathcal X_{k(x_j)} }\frac {O_{j,k}(x_j,x_k)} {O_j(x_j)} g(x_k) \mathrm dx_k\right | \leq C \left (\int_{\mathcal X_{k} }O^{1-\delta}_{k}(x_k)g^2(x_k) \mathrm dx_k\right )^{1/2}.
\end{align}
Furthermore, it holds with probability tending to one, that
\begin{align} \label{helpdec1a}
&\int_{\mathcal X_{j} }\left (\int_{\mathcal X_{k(x_j)} }\frac {\widetilde O_{j,k}(x_j,x_k)} {\widetilde O_j(x_j)}g(x_k) \mathrm dx_k\right )^2 O^{1-\delta}_j(x_j)\mathrm dx_j \leq C \int_{\mathcal X_{k} }O_{k}(x_k)g^2(x_k) \mathrm {d}x_k,\\ \label{helpdec2a}
&\sup_{x_j \in \mathcal X_{j}}  O_j^{1/2}(x_j)\left |\int_{\mathcal X_{k(x_j)} }\frac {\widetilde O_{j,k}(x_j,x_k)} {\widetilde O_j(x_j)} g(x_k) \mathrm dx_k\right | \leq C \left (\int_{\mathcal X_{k} }O^{1-\delta}_{k}(x_k)g^2(x_k) \mathrm dx_k\right )^{1/2}.
\end{align}
Below we will also use that a function with bounded $\|..\|_{O,\infty}$-norm is mapped by $T$ and $\hat T$ into a function with bounded $\|..\|_{O,\infty}$-norm. This follows from 
\begin{align} \label{helpdec4}
&\sup_{x_j \in \mathcal X_{j}}  O_j^{1/2}(x_j)\left |\int_{\mathcal X_{k(x_j)} }\frac {O_{j,k}(x_j,x_k)} {O_j(x_j)} g(x_k) \mathrm dx_k\right | \leq C^* \sup_{x_k \in \mathcal X_{k}}  O_k^{1/2}(x_k)| g(x_k)|,\\
\label{helpdec4a}
&\sup_{x_j \in \mathcal X_{j}}  O_j^{1/2}(x_j)\left |\int_{\mathcal X_{k(x_j)} }\frac {\widetilde O_{j,k}(x_j,x_k)} {\widetilde O_j(x_j)} g(x_k) \mathrm dx_k\right | \leq C^* \sup_{x_k \in \mathcal X_{k}}  O_k^{1/2}(x_k)| g(x_k)|
\end{align}
with probability to one. 
We now show \eqref{helpdec1}--\eqref{helpdec4a}. The bound \eqref{helpdec4} follows directly from the last inequality in Condition [A7]. For the proof of \eqref{helpdec1} note that the left hand side of \eqref{helpdec1} can be bounded by a constant times 
 \begin{align*}
 &\int_{\mathcal X_{x_j,x_k} }\frac {1} {O_j^{1+\delta}(x_j)O_{k}(x_k)} \mathrm{d}x_j\ \mathrm{d}x_k  \int_{\mathcal X_{k} }O_{k}(x_k)g^2(x_k) \mathrm dx_k. \end{align*} 
 Thus,  \eqref{helpdec1} follows by application of the first  inequality in Condition [A7]. For the proof of \eqref{helpdec2} note that the left hand side of \eqref{helpdec2} can be bounded by a constant times 
\begin{align*}
&\left ( \sup_{x_k \in \mathcal X_{k}}  \int_{\mathcal X_{j(x_k)} }\frac {1} {O_j^{1-\delta}(x_j)O_{k}(x_k)} \mathrm{d} x_j \int_{\mathcal X_{k} }O^{1-\delta}_{k}(x_k)g^2(x_k) \mathrm dx_k\right )^{1/2}.
 \end{align*} 
 Here, \eqref{helpdec2} follows by application of the second  inequality in Condition A7.
For the proof of \eqref{helpdec1a}, \eqref{helpdec2a} and \eqref{helpdec4a} one uses \eqref{equniform} to show that the left hand sides of the equations in Condition [A7] are of order $O_P(1)$ if one replaces $O_j$ and $O_k$ by $\widetilde O_j$ and $\widetilde O_k$, respectively. Thus, one can show \eqref{helpdec1a}, \eqref{helpdec2a} and \eqref{helpdec4a} by using the same arguments as in the proofs of \eqref{helpdec1}, \eqref{helpdec2} and \eqref{helpdec4}.

We now want to show that 
\begin{align}  \label{claimadd274}\|\overline \delta_0^r(x_0)+...+\overline \delta_d^r(x_d)- \sum_{l=0}^\infty  T^l  \tau^r (x)\|_{O,\infty} = o_P(n^{-2/5}), \end{align}
where 
\begin{align*} \tau^r &=&  \Psi_d \cdot ... \cdot  \Psi_1 ( \widehat \mu^r_0 - \overline \mu_0^{*,r}) + ... + \Psi_d  ( \widehat \mu^r_{d-1}- \overline \mu_{d-1}^{*,r} ) +  ( \widehat \mu^r_{d}- \overline \mu_d^{*,r} )\end{align*} 
and where for $\delta= (\delta^*, \delta_0, ...,\delta_d)^\intercal \in \mathbb{R} \times {\cal L} $ we define $\| \delta^* + \delta_0+...+ \delta_d\|_{O,\infty}$ as $\| ( \delta_0, ...,\delta_d)^\intercal\|_{O,\infty} \vee |\delta^*|$.

Using \eqref{helpdec1}--\eqref{helpdec4a}, $\|  T \| < 1$ and the fact that $\| \hat T \|_n < \gamma$
for some $\gamma < 1$ with probability tending to one, one gets that for \eqref{claimadd274} it suffices to show  that for all choices of $c >0$
\begin{align}  \label{claimadd274a}\|\sum_{l=0}^{c \log n}  \widehat T^l  \widehat \tau^r (x)- T^l  \tau^r (x)\|_{O,\infty} = o_P(n^{-2/5}). \end{align}
For the proof of this claim it suffices to show that the norm of each summand is of order $ o_P(n^{-2/5}(\log n)^{-1})$. This can be shown by using condition B1, \eqref{helpdec1}--\eqref{helpdec4a}, and 

\begin{align} \label{helpdec1aa}
&\int_{\mathcal X_{j} }\left (\int_{\mathcal X_{k(x_j)} }\left [\frac {\widetilde O_{j,k}(x_j,x_k)} {\widetilde O_j(x_j)}- \frac {\widetilde O_{j,k}(x_j,x_k)} {\widetilde O_j(x_j)} \right]g(x_k) \mathrm dx_k\right )^2 O^{1-\delta}_j(x_j)\mathrm dx_j\\ \nonumber
& \qquad = o_P((\log n)^{-1}) \int_{\mathcal X_{k} }O_{k}(x_k)g^2(x_k) \mathrm {d}x_k,\\ \label{helpdec2aa}
&\sup_{x_j \in \mathcal X_{j}}  O_j^{1/2}(x_j)\left |\int_{\mathcal X_{k(x_j)} }\left [\frac {\widetilde O_{j,k}(x_j,x_k)} {\widetilde O_j(x_j)}- \frac {\widetilde O_{j,k}(x_j,x_k)} {\widetilde O_j(x_j)} \right] g(x_k) \mathrm dx_k\right | \\ \nonumber
& \qquad  = o_P((\log n)^{-1})  \left (\int_{\mathcal X_{k} }O^{1-\delta}_{k}(x_k)g^2(x_k) \mathrm dx_k\right )^{1/2},
\\ \label{helpdec4aa}
&\sup_{x_j \in \mathcal X_{j}}  O_j^{1/2}(x_j)\left |\int_{\mathcal X_{k(x_j)} }\left [\frac {\widetilde O_{j,k}(x_j,x_k)} {\widetilde O_j(x_j)}- \frac {\widetilde O_{j,k}(x_j,x_k)} {\widetilde O_j(x_j)} \right]g(x_k) \mathrm dx_k\right |\\ \nonumber
& \qquad  = o_P((\log n)^{-1})  \sup_{x_k \in \mathcal X_{k}}  O_k^{1/2}(x_k)| g(x_k)|.
\end{align}
Claims \eqref{helpdec1aa}--\eqref{helpdec4aa} can be shown similarly as \eqref{helpdec1}--\eqref{helpdec4a} by using additionally condition [B1].

For $r=B$ we note  that $\overline \mu_{k}^{*,r,n} - \overline \mu_{k}^{*,r} = o_P(n^{-2/5})$ because of [B1] and [B4], see also \eqref{equniform} and that the sum of the elements of $(I- \pi)^{-1}(\widehat \mu^B- \overline \mu^{B,*})$ is equal to $\sum_{l=0}^\infty  T^l  \tau^B (x)$. For $r= C$ one checks easily that $\|\sum_{l=0}^\infty  T^l  \tau^C\|_{O,\infty} = o_P(n^{-2/5})$. For the statement of the proposition it remains to show that $\|\sum_{l=1}^\infty  T^l  \tau^A\|_{O,\infty} = o_P(n^{-2/5})$ and that $
\|
\tau^r - ( \widehat \mu^r_0 - \overline \mu_0^{*,r} + ... +  \widehat \mu^r_{d}- \overline \mu_d^{*,r} )
\|_{O,\infty}= o_P(n^{-2/5})
$. For the proof of these two claims one applies condition B5. 

\subsection{Proof of Proposition \ref{prop:expansion}}
The statement of Proposition \ref{prop:expansion} follows immediately from [B3] and Proposition 1.

\subsection{Proof of Theorem \ref{theoremfbar}}

The main tool to prove this theorem is the Newton-Kantorovich theorem, see for example \citet{Deimling:85}.
Since this theorem is central in our considerations we will state it here. 
\begin{theorem}[Newton-Kantorovich theorem] \label{thmNK}Consider Banach spaces ${\cal X}, {\cal Y}$ and a  map $F:B_r(x_0)=\{x: \|x-x_0\| \leq r\}  \subset {\cal X} \mapsto {\cal Y}$ for $x_0 \in {\cal X}$ and $r > 0$.
We assume that the Fr\'echet derivative $F^\prime$ exists for $x \in B_r(x_0)$, that it is invertible and that the following conditions are satisfied
\begin{enumerate}
\item $\|F^\prime(x_0)^{-1}F(x_0)\|\leq \gamma$, 
\item $\|F^\prime(x_0)^{-1}\|\leq \beta $,
\item $\|F^\prime(x)-F^\prime(x^*)\|\leq l \|x-x^*\|$ for all $x,x^*\in B_r(x_0)$,
\item $2\gamma\beta l <1$  and  $2\gamma<r$.
\end{enumerate}
Then the equation
\[
F(x)=0
\]
has a unique solution $x^*$ in $\overline B_{2\gamma}(x_0)$ and furthermore,
$x^*$ can be approximated by Newtons iterative method 
\[
x_{k+1}=x_k - F^\prime(x_k)^{-1} F(x_k), 
\]
and it holds that
\[
\|x_k-x^*\| \leq \frac {\gamma}{2^{k-1}}q^{2^k-1}, \quad \mbox{with} \ q=2\gamma\beta l<1.
\]
\end{theorem}
We now come to the proof of Theorem \ref{theoremfbar}.
\begin{proof}[Proof of Theorem \ref{theoremfbar}]
Equation \eqref{backfitting:deltahat} can be rewritten as
\[
\widehat {\mathcal F}(\widehat \delta^*, \widehat \delta_0,\dots, \widehat \delta_d)=0,
\]
where
\[
\widehat  {\mathcal F}(f^*,f_0,\dots, f_d)(x)= \left(\widehat  {\mathcal F}_k(f^*,f_0,\dots, f_d)(x)\right)_{k=*,0,\dots,d}.
\]
with
\begin{align*}&\widehat {\mathcal F}_*(f^*, f_0,\dots, f_d)(x)=  \int_{\mathcal X} \bigg[(1+f^*) \prod_{j=0}^d  \left\{1+f_j(x_j)\right\}-1\bigg]\\ & \qquad \times \prod_{j=0}^d  \alpha_j(x_j) \widehat E(x) \mathrm dx- \int_{{\cal X}_k} \widehat \Delta_k(x_k)\mathrm d x_k,
\\
&\widehat {\mathcal F}_k(f^*, f_0,\dots, f_d)(x)=  \int_{\mathcal X_{x_k}} \bigg[(1+f^*) \prod_{j=0}^d  \left\{1+f_j(x_j)\right\}-1\bigg]\\ & \qquad \times \prod_{j=0}^d  \alpha_j(x_j) \widehat E(x) \mathrm dx_{-k}- \widehat \Delta_k(x_k)- \widehat {\mathcal F}_*(f^*, f_0,\dots, f_d)(x)
\end{align*}
for $k=0,\dots,d$. Note that $\int_{{\cal X}_k} \widehat \Delta_k(x_k)\mathrm d x_k$ does not depend on $k$.

We define an additional   operator $\mathcal F$ by the following equations
\[
\mathcal F(f^*,f_0,\dots, f_d)(x)= \left(\mathcal F_k(f^*,f_0,\dots, f_d)(x)\right)_{k=*,0,\dots,d}
\]
with 
 \begin{align*}&
\mathcal F_*(f^*,f_0,\dots, f_d)(x)= \int_{\mathcal X} \left[\left\{1+f^*\right\} \prod_{j=0}^d  \left\{1+f_j(x_j)\right\}-1\right]\\
& \qquad \times \prod_{j=0}^d  \alpha_j(x_j) E(x) \mathrm dx,\\
&
\mathcal F_k(f^*,f_0,\dots, f_d)(x)= \int_{\mathcal X_{x_k}} \left[\left\{1+f^*\right\} \prod_{j=0}^d  \left\{1+f_j(x_j)\right\}-1\right]\\
& \qquad \times \prod_{j=0}^d  \alpha_j(x_j) E(x) \mathrm dx_{-k}- \mathcal F_*(f^*,f_0,\dots, f_d)(x)
\end{align*}
for ${k=0,\dots,d}$.

Note that $\mathcal F(0)=0$.
The  Fr\'echet derivatives of $\widehat {\mathcal F}$ and $\mathcal F$ in $0$ are
\begin{align*}
\widehat {\mathcal F}_*^{\prime} (0)(f)&= \int_{\mathcal X}\left (f^*+\sum_{j=0}^d f_j(x_j)\right )\alpha(x)\widehat E(x) \mathrm dx , \\
 {\mathcal F}_*^{\prime} (0)(f)&= \int_{\mathcal X}\left (f^*+\sum_{j=0}^d f_j(x_j)\right )\alpha(x) E(x) \mathrm dx,\\
\widehat {\mathcal F}_k^{\prime} (0)(f)&= \int_{\mathcal X_{x_k}}\left (f^*+\sum_{j=0}^d f_j(x_j)\right )\alpha(x)\widehat E(x) \mathrm dx_{-k} \\ & \qquad - \widehat {\mathcal F}_*^{\prime} (0)(f), \\
 {\mathcal F}_k^{\prime} (0)(f)&= \int_{\mathcal X_{x_k}}\left (f^*+\sum_{j=0}^d f_j(x_j)\right )\alpha(x) E(x) \mathrm dx_{-k}  \\ & \qquad -  {\mathcal F}_*^{\prime} (0)(f)\end{align*}
for $k=0,\dots,d$.

The main idea of our proof is to apply the Newton-Kantorovich theorem, Theorem \ref{thmNK}, with the mapping $F=\widehat {\mathcal F }$ and  norm $\|(f_0,...,f_d)\|_{O,\infty}\vee |f^*|$ which in abuse of notation we also denote by $\|(f^*, f_0,...,f_d)\|_{O,\infty}$. As  starting point $x_0$ we choose $x_0=(\overline \delta^*,\overline \delta)$. In our application of Theorem \ref{thmNK}, the spaces ${\cal X}$ and ${\cal Y}$ are equal to $\mathbb{R} \times \{(f_0,...,f_d)^\intercal : \|(f_0,...,f_d)^\intercal\|_{O,\infty} < \infty, \int f_j(x_j) \widetilde O_j(x_j) \mathrm {d} x_j = 0 $ for $ j=0,...,d\}$. We consider ${\mathcal F}$ and $\widehat {\mathcal F}$ as operators from ${\cal X}$ to ${\cal X}$. Note that  ${\cal  FX} \subset {\cal  X}$ because of [B6] and the last assumption of [A7]. Note that we get from [B6] and the last assumption of [A7] that $$ \int_{x \in \mathcal {X} }\frac {1} {\prod_{k =0}^d O_{k}(x_k)^{1/2}} { O(x)}  \mathrm dx = O(1).$$ 
Similarly, one uses [B6] and the last assumption of [A7] to show  that ${\cal  \widehat FX} \subset {\cal  X}$, with probability tending to one. 

We will show that 
\begin{align}\label{cond1}
\norm{\widehat {\mathcal F}\left((\overline \delta^*,\overline \delta) \right)}_{O,\infty}=o_p(n^{-2/5}),
\end{align}
and  that $\widehat { \cal F}^\prime$ is locally Lipschitz around $0$, i.e., that there exist constants $r^*, C$ such that with probability tending to one
\begin{align}\label{cond2}
\norm{\widehat {\mathcal F}^\prime(g)(f)-\widehat {\mathcal F}^\prime(g^*)(f)}_{O,\infty}\leq C \norm{g-g^*}_{O,\infty} \norm{f}_{O,\infty} \quad \text{for all} \quad g, g^*\in B_{r^*}(0).
\end{align}
Furthermore, we will show, that
\begin{align}\label{cond3}
&\mathcal F^\prime(0)  \text{ is invertible, with}  \ \norm{\mathcal F^\prime(0)^{-1}}_{O,\infty}<C^*, \quad \mbox{for some }C^*>0 .
\end{align}
We now argue that by application of the Newton-Kantorovich theorem \eqref{cond1}-\eqref{cond3} imply that
\begin{align}\label{condclaim}\norm{(\overline \delta^*,\overline \delta)- (\widehat \delta^*,\widehat \delta)}_{O,\infty}=o_p(n^{-2/5}).\end{align}
This implies the statement of the theorem.

We now show  that  \eqref{cond1}-\eqref{cond3}  imply \eqref{condclaim}. Since $\norm{(\overline \delta^*,\overline \delta)}_{O,\infty}=o_P(1)$,  the inequality \eqref{cond2}
also holds with a constant $r$ for all $g, g^*\in B_{r}\left ((\overline \delta^*,\overline \delta)\right)$
with probability tending to one. This gives condition (c) of the Newton-Kantorovich theorem. 

Furthermore, by application of \eqref{equniform} we get that $\norm{\widehat{\mathcal F}^\prime(0)-\mathcal F^\prime(0)}_{O,\infty}=o_P(1)$.
This together with 
 $\norm{(\overline \delta^*,\overline \delta)}_{O,\infty}=o_P(1)$  and \eqref{cond2}
gives 
\[
\norm{\widehat{\mathcal F}^\prime\left((\overline \delta^*,\overline \delta)\right)- \mathcal F^\prime(0)}_{O,\infty}=o_p(1).
\]
Therefore with probability tending to one, condition \eqref{cond3} also holds if $\mathcal F^\prime(0)$ is replaced by $\widehat {\mathcal F}^\prime\left ((\overline \delta^*,\overline \delta)\right )$.
Thus, we get from \eqref{cond1}-\eqref{cond3} that conditions (a)--(d) of the Newton-Kantorovich theorem are fulfilled with probability tending to one, with $\gamma=o_P(n^{-2/5})$. This shows \eqref{condclaim}.

It remains to show \eqref{cond1}, \eqref{cond2} and \eqref{cond3}. For the proof of \eqref{cond1} note that
 $\norm {(\overline \delta^*,\overline \delta)}_{O,\infty}=o_p(n^{-1/5})$ and that $\widehat {\mathcal F}^\prime$
is Lipschitz. A first order Taylor expansion yields
\[
\widehat {\mathcal F}(\overline \delta)=
\widehat {\mathcal F}(0) +  \widehat {\mathcal F}^\prime(0) \left ((\overline \delta^*,\overline \delta) \right)+o_p(n^{-2/5}).
\]
Equation \eqref{cond1} follows from $\widehat {\mathcal F}(0) +  \widehat {\mathcal F}^\prime(0) \left ((\overline \delta^*,\overline \delta) \right)=\widehat {\mathcal F}(0) - \widehat {\mathcal F}(0)=0 $.

Claim  \eqref{cond2} follows directly from assumption [B6].

For the proof of \eqref{cond3} we 
have to show  that $\mathcal F^\prime(0)$ is invertible. For the proof of this claim we start by showing that it is bijective.
For the proof of injectivity,
assume that $\mathcal F^\prime(0)(f)=0$ for some $f=(f^*,f_0,...,f_d)^\intercal\in {\cal X}$. We will show that this implies that $f=0$. It holds that
\begin{align*}
&\int_{\mathcal X}  \left ( f^* + \sum_{ j=0}^{d}f_j(x_j)\right ) \alpha(x)E(x)\mathrm dx=0,\\ &\int_{\mathcal X_{x_k}}  \left ( f^* + \sum_{ j=0}^{d}f_j(x_j)\right ) \alpha(x)E(x)\mathrm dx_{-k}=0, \quad \text{for all} \ k=0,\dots,d.
\end{align*}

With $\bar f_j (x_j) = f_j (x_j)  - \int f_j(u_j) v(u) \mathrm{d} u $ and $\bar f^* = f^* + \sum_{j=0}^d \int f_j(u_j) v(u) \mathrm{d} u$ for $v(u) = \alpha (u)  E(u) / \int \alpha(s)  E(s) \mathrm{d}s$ this implies that \[
\bar f^* \int_{\mathcal X}  \alpha(x)E(x)\mathrm dx = 0\] and thus it holds that  $\bar f^*=0$. Furthermore, we get that for $ k=0,\dots,d$
\[
0= \int_{{\mathcal X}_k} \bar  f_k(x_k)\int_{{\mathcal X}_{x_k}} \left (\sum_{ j=0}^{d}\bar f_j(x_j)\right )\alpha(x)E(x)\mathrm dx_{-k} \mathrm dx_k= \int_{{\mathcal X}} \bar  f_k(x_k) \left (\sum_{ j=0}^{d}\bar f_j(x_j)\right )\alpha(x)E(x)\mathrm dx.\]
By summing these terms up over $k$, we get that 
\[
\int_{\mathcal X}  \left\{\sum_{ j=0}^{d} \bar f_j(x_j)\right\}^2\alpha(x)E(x)\mathrm dx=0,
\]
which implies that
\[
\sum_{ j=0}^{d}\bar f_j(x_j)=0, \quad \ \text{a.e.} \ \text{on} \ \mathcal X.
\]
By application of \eqref{normbound} this implies that $\bar f_j$ and $f_j$ are constant functions. Because of $\int f_j(u_j) \alpha(u) \widehat E(u) \mathrm{d} u= 0$ this implies $f=0$. \\

Now we check that $\mathcal F^\prime(0)$ is surjective.
Consider $g=(g^*,g_0,...,g_d)^\intercal \in {\cal X}$ 
such 
that \begin{equation} \label{eqortho} \langle \mathcal F^\prime(0)(f),g\rangle = {\cal F}_*^{\prime} (0)(f) g^*+\sum_{k=0}^d \int_{{\cal X}_k} {\cal F}_k^{\prime} (0)(f)(x_k) g_k(x_k) \mathrm d x_k=0\end{equation}  for all $f=(f^*,f_0,...,f_d)^\intercal\in {\cal X}$. We will show that then $g=0$. This implies that $g=0$ is the only element in ${\cal Y}$ that is perpendicular to the range space of  $\mathcal F^\prime(0)$. 
Since  $\mathcal F^\prime(0)$ is linear, this shows that $\mathcal F^\prime(0)$ is surjective.

From \eqref{eqortho} one gets with the choice  $f_k = g_k$ and  $f^*= g^*$ that 
\begin{align*}
0 &=&(g^*)^2
\int_{\mathcal X} \alpha(x)E(x)\mathrm dx
+
\sum_{k=0}^d  \int_{{\mathcal X}_k} g_k(x_k)\int_{{\mathcal X}_{x_k}} \left ( \sum_{ j=0}^{d} g_j(x_j)\right )\alpha(x)E(x)\mathrm dx_{-k} \mathrm dx_k\\
&= &(g^*)^2
\int_{\mathcal X} \alpha(x)E(x)\mathrm dx
+
  \int_{{\mathcal X}}  \left ( \sum_{ j=0}^{d} g_j(x_j)\right )^2\alpha(x)E(x)\mathrm dx.
\end{align*}
With exactly the same arguments as for the injectivity we conclude that $g=0$ and that $g^*=0$. Thus, we have shown that
$\mathcal F^\prime(0)$ is invertible. \\

It remains to show that $\mathcal F^\prime(0)^{-1}$ is bounded. By the  bounded inverse theorem for this claim  it suffices to show that $\mathcal F^\prime(0)$  is bounded. Boundedness of $\mathcal F^\prime(0)$ can be shown by application of 
\eqref{helpdec4a}. This concludes the proof of Theorem  \ref{theoremfbar}.
\end{proof}

\section{Simulation results}

\subsubsection{Simulation study of \cite{Honda:05} and \cite{Lin:etal:16}}
In this section we present the simulation results of our proposed smooth
backfitting estimator in the setting of \cite{Honda:05} and \cite{Lin:etal:16}.
We refer to\cite{Lin:etal:16} for notation and simulation setting. We only
present the results for Model 1 and Model 2 (not Model 3) of \cite{Lin:etal:16},
because in Model 3, \cite{Lin:etal:16} estimatev linear $\eta's$ by over-smoothing
(bandwidth is double the size of the support) their local linear estimator.
Hence, if over-smoothing is done because of the knowledge of an underlying
linear function, then using the (linear) Cox proportional hazard model is
actually more appropriate.

In each model, the last two rows of Table \ref{tab:lin} show the results of a simulation
(with 1000 repetitions) run by us; the first five  rows are copied from Table \ref{tab:lin} in
\cite{Lin:etal:16}. We conclude that our proposed estimator, in the setting
considered, shares top performance together with \cite{Lin:etal:16}, with all
other estimators showing inferior performance.

\begin{table}[ht]
\footnotesize
\begin{tabular}{|l|ccc|ccc|c}
\hline
\hline
&\multicolumn{3}{c}{$z=0.5$} &\multicolumn{3}{|c|} {$z=-0.5)$} \\
\hline
Method &Bias &SD &RMSE& Bias &SD& RMSE\\
\hline
 \multicolumn{7}{|c|}{Model 1}\\
 \hline
\cite{Huang:99} & 0.0409 &0.2006& 0.2048 &0.0467 &0.2181 &0.2230\\
\cite{Fan:etal:97} &-0.1006 &0.2502& 0.2697& -0.1267 &0.1879& 0.2265\\
\cite{Linton:etal:03}& -0.0350 &0.2949 &0.2970 &0.0400& 0.2702 &0.2731\\
\cite{Honda:05} &0.1130 &0.2588 &0.2824& -0.1090& 0.2569& 0.2791\\
 \cite{Lin:etal:16} (copied)& -0.0381& 0.1380 &0.1432 &0.0881& 0.1558 &0.1790\\
 \cite{Lin:etal:16} (re-run) &0.0389& 0.0990& 0.1061& 0.0394 &0.1042 &0.1111\\
SBF &0.0803 &0.1018 &0.1293&0.08148& 0.10253& 0.1306\\
 \hline
 \hline
 \multicolumn{7}{|c|}{Model 2}\\
 \hline
\cite{Huang:99} &0.0457 &0.2066 &0.2116& -0.0014 &0.2116 &0.2116\\
\cite{Fan:etal:97} &0.0168& 0.2995& 0.2999 &0.1662& 0.2506 &0.3007\\
\cite{Linton:etal:03}&-0.0560 &0.3209 &0.3257& 0.0780& 0.3049& 0.3147\\
\cite{Honda:05} &0.1040& 0.2302 &0.2526& -0.1060& 0.2408& 0.2631\\
 \cite{Lin:etal:16} (copied)& 0.0621& 0.0606 &0.0868 &-0.0372& 0.0901& 0.0975\\
 \cite{Lin:etal:16} (re-run)& -0.0660& 0.1488 &0.1630& -0.0697 &0.1450& 0.1606\\
SBF &-0.0038 &0.1463 &0.1463& -0.0098 &0.1376& 0.1375\\
 \hline
 \hline
\end{tabular}
\caption{Summaries of the simulation results in the setting of \cite{Honda:05}
and \cite{Lin:etal:16}. In each model, the last two rows show the results of
a simulation (with 1000 repetitions) run by us; the first five rows are copied
from Table 1 in \cite{Lin:etal:16}.}
\label{tab:lin}
\end{table}

\subsection{Additional simulation results in the setting of the main manuscript}
In this section we present additional simulation results for the simulation
study conducted in the main manuscript. We refer to the main manuscript
for notation and simulation setting.
\subsection{Dimension $d=2$}
\begin{table}[ht]
\label{tab:breakdownd:d3}
\centering
\begin{tabular}{|c|ccc|ccc|}
\hline
\multicolumn{7}{|c|}{Number of breakdowns in \citet{Lin:etal:16}for $d=2$} \\ \multicolumn{7}{|c|} {(out of 200 simulations)} \\
\hline
 &\multicolumn{3}{c|}{Model 1}&\multicolumn{3}{c|}{Model 2}\\
  \hline
&$\rho=0$ &$\rho=0.5$&$\rho=0.8$&$\rho=0$ &$\rho=0.5$&$\rho=0.8$\\
  \hline
n=200 &   0 &   1 &   2 &   0 &   0 &   0 \\ 
  n=500 &   0 &   0 &   0&   0 &   0 &   0 \\    \hline
\end{tabular}
\caption{Number of breakdowns in the algorithm of \citet{Lin:etal:16}out of 200 simulations for dimension $d=2$.}
\end{table}

\begin{table}[ht]
\footnotesize
\centering
\begin{tabular}{|c|c|c|ccc|ccc|}
\hline
\multicolumn{9}{|c|}{ISE values for Model 1 (d=2)} \\
\hline
&&&\multicolumn{3}{c|}{\citet{Lin:etal:16}}&\multicolumn{3}{c|}{SBF}\\
  \hline
sample size &correlation & component & mean & median & sd & mean & median & sd \\ 
  \hline
&$\rho=0 $&$k=1$ & 0.074 & 0.059 & 0.059 & 0.052 & 0.039 & 0.048 \\ 
&& $k=2$ &0.087 & 0.068 & 0.063 & 0.058 & 0.045 & 0.048 \\  \hline
n=200&$\rho=0.5$  &$k=1 $&0.077 & 0.060 & 0.055 & 0.061 & 0.047 & 0.048 \\ 
&  &$k=2 $&0.087 & 0.068 & 0.069 & 0.058 & 0.044 & 0.044 \\ \hline
  &$\rho=0.8$&$k=1 $&  0.128 & 0.094 & 0.112 & 0.089 & 0.075 & 0.062 \\
  &&$k=2 $& 0.136 & 0.096 & 0.120 & 0.088 & 0.071 & 0.057 \\   
   \hline
    \hline
&$\rho=0 $&$k=1$ & 0.026 & 0.022 & 0.019 & 0.019 & 0.015 & 0.016 \\ 
&& $k=2$ &0.033 & 0.027 & 0.026 & 0.026 & 0.021 & 0.020 \\ \hline
n=500&$\rho=0.5$  &$k=1 $&0.027 & 0.021 & 0.020 & 0.023 & 0.018 & 0.017 \\ 
&  &$k=2 $&0.032 & 0.026 & 0.023 & 0.029 & 0.025 & 0.021\\ \hline
  &$\rho=0.8$&$k=1 $&  0.041 & 0.030 & 0.040 & 0.048 & 0.043 & 0.028 \\ 
  &&$k=2 $&  0.051 & 0.037 & 0.048 & 0.055 & 0.051 & 0.030 \\    \hline
\end{tabular}
\caption{Summaries for integrated squared errors. Simulations where the algorithm of  \citet{Lin:etal:16}broke down are taken out. }
\end{table}

\begin{table}[ht]
\footnotesize
\centering
\begin{tabular}{|c|c|c|ccc|ccc|}
\hline
\multicolumn{9}{|c|}{ISE values for Model 2 (d=2)} \\
\hline
&&&\multicolumn{3}{c|}{\citet{Lin:etal:16}}&\multicolumn{3}{c|}{SBF}\\
  \hline
sample size &correlation & component & mean & median & sd & mean & median & sd \\ 
  \hline
&$\rho=0 $&$k=1$ & 0.120 & 0.109 & 0.068 & 0.132 & 0.105 & 0.085 \\ 
&& $k=2$  &0.080 & 0.065 & 0.059 & 0.066 & 0.052 & 0.058 \\ \hline
n=200&$\rho=0.5$  &$k=1 $&0.151 & 0.127 & 0.088 & 0.159 & 0.124 & 0.113 \\ 
&  &$k=2 $&0.079 & 0.064 & 0.060 & 0.072 & 0.063 & 0.047 \\  \hline
  &$\rho=0.8$&$k=1 $& 0.303 & 0.236 & 0.242 & 0.172 & 0.153 & 0.105 \\ 
  &&$k=2 $&  0.140 & 0.089 & 0.146 & 0.090 & 0.080 & 0.062 \\   
  \hline
    \hline
&$\rho=0 $&$k=1$ & 0.086 & 0.078 & 0.041 & 0.105 & 0.092 & 0.053 \\ 
&& $k=2$  &0.031 & 0.025 & 0.021 & 0.027 & 0.022 & 0.020 \\ 
n=500&$\rho=0.5$  &$k=1 $&0.104 & 0.099 & 0.047 & 0.111 & 0.098 & 0.055 \\ 
&  &$k=2 $&0.030 & 0.025 & 0.019 & 0.036 & 0.031 & 0.021 \\  \hline
  &$\rho=0.8$&$k=1 $&  0.199 & 0.185 & 0.092 & 0.135 & 0.124 & 0.062 \\ 
  &&$k=2 $&0.046 & 0.036 & 0.032 & 0.061 & 0.056 & 0.033 \\   \hline
\end{tabular}
\caption{Summaries for integrated squared errors. Simulations where the algorithm of  \citet{Lin:etal:16}broke down are taken out. }
\end{table}

\begin{figure}[h!]

%\centering
%\begin{subfigure}[t]{0.213cm}
%\centering
\includegraphics[width=13cm]{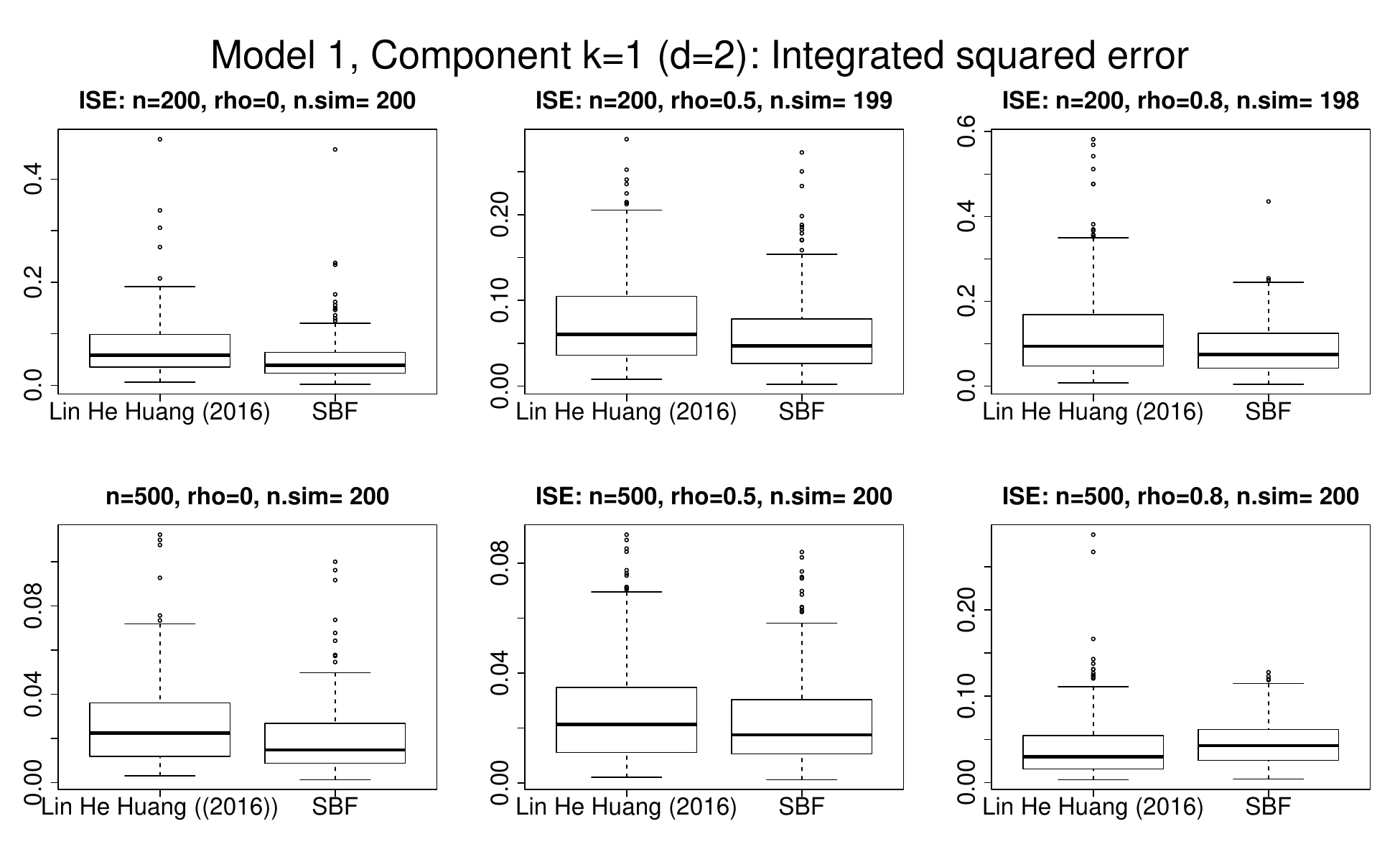} 
%\caption{Generic} \label{fig:timing1}
%\end{subfigure}

%\begin{subfigure}[t]{0.213cm}
%\centering
\includegraphics[width=13cm]{boxplot_ise312.pdf} 
%\caption{Competitors} \label{fig:timing2}
%\end{subfigure}
\caption{Boxplots of the integrated squared errors. Simulations where the algorithm of  \citet{Lin:etal:16}broke down are taken out. The value $n.sim$ is the number of successful simulations, i.e.,  200 minus number of break downs.}
\label{fig:boxplot:ise:d31}
\end{figure}

\begin{figure}[h!]
\includegraphics[width=13cm]{boxplot_ise321.pdf} 
%\caption{Generic} \label{fig:timing1}
%\end{subfigure}

%\begin{subfigure}[t]{0.213cm}
%\centering
\includegraphics[width=13cm]{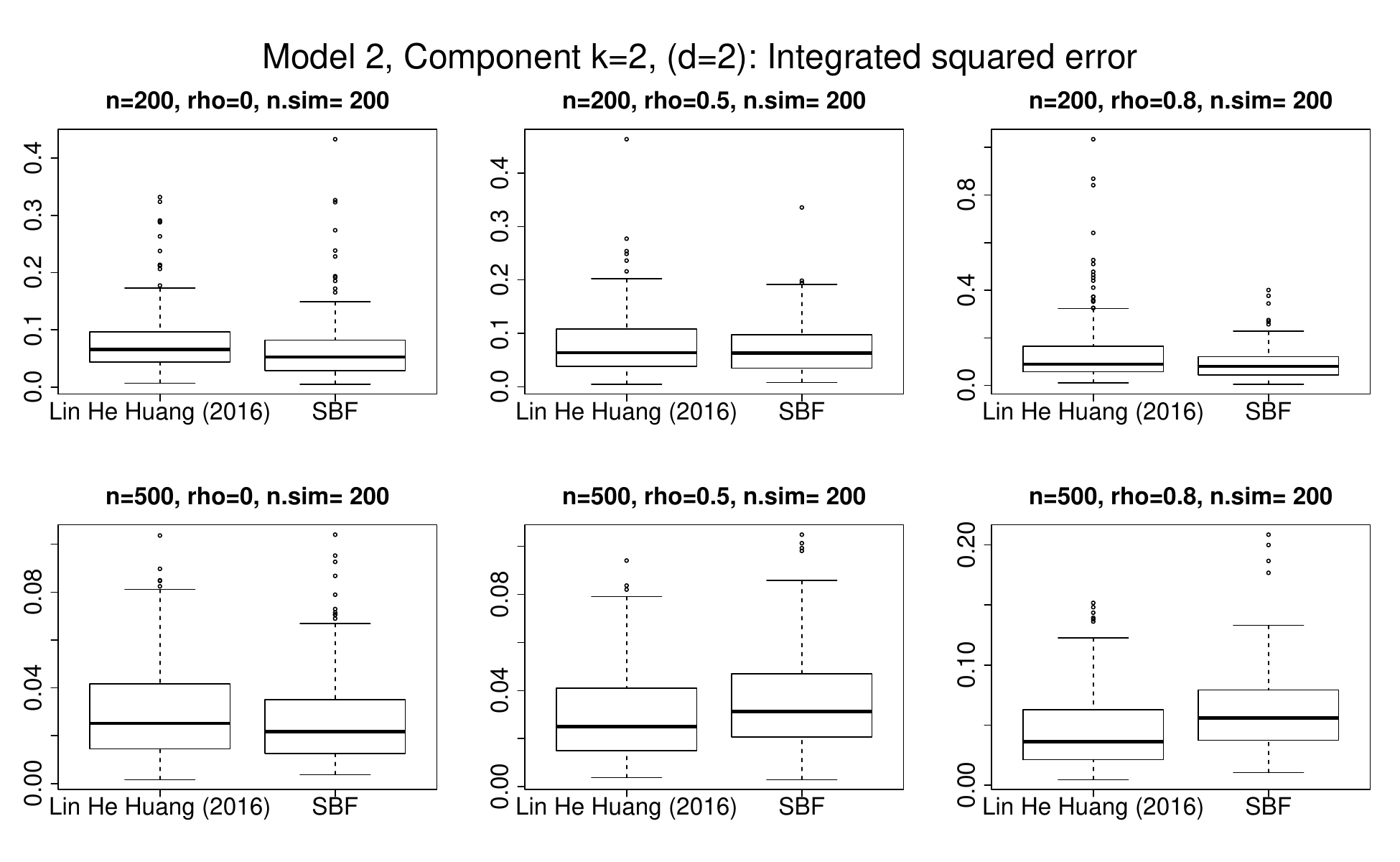} 
\caption{Boxplots of the integrated squared errors. Simulations where the algorithm of  \citet{Lin:etal:16}broke down are taken out. The value $n.sim$ is the number of successful simulations, i.e.,  200 minus number of break downs.}
\label{fig:boxplot:ise:d32}
\end{figure}

%
%%\begin{subfigure}[t]{0.213cm}
%%\centering
%\begin{figure}[h!]
%\includegraphics[width=13cm]{plot3208.pdf} 
%\caption{Simualtion 17/200 is one example where the \citet{Lin:etal:16} estimator performes bad under high correlation already in dimension $d=2$}
%\label{fig:3208}
%\end{figure}

\clearpage

\subsubsection{Dimension $d=9$}

\begin{table}[h!]
\label{tab:breakdownd:d10}
\centering
\begin{tabular}{|c|ccc|ccc|}
\hline
\multicolumn{7}{|c|}{Number of breakdowns in \citet{Lin:etal:16}  for $d=9$} \\ \multicolumn{7}{|c|} {(out of 200 simulations)} \\
\hline
 &\multicolumn{3}{c|}{Model 1}&\multicolumn{3}{c|}{Model 2}\\
  \hline
&$\rho=0$ &$\rho=0.5$&$\rho=0.8$&$\rho=0$ &$\rho=0.5$&$\rho=0.8$\\
  \hline
n=200 &   4 &   5 &  48 &0 &   0 &   9 \\ 
  n=500 &   0 &   1 &   1 &   0 &   0 &   0 \\    \hline
\end{tabular}
\caption{Number of breakdowns in the algorithm of \citet{Lin:etal:16}out of 200 simulations for dimension $d=9$.}
\end{table}

\begin{table}[ht]
\footnotesize
\centering
\begin{tabular}{|c|c|c|ccc|ccc|}
\hline
\multicolumn{9}{|c|}{ISE values for Model 1 (d=9)} \\
\hline
&&&\multicolumn{3}{c|}{\citet{Lin:etal:16}}&\multicolumn{3}{c|}{SBF}\\
  \hline
sample size &correlation & component & mean & median & sd & mean & median & sd \\ 
  \hline
&$\rho=0 $&$k=odd$ & 0.119 & 0.095 & 0.100 & 0.056 & 0.043 & 0.049 \\ 
&& $k=even$ & 0.166 & 0.124 & 0.154 & 0.063 & 0.048 & 0.057 \\  \hline
n=200&$\rho=0.5$  &$k=odd $& 0.164 & 0.121 & 0.151 & 0.072 & 0.057 & 0.056 \\ 
&  &$k=even $& 0.228 & 0.155 & 0.256 & 0.078 & 0.061 & 0.060 \\ \hline
  &$\rho=0.8$&$k=odd $&0.362 & 0.256 & 0.534 & 0.117 & 0.096 & 0.090 \\
  &&$k=even $& 0.443 & 0.282 & 0.476 & 0.137 & 0.113 & 0.106 \\
   \hline
    \hline
&$\rho=0 $&$k=odd$ & 0.032 & 0.026 & 0.024 & 0.022 & 0.016 & 0.019 \\ 
&& $k=even$ & 0.041 & 0.032 & 0.032 & 0.027 & 0.022 & 0.021 \\ \hline
n=500&$\rho=0.5$  &$k=odd $& 0.041 & 0.032 & 0.031 & 0.032 & 0.027 & 0.022 \\ 
&  &$k=even $& 0.048 & 0.035 & 0.040 & 0.038 & 0.032 & 0.026 \\  \hline
  &$\rho=0.8$&$k=odd $ & 0.091 & 0.063 & 0.087 & 0.061 & 0.052 & 0.041 \\ 
  &&$k=even $& 0.107 & 0.071 & 0.112 & 0.077 & 0.068 & 0.047 \\   \hline
\end{tabular}
\caption{Summaries for integrated squared errors. Simulations where the algorithm of  \citet{Lin:etal:16}broke down are taken out. }
\end{table}

\begin{table}[ht]
\footnotesize
\centering
\begin{tabular}{|c|c|c|ccc|ccc|}
\hline
\multicolumn{9}{|c|}{ISE values for Model 2 (d=9)} \\
\hline
&&&\multicolumn{3}{c|}{\citet{Lin:etal:16}}&\multicolumn{3}{c|}{SBF}\\
  \hline
sample size &correlation & component & mean & median & sd & mean & median & sd \\ 
  \hline
&$\rho=0 $&$k=odd$ & 0.302 & 0.284 & 0.122 & 0.216 & 0.194 & 0.115 \\ 
&& $k=even$  & 0.134 & 0.116 & 0.083 & 0.104 & 0.086 & 0.077 \\  \hline
n=200&$\rho=0.5$&$k=odd$ & 0.372 & 0.343 & 0.167 & 0.218 & 0.184 & 0.143 \\ 
&  &$k=even $& 0.133 & 0.107 & 0.099 & 0.126 & 0.106 & 0.085 \\ \hline
  &$\rho=0.8$&$k=odd $& 0.821 & 0.664 & 0.561 & 0.252 & 0.210 & 0.163 \\ 
  &&$k=even $& 0.297 & 0.196 & 0.331 & 0.211 & 0.184 & 0.124 \\ 
  \hline
    \hline
&$\rho=0 $&$k=odd$ & 0.286 & 0.281 & 0.063 & 0.166 & 0.156 & 0.062 \\ 
&& $k=even$  & 0.109 & 0.104 & 0.042 & 0.057 & 0.051 & 0.033 \\ 
\hline
n=500&$\rho=0.5$ &$k=odd$& 0.311 & 0.310 & 0.076 & 0.157 & 0.144 & 0.073 \\
&  &$k=even $& 0.080 & 0.072 & 0.042 & 0.072 & 0.063 & 0.042 \\  \hline
  &$\rho=0.8$&$k=odd $& 0.571 & 0.540 & 0.205 & 0.189 & 0.166 & 0.096 \\ 
  &&$k=even $ & 0.085 & 0.064 & 0.075 & 0.150 & 0.141 & 0.064 \\   \hline
\end{tabular}
\caption{Summaries for integrated squared errors. Simulations where the algorithm of  \citet{Lin:etal:16}broke down are taken out. }
\end{table}

\begin{figure}[h!]

%\centering
%\begin{subfigure}[t]{0.213cm}
%\centering
\includegraphics[width=13cm]{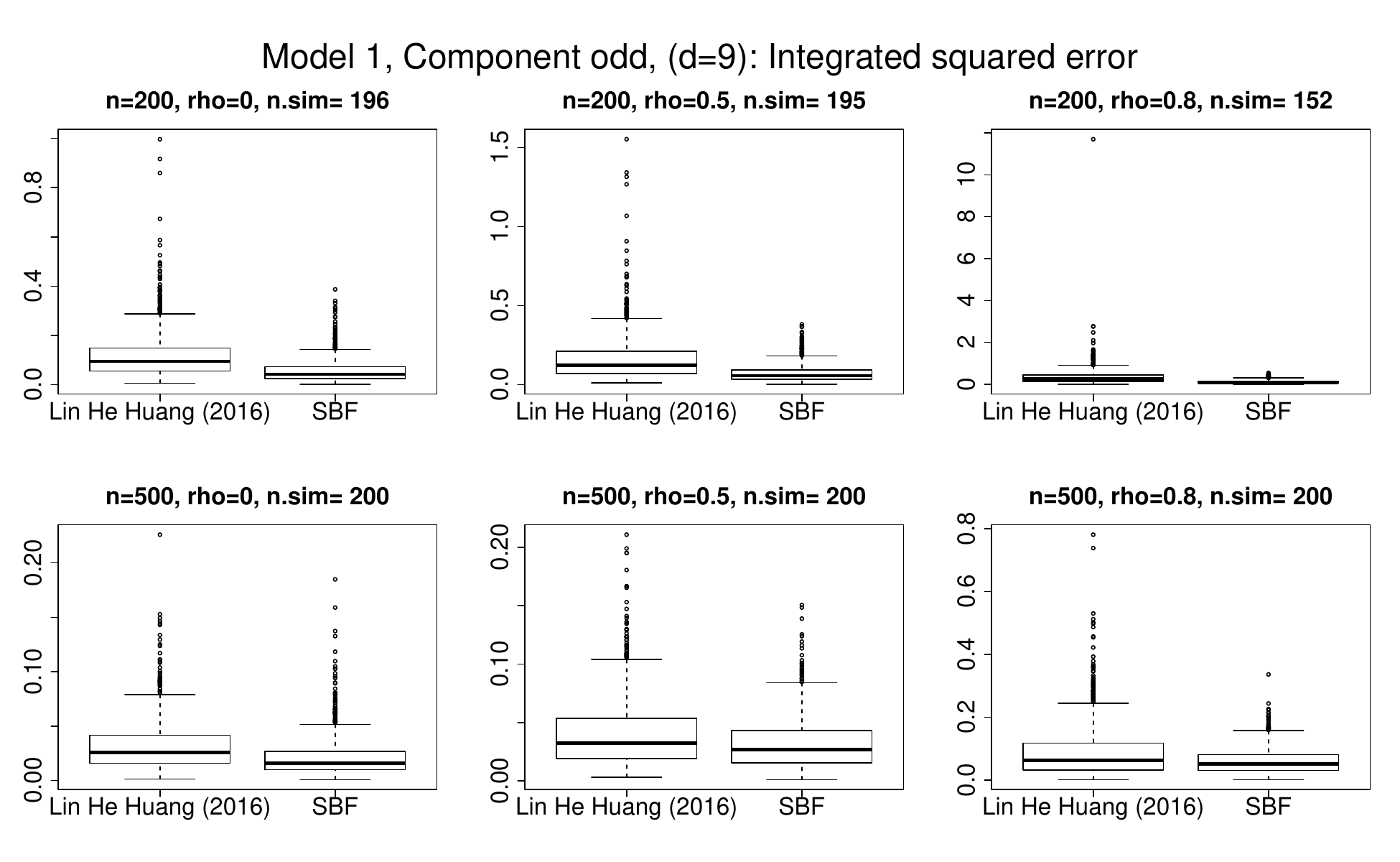} 
%\caption{Generic} \label{fig:timing1}
%\end{subfigure}

%\begin{subfigure}[t]{0.213cm}
%\centering
\includegraphics[width=13cm]{boxplot_ise1012.pdf} 
%\caption{Competitors} \label{fig:timing2}
%\end{subfigure}
\caption{Boxplots of the integrated squared errors in Model 1. Simulations where the algorithm of  \citet{Lin:etal:16}broke down are taken out. The value $n.sim$ is the number of successful simulations, i.e.,  200 minus number of break downs.}
\label{fig:boxplot:ise:d101}
\end{figure}

\begin{figure}[h!]
\includegraphics[width=13cm]{boxplot_ise1021.pdf} 
%\caption{Generic} \label{fig:timing1}
%\end{subfigure}

%\begin{subfigure}[t]{0.213cm}
%\centering
\includegraphics[width=13cm]{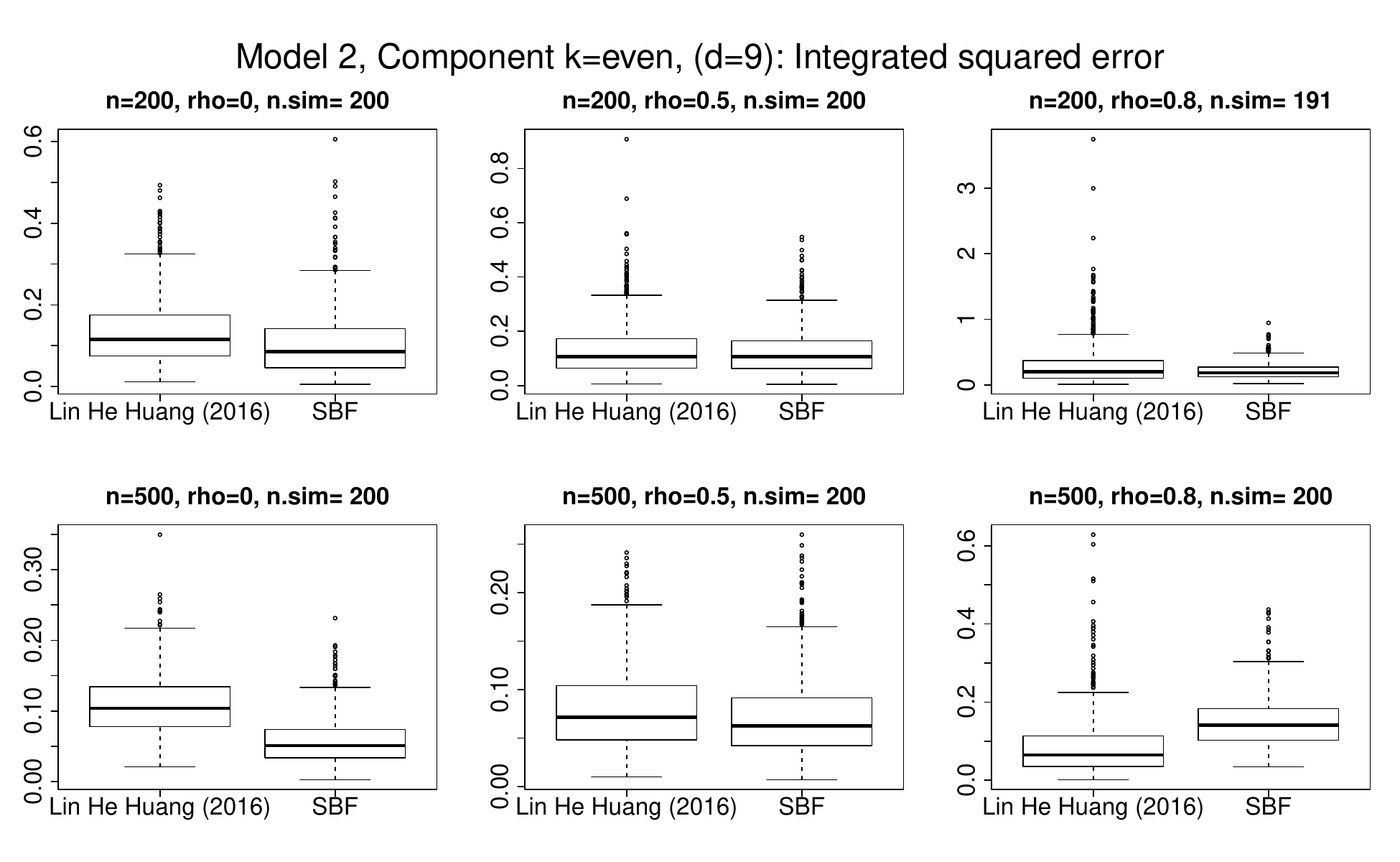} 
\caption{Boxplots of the integrated squared errors in Model 2. Simulations where the algorithm of  \citet{Lin:etal:16}broke down are taken out. The value $n.sim$ is the number of successful simulations, i.e.,  200 minus number of break downs.}
\label{fig:boxplot:ise:d102}
\end{figure}

%
%%\begin{subfigure}[t]{0.213cm}
%%\centering
%\begin{figure}
%\includegraphics[width=13cm]{plot10108.pdf} 
%\caption{Simualtion 41/200 components $k=5,6$ is one example where the \citet{Lin:etal:16} estimator performes bad under high correlation already in dimension $d=9$}
%\label{fig:10208}
%\end{figure}

\clearpage
\subsubsection{Dimension $d=30$}

\begin{table}[h!]
\label{tab:breakdownd:d30}
\centering
\begin{tabular}{|c|ccc|ccc|}
\hline
\multicolumn{7}{|c|}{Number of breakdowns in \citet{Lin:etal:16} for $d=30$} \\ \multicolumn{7}{|c|} {(out of 200 simulations)} \\
\hline
 &\multicolumn{3}{c|}{Model 1}&\multicolumn{3}{c|}{Model 2}\\
  \hline
&$\rho=0$ &$\rho=0.5$&$\rho=0.8$&$\rho=0$ &$\rho=0.5$&$\rho=0.8$\\
  \hline
n=200 & 106 & 176 & 200 &   2 &  67 & 96\\ 

  n=500 &   2 &   2 &  47 &   0 &   0&  7 \\ 
  \hline
\end{tabular}
\caption{Number of breakdowns in the algorithm of \citet{Lin:etal:16} out of 200 simulations for dimension $d=30$.}
\end{table}

\begin{table}[ht]
\footnotesize
\centering
\begin{tabular}{|c|c|c|ccc|ccc|}
\hline
\multicolumn{9}{|c|}{ISE values for Model 1 (d=30)} \\
\hline
&&&\multicolumn{3}{c|}{\citet{Lin:etal:16}}&\multicolumn{3}{c|}{SBF}\\
  \hline
sample size &correlation & component & mean & median & sd & mean & median & sd \\ 
  \hline
&$\rho=0 $&$k=odd$ &0.406 & 0.295 & 0.368 & 0.078 & 0.061 & 0.065 \\ 
&& $k=even$ & 0.741 & 0.571 & 0.635 & 0.098 & 0.082 & 0.072 \\ \hline
n=200&$\rho=0.5$  &$k=odd $& 0.664 & 0.481 & 0.625 & 0.107 & 0.091 & 0.079 \\ 
&  &$k=even $ & 1.092 & 0.759 & 1.024 & 0.128 & 0.099 & 0.104 \\\hline
  &$\rho=0.8$&$k=odd $&NA& NA & NA & NA & NA & NA \\
  &&$k=even $& NA & NA& NA & NA & NA&NA\\
   \hline
    \hline
&$\rho=0 $&$k=odd$ &0.073 & 0.055 & 0.059 & 0.030 & 0.023 & 0.023 \\ 
&& $k=even$ & 0.146 & 0.108 & 0.127 & 0.048 & 0.041 & 0.031 \\ \hline
n=500&$\rho=0.5$  &$k=odd $& 0.106 & 0.079 & 0.095 & 0.051 & 0.043 & 0.035 \\
&  &$k=even $& 0.185 & 0.135 & 0.180 & 0.065 & 0.057 & 0.040 \\  \hline
  &$\rho=0.8$&$k=odd $ & 0.283 & 0.182 & 0.285 & 0.099 & 0.088 & 0.059  \\ 
  &&$k=even $& 0.382 & 0.234 & 0.478 & 0.122 & 0.111 & 0.067 \\   \hline
\end{tabular}
\caption{Summaries for integrated squared errors. Simulations where the algorithm of  \citet{Lin:etal:16}broke down are taken out. }
\end{table}

\begin{table}[ht]
\footnotesize
\centering
\begin{tabular}{|c|c|c|ccc|ccc|}
\hline
\multicolumn{9}{|c|}{ISE values for Model 2 (d=30)} \\
\hline
&&&\multicolumn{3}{c|}{\citet{Lin:etal:16}}&\multicolumn{3}{c|}{SBF}\\
  \hline
sample size &correlation & component & mean & median & sd & mean & median & sd \\ 
  \hline
&$\rho=0 $&$k=odd$ &0.766 & 0.724 & 0.329 & 0.431 & 0.397 & 0.200 \\ 
&& $k=even$  & 0.312 & 0.264 & 0.221 & 0.210 & 0.176 & 0.149 \\ \hline
n=200&$\rho=0.5$&$k=odd$ & 1.153 & 0.939 & 0.750 & 0.341 & 0.288 & 0.220 \\ 
&  &$k=even $& 0.526 & 0.354 & 0.554 & 0.213 & 0.175 & 0.150 \\ \hline
  &$\rho=0.8$&$k=odd $& 3.089 & 2.634 & 2.437 & 0.379 & 0.380 & 0.203 \\
  &&$k=even $& 2.176 & 1.653 & 1.875 & 0.340 & 0.312 & 0.179 \\
  \hline
    \hline
&$\rho=0 $&$k=odd$ & 0.643 & 0.637 & 0.101 & 0.329 & 0.316 & 0.101 \\
&& $k=even$  & 0.258 & 0.254 & 0.068 & 0.136 & 0.127 & 0.067 \\ 
\hline
n=500&$\rho=0.5$ &$k=odd$& 0.711 & 0.699 & 0.143 & 0.211 & 0.186 & 0.116 \\
&  &$k=even $& 0.224 & 0.214 & 0.098 & 0.119 & 0.103 & 0.075 \\ \hline
  &$\rho=0.8$&$k=odd $& 1.250 & 1.152 & 0.512 & 0.270 & 0.241 & 0.144\\ 
  &&$k=even $ & 0.253 & 0.186 & 0.241 & 0.220 & 0.205 & 0.098\\   \hline
\end{tabular}
\caption{Summaries for integrated squared errors. Simulations where the algorithm of  \citet{Lin:etal:16}broke down are taken out. }
\end{table}

\begin{figure}[h!]

%\centering
%\begin{subfigure}[t]{0.213cm}
%\centering
\includegraphics[width=13cm]{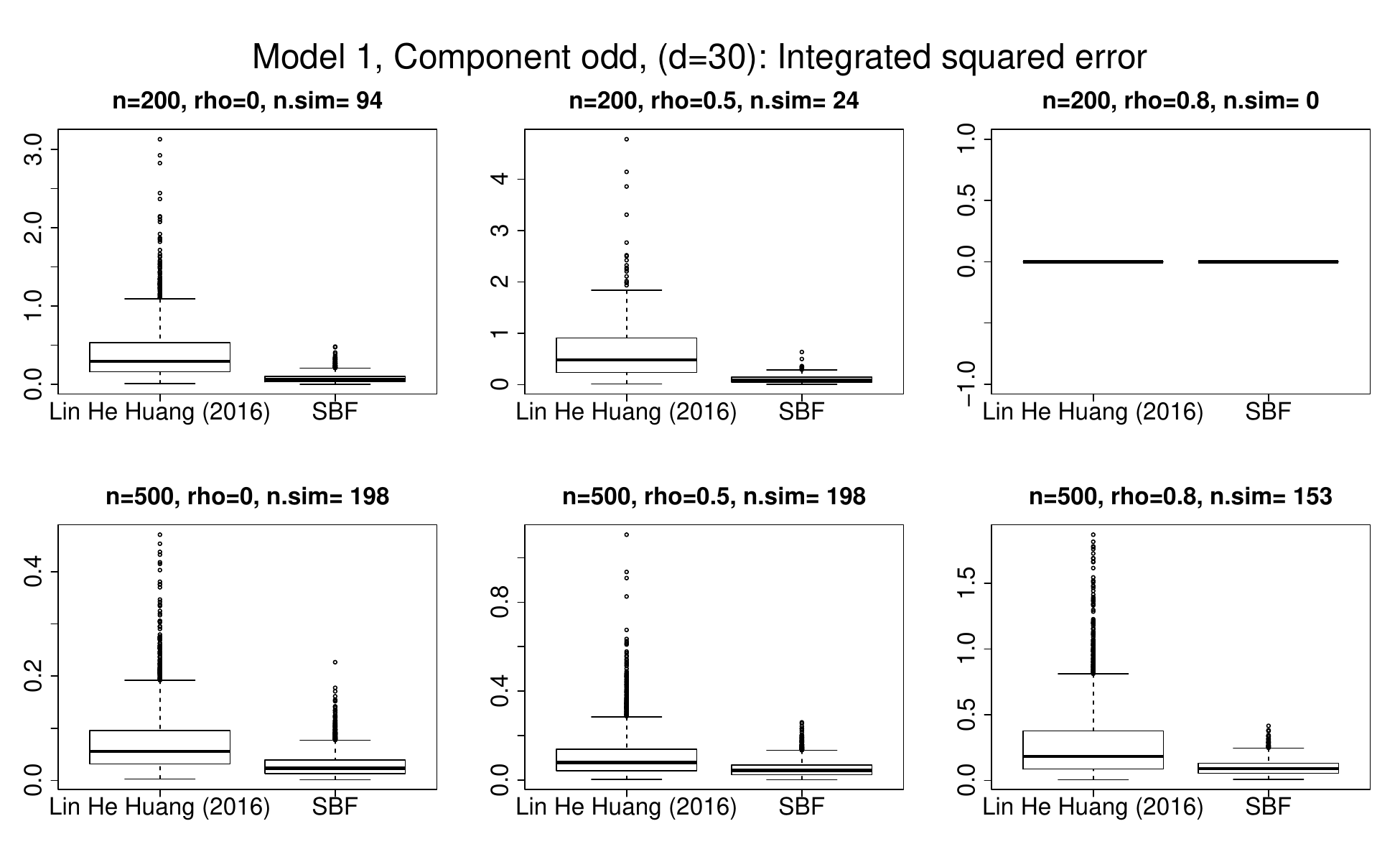} 
%\caption{Generic} \label{fig:timing1}
%\end{subfigure}

%\begin{subfigure}[t]{0.213cm}
%\centering
\includegraphics[width=13cm]{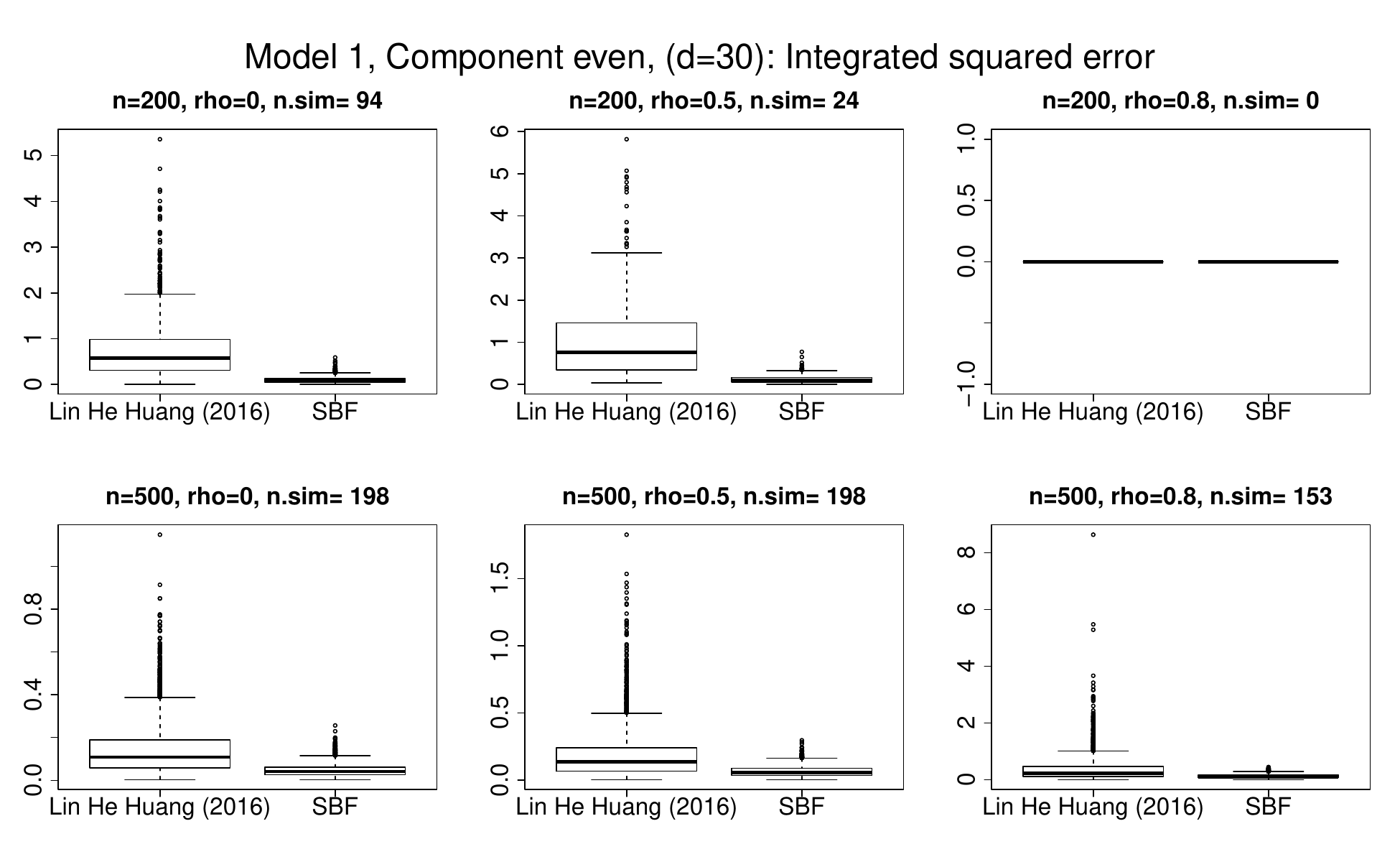} 
%\caption{Competitors} \label{fig:timing2}
%\end{subfigure}
\caption{Boxplots of the integrated squared errors in Model 1. Simulations where the algorithm of  \citet{Lin:etal:16}broke down are taken out. The value $n.sim$ is the number of successful simulations, i.e.,  200 minus number of break downs.}
\label{fig:boxplot:ise:d301}
\end{figure}

\begin{figure}[h!]
\includegraphics[width=13cm]{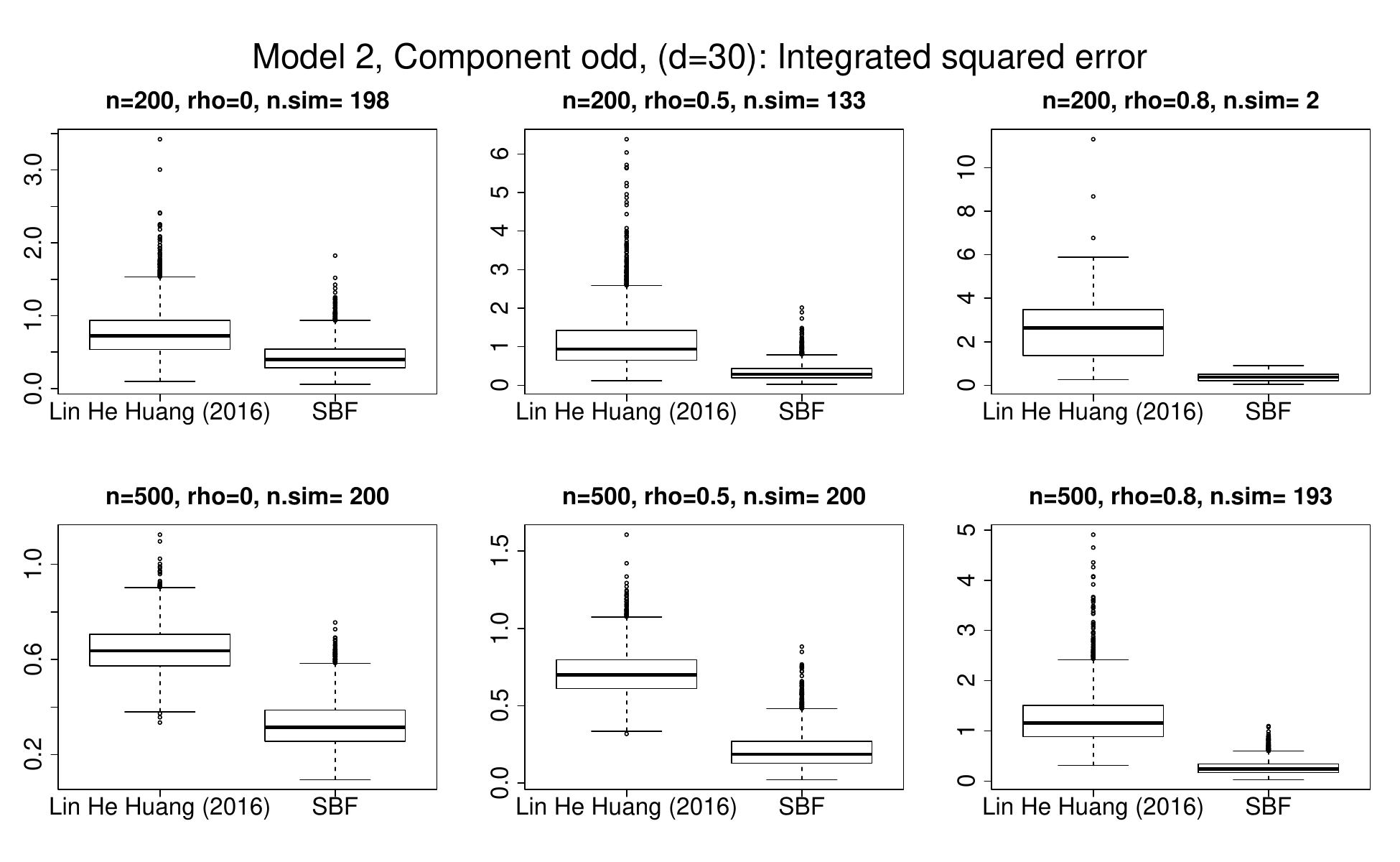} 
%\caption{Generic} \label{fig:timing1}
%\end{subfigure}

%\begin{subfigure}[t]{0.213cm}
%\centering
\includegraphics[width=13cm]{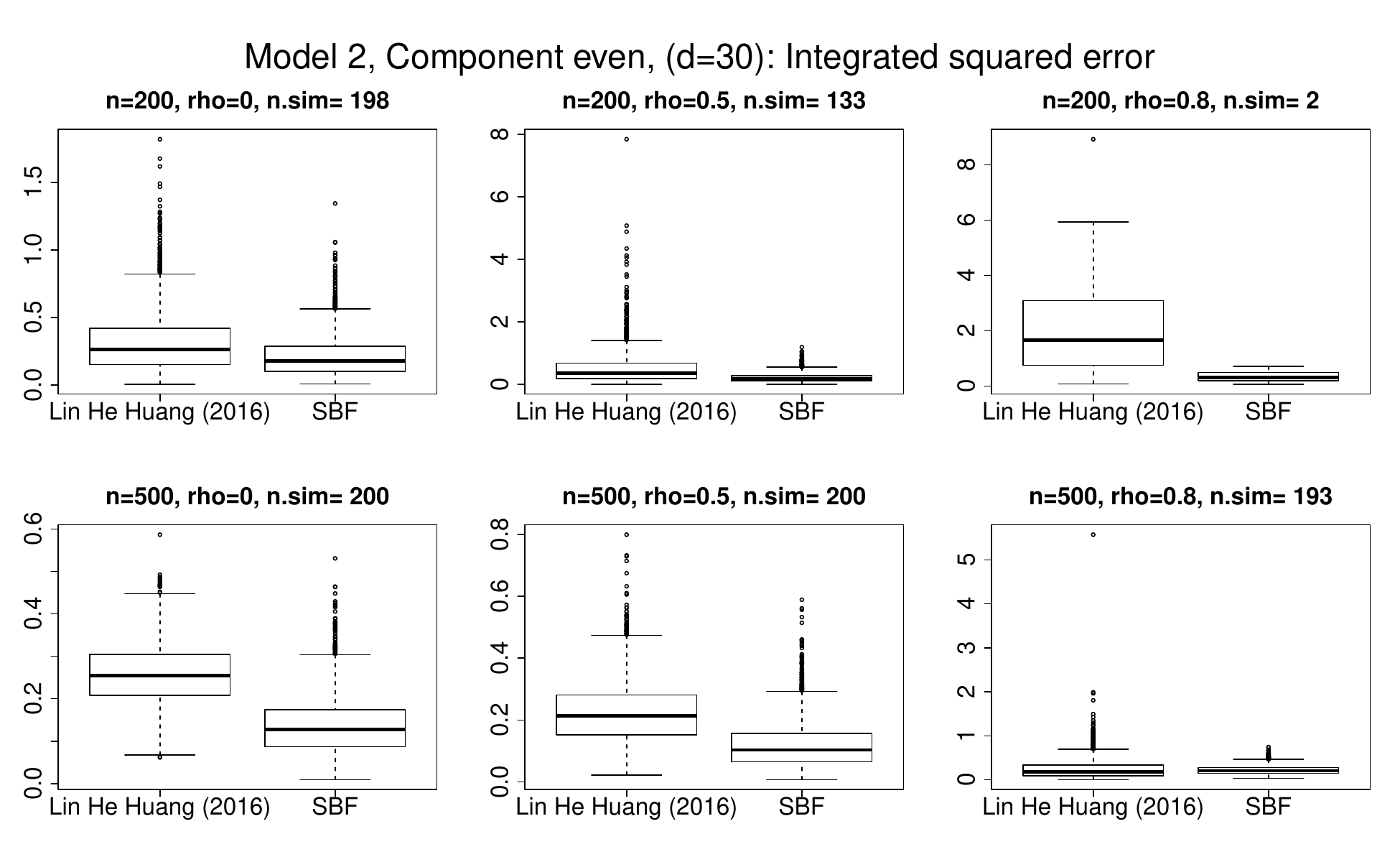} 
\caption{Boxplots of the integrated squared errors in Model 2. Simulations where the algorithm of  \citet{Lin:etal:16}broke down are taken out. The value $n.sim$ is the number of successful simulations, i.e.,  200 minus number of break downs.}
\label{fig:boxplot:ise:d302}
\end{figure}
%
%
%%\begin{subfigure}[t]{0.213cm}
%%\centering
%\begin{figure}
%\includegraphics[width=13cm]{plot30205.pdf} 
%\caption{Simualtion 36/200 components $k=29, 30$ is one example where the \citet{Lin:etal:16} estimator gives wrong impression of the underlying function}
%\label{fig:3005}
%\end{figure}

\clearpage
\subsubsection{Dimension $d=99$}

\begin{table}[h!]
\label{tab:breakdownd:d99}
\centering
\begin{tabular}{|c|cc|cc|}
\hline
\multicolumn{5}{|c|}{Number of breakdowns in \citet{Lin:etal:16}for $d=99$} \\ \multicolumn{5}{|c|} {(out of 200 simulations)} \\
\hline
 &\multicolumn{2}{c|}{Model 1}&\multicolumn{2}{c|}{Model 2}\\
  \hline
&$\rho=0$ &$\rho=0.5$&$\rho=0$ &$\rho=0.5$\\
  \hline
n=200 & 200 &200 & 200 & 200 \\ 
  n=500 & 180 & 200 &  13 & 194 \\ 
  \hline
\end{tabular}
\caption{Number of breakdowns in the algorithm of \citet{Lin:etal:16}out of 200 simulations for dimension $d=99$.}
\end{table}

\begin{table}[ht]
\footnotesize
\centering
\begin{tabular}{|c|c|c|ccc|ccc|}
\hline
\multicolumn{9}{|c|}{ISE values for Model 1 (d=99)} \\
\hline
&&&\multicolumn{3}{c|}{\citet{Lin:etal:16}}&\multicolumn{3}{c|}{SBF}\\
  \hline
sample size &correlation & component & mean & median & sd & mean & median & sd \\ 
    \hline
n=500&$\rho=0 $&$k=odd$& 0.394 & 0.290 & 0.359 & 0.064 & 0.052 & 0.046 \\ 
& &$k=even$ & 0.744 & 0.578 & 0.525 & 0.136 & 0.122 & 0.075 \\  \hline
\end{tabular}
\caption{Summaries for integrated squared errors. Simulations where the algorithm of  \citet{Lin:etal:16}broke down are taken out. }
\end{table}

\begin{table}[ht]
\footnotesize
\centering
\begin{tabular}{|c|c|c|ccc|ccc|}
\hline
\multicolumn{9}{|c|}{ISE values for Model 2 (d=99)} \\
\hline
&&&\multicolumn{3}{c|}{\citet{Lin:etal:16}}&\multicolumn{3}{c|}{SBF}\\
  \hline
sample size &correlation & component & mean & median & sd & mean & median & sd \\ 
    \hline
n=500&$\rho=0 $&$k=odd$ &1.361 & 1.289 & 0.476 & 0.738 & 0.710 & 0.246 \\ 
&& $k=even$ & 0.570 & 0.497 & 0.378 & 0.350 & 0.326 & 0.180 \\ 
\hline
n=500&$\rho=0.5$ &$k=odd$& 1.657 & 1.505 & 0.785 & 0.330 & 0.288 & 0.194 \\
&  &$k=even $& 0.659 & 0.485 & 0.643 & 0.208 & 0.179 & 0.134 \\  \hline
  \end{tabular}\caption{Summaries for integrated squared errors. Simulations where the algorithm of  \citet{Lin:etal:16}broke down are taken out. }
\end{table}

\begin{figure}[h!]

%\centering
%\begin{subfigure}[t]{0.213cm}
%\centering
\includegraphics[width=13cm]{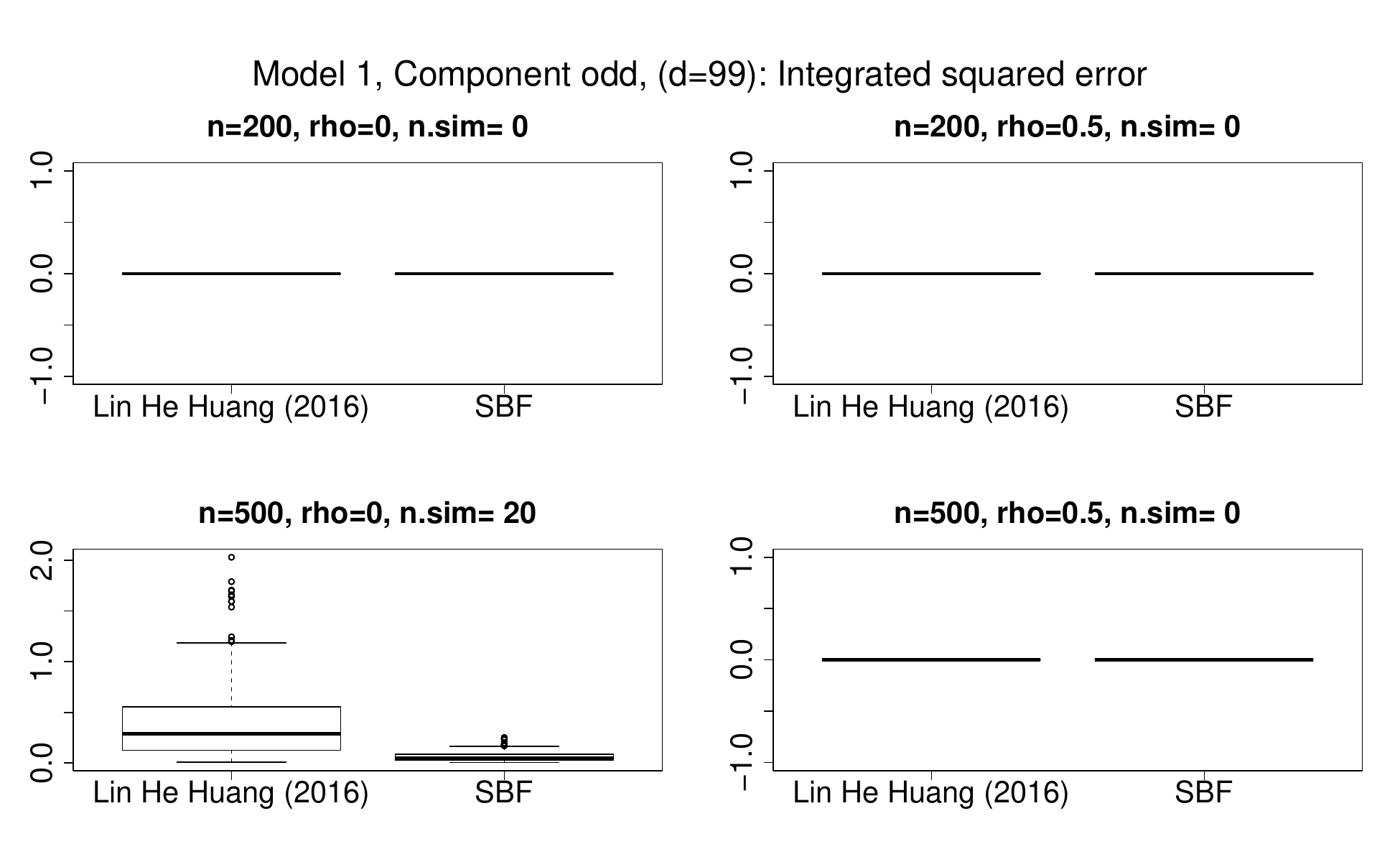} 
%\caption{Generic} \label{fig:timing1}
%\end{subfigure}

%\begin{subfigure}[t]{0.213cm}
%\centering
\includegraphics[width=13cm]{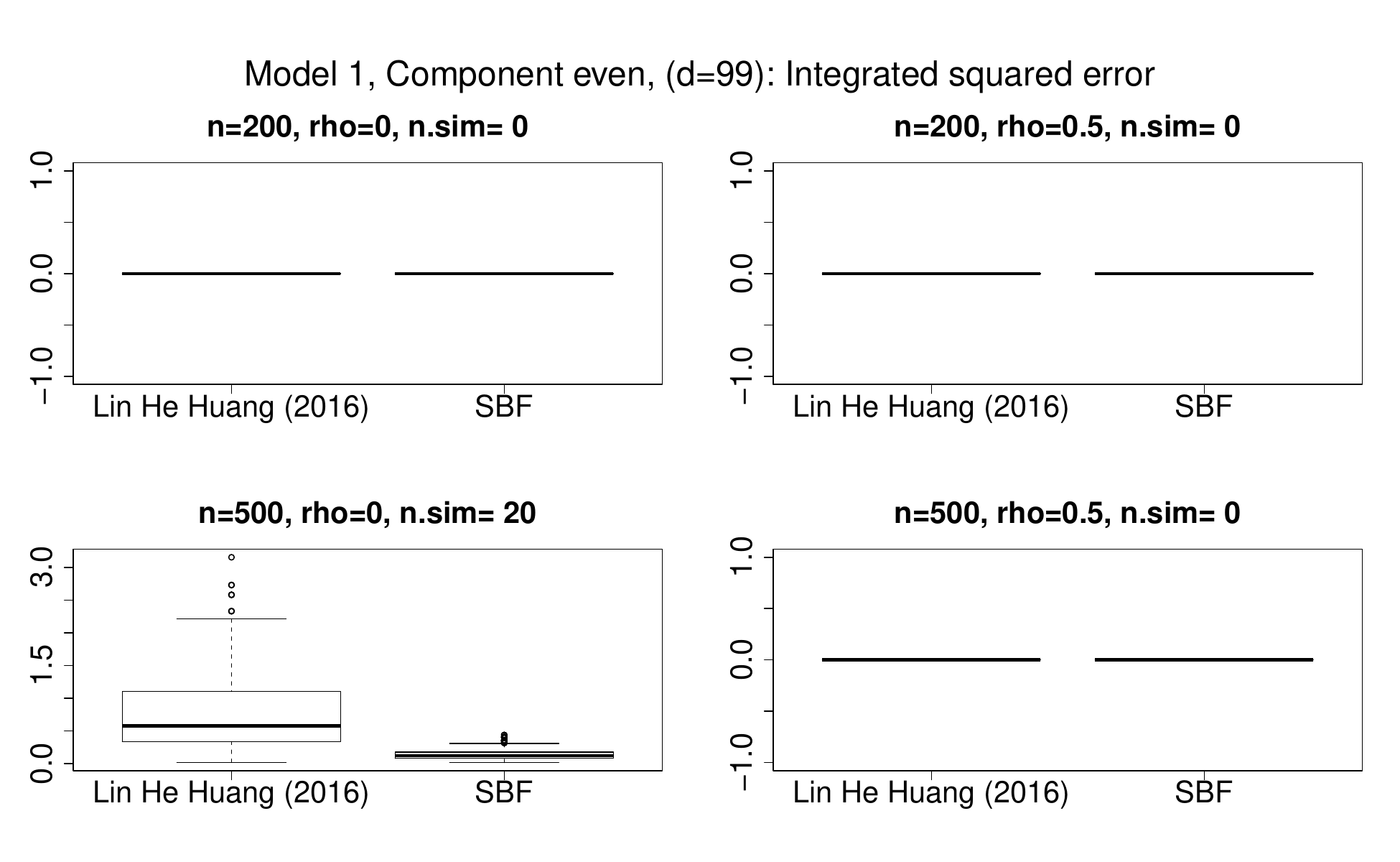} 
%\caption{Competitors} \label{fig:timing2}
%\end{subfigure}
\caption{Boxplots of the integrated squared errors in Model 1. Simulations where the algorithm of  \citet{Lin:etal:16}broke down are taken out. The value $n.sim$ is the number of successful simulations, i.e.,  200 minus number of break downs.}
\label{fig:boxplot:ise:d301}
\end{figure}

\begin{figure}[h!]
\includegraphics[width=13cm]{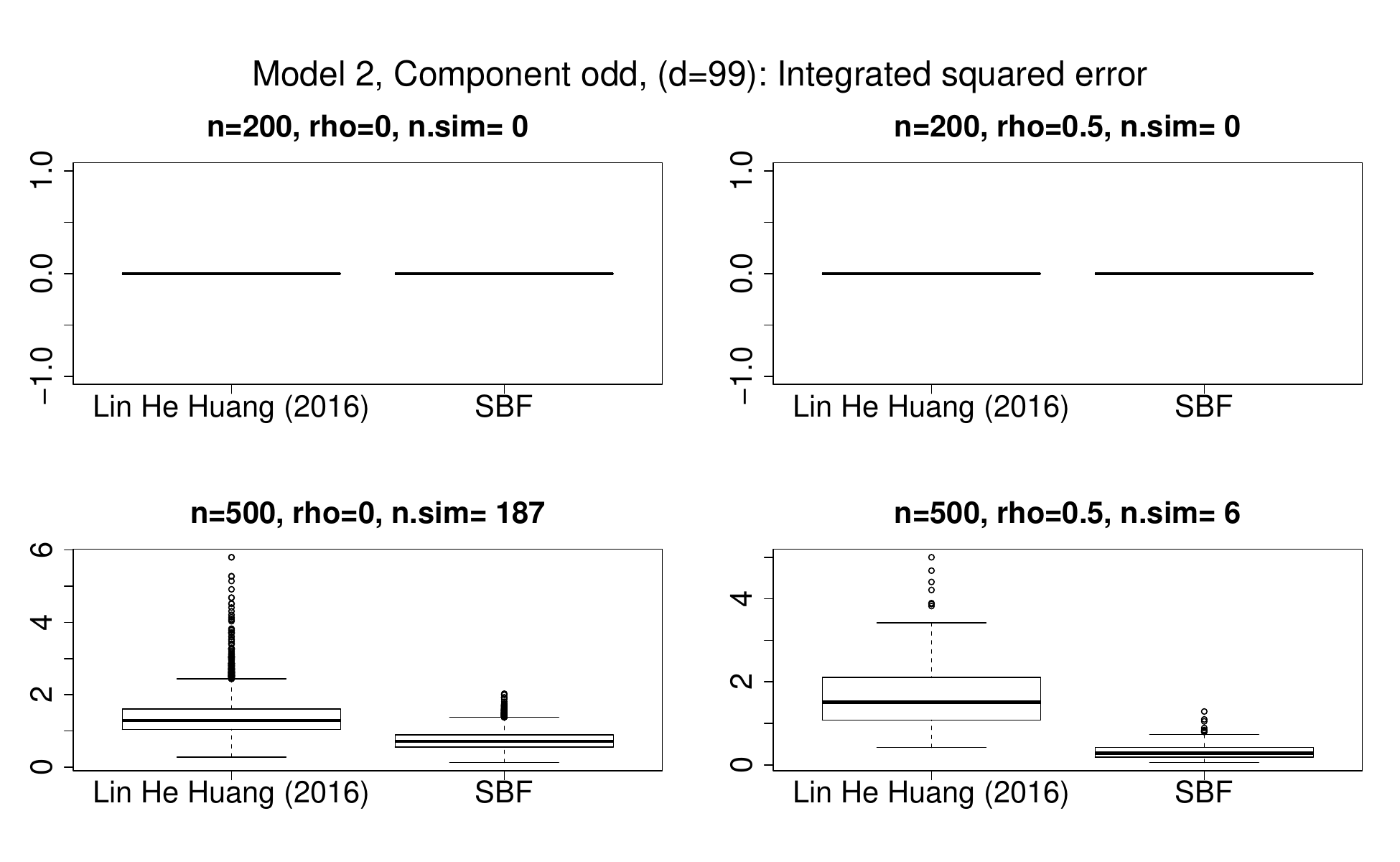} 
%\caption{Generic} \label{fig:timing1}
%\end{subfigure}

%\begin{subfigure}[t]{0.213cm}
%\centering
\includegraphics[width=13cm]{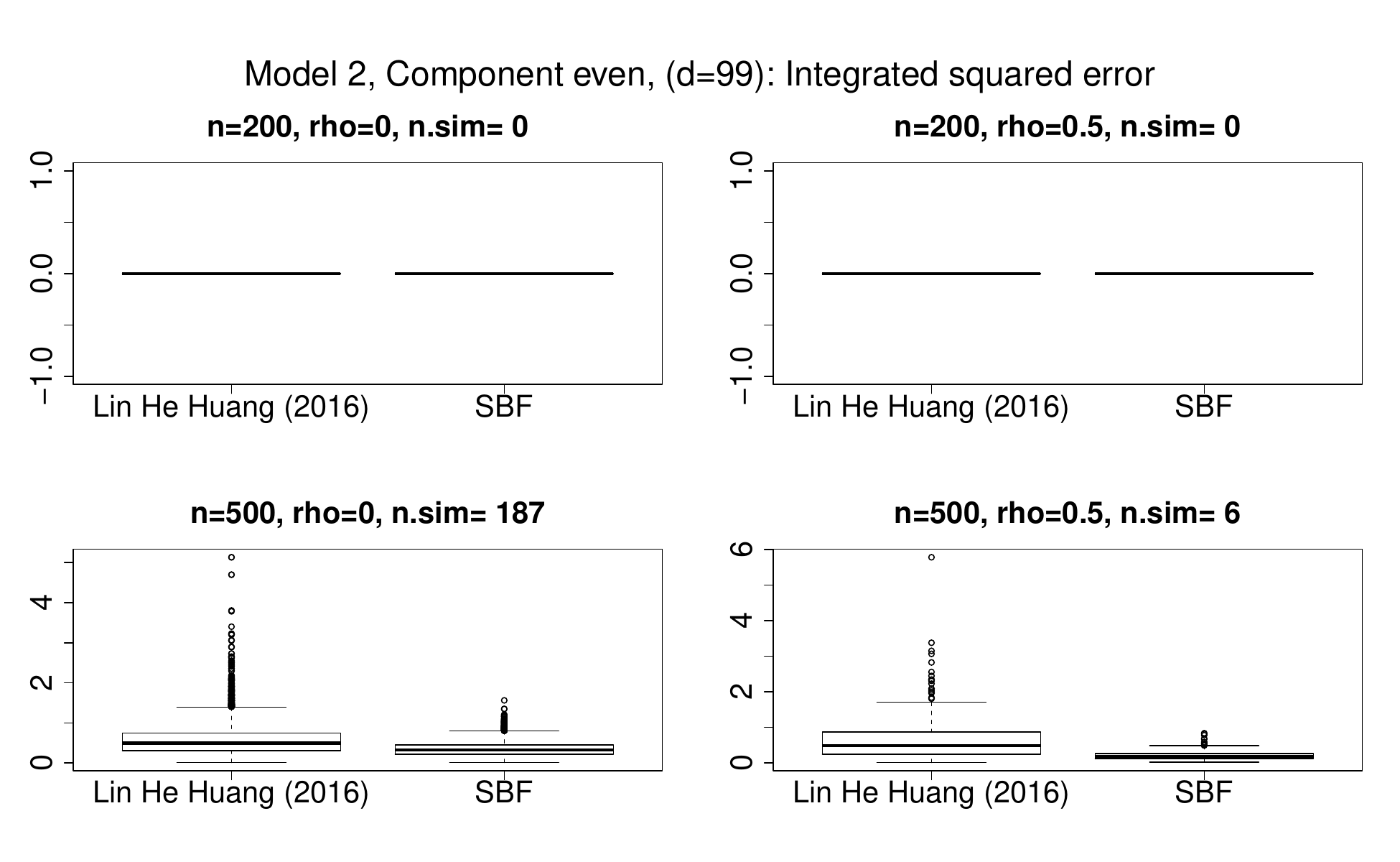} 
\caption{Boxplots of the integrated squared errors in Model 2. Simulations where the algorithm of  \citet{Lin:etal:16}broke down are taken out. The value $n.sim$ is the number of successful simulations, i.e.,  200 minus number of break downs.}
\label{fig:boxplot:ise:d302}
\end{figure}

\clearpage

\section{Bandwidth selection}\label{bandwidth}
A crucial problem in practice is finding the right amount of smoothing when using nonparametric approaches.
In the application described in this paper we have considered the maybe most straightforward way to estimate the optimal bandwidth --  the data-driven cross-validation method.

The data-driven cross-validation method in density estimation goes back to \citet{Rudemo:82} and \citet{Bowman:84}. Nowadays,  a slightly modified version (see \citet{Hall:83})  is used which aims to minimize the integrated squared error. 
 In our framework, the cross-validation bandwidth has been proposed in \citet{Nielsen:Linton:95}. 
%For the practical purposes in Section 5, in contrast to the previous sections, we have allowed the bandwidth to be different in each direction.
Cross-validation arises from the idea to minimize the integrated squared error
\[
n^{-1}\sum_{i=1}^n \int_0^{R_0} \left[\widehat \alpha\{X_i(s)\}-\alpha \{X_i(s)\}\right]^2 Y_i(s)  \mathrm ds.
\]
By expanding the square, only two of the three terms depend on the bandwidth and are thus considered.
While $\int \widehat \alpha(X_i(s))^2  \mathrm ds$ is feasible,
we have to estimate $\sum_i \int \widehat \alpha(X_i(s))\alpha(X_i(s))Y_i(s)  \mathrm ds$.
In cross-validation this is done by the unbiased  leave-one-out estimator
\[
\int \widehat \alpha^{[i]}\{X_i(s)\}\mathrm d N_i(s),
\]
where $\widehat \alpha^{[i]}$ is the leave-one-out version,  which arises from the definition of structured estimator $\widehat \alpha$ 
by setting $N_i=0$.
Finally we define the cross-validated bandwidth, $b_{CV}$, as 
\begin{align}\label{CV}
b_{CV}=\arg\min_b \sum_{i=1}^n \int \widehat \alpha(X_i(s))^2 \mathrm ds - 2\sum_{i=1}^n \int \widehat \alpha^{[i]}\{X_i(s)\}\mathrm dN_i(s).
\end{align}
Theoretical properties of cross-validation in hazard estimation in the one dimensional case are derived in \citet{Mammen:etal:15}.
To our knowledge there is no theoretical analysis of cross-validation in the multivariate hazard case of this paper.
An extensive simulation study of the multivariate case can be found in \citet{Gamiz:etal:13a}.

\clearpage

\bibliographystyle{chicago}
%\bibliography{mybib}
\bibliography{mybib}
\end{document}